\documentclass[11pt,letterpaper,twoside,reqno,nosumlimits]{amsart}

\allowdisplaybreaks

\usepackage[usenames,dvipsnames]{xcolor}
\usepackage{fancyhdr}
\usepackage{amsmath,amsfonts,amsbsy,amsgen,amscd,amssymb,amsthm}
\usepackage{url}

\usepackage{mathtools}

\usepackage[font=small,margin=0.25in,labelfont={bf},labelsep={colon}]{caption}

\usepackage{subcaption}

\usepackage{tikz}
\usepackage{microtype}
\usepackage{enumitem}

\definecolor{dark-gray}{gray}{0.3}
\definecolor{dkgray}{rgb}{.4,.4,.4}
\definecolor{dkblue}{rgb}{0,0,.5}
\definecolor{medblue}{rgb}{0,0,.75}
\definecolor{rust}{rgb}{0.5,0.1,0.1}

\usepackage[colorlinks=true]{hyperref}

\hypersetup{urlcolor=rust}
\hypersetup{citecolor=dkblue}
\hypersetup{linkcolor=dkblue}

\usepackage{setspace}

\usepackage{graphicx}
\usepackage{booktabs,longtable,tabu} 
\usepackage{multirow} 

\usepackage{float}

\usepackage[full]{textcomp}

\usepackage[scaled=.98,sups]{XCharter}
\usepackage[scaled=1.04,varqu,varl]{inconsolata}
\usepackage[type1]{cabin}
\usepackage[charter,vvarbb,scaled=1.07]{newtxmath}
\usepackage[cal=boondoxo]{mathalfa}
\usepackage[T1]{fontenc}

\usepackage{bm} 

\graphicspath{{figures/}}

\theoremstyle{definition}

\numberwithin{equation}{section} 
\numberwithin{figure}{section}
\numberwithin{table}{section}

\floatstyle{plaintop}
\newfloat{recipe}{thp}{lor}
\floatname{recipe}{Recipe}
\numberwithin{recipe}{section}

\providecommand{\mathbold}[1]{\bm{#1}}  
                                
\renewcommand{\phi}{\varphi}

\providecommand{\mathbbm}{\mathbb} 

\newcommand{\R}{\mathbbm{R}}

\newcommand{\diff}[1]{\mathrm{d}{#1}}
\newcommand{\idiff}[1]{\, \diff{#1}}

\newcommand{\vct}[1]{\mathbold{#1}}

\newcommand{\triplenorm}[1]{{\left\vert\kern-0.25ex\left\vert\kern-0.25ex\left\vert #1
    \right\vert\kern-0.25ex\right\vert\kern-0.25ex\right\vert}}

\usepackage[numbers]{natbib}

\newcommand{\om}{\omega}
\newcommand{\vom}{\vct{\omega}}
\newcommand{\vu}{\vct{u}}

\begin{document}

\title[Nearly self-similar blowup of the generalized Navier--Stokes]{{\small Nearly self-similar blowup of generalized axisymmetric Navier--Stokes equations}}

\author[T. Y. Hou]{Thomas Y. Hou\\ \\
\small Dedicated to Russel Caflisch on the occasion of his $70$th birthday}
\address{Applied and Computational Mathematics, California Institute of Technology, Pasadena, CA 91125, USA}
\email{hou@cms.caltech.edu}

\date{\today}

\begin{abstract}
We numerically investigate the nearly self-similar blowup of the generalized axisymmetric Navier--Stokes equations. First, we rigorously derive the axisymmetric Navier--Stokes equations with swirl in both odd and even dimensions, marking the first such derivation for dimensions greater than three. Building on this, we generalize the equations to arbitrary positive real-valued dimensions, preserving many known properties of the 3D axisymmetric Navier--Stokes equations. To address scaling instability, we dynamically vary the space dimension to balance advection scaling along the $r$ and $z$ directions. A major contribution of this work is the development of a novel two-scale dynamic rescaling formulation, leveraging the dimension as an additional degree of freedom. This approach enables us to demonstrate a one-scale self-similar blowup with solution-dependent viscosity. Notably, the self-similar profile satisfies the axisymmetric Navier--Stokes equations with constant viscosity. We observe that the effective dimension is approximately 3.188 and appears to converge toward 3 as background viscosity diminishes. Furthermore, we introduce a rescaled Navier--Stokes model derived by dynamically rescaling the axial velocity in 3D. 
This model retains essential properties of 3D Navier--Stokes. 
Our numerical study shows that this rescaled Navier--Stokes model with two constant viscosity coefficients exhibits a nearly self-similar blowup with maximum vorticity growth on the order of $O(10^{30})$.

\end{abstract}

\maketitle


\vspace{0.05in}

\vspace{0.05in}
Communicated by Eitan Tadmor

\section{Introduction}

The question of global well-posedness of the three-dimensional (3D) incompressible Euler and Navier--Stokes equations is one of the most fundamental and challenging open problems in fluid dynamics and nonlinear partial differential equations (PDEs) \cite{fefferman2006existence}. The primary difficulty arises from the vortex stretching mechanism, which plays a critical role in the potential formation of singularities. Despite extensive theoretical efforts and numerical investigations, determining whether smooth solutions to the 3D Navier--Stokes equations can develop finite-time singularities remains an unresolved issue.
Recent developments in the study of singularity formation for the 3D incompressible Euler equations with smooth initial data in the presence of boundaries or with H\"older continuous initial vorticity, have provided valuable insights, see e.g. \cite{elgindi2019finite,Elg19,ChenHou2023a,ChenHou2023b,ChenHou2025,elgindi2023instability,elgindi2023invertibility,chen2023remarks,cordoba2023finite,cordoba2023blow}. However, the case of smooth initial data in the interior domain continues to present significant challenges. In two recent studies by the author \cite{Hou-euler-2022,Hou-nse-2022}, a new potential blowup scenario for the axisymmetric Navier--Stokes equations was proposed. This new blowup scenario is a tornado-like traveling wave solution centered at $(R(t),Z(t))$ traveling toward the origin with an increase in maximum vorticity by a factor of $O(10^{7})$. We observed a mild scaling instability in the sense that 
$R(t)/Z(t)$ experienced a mild logarithmic growth in the late stage, which prevented us from getting arbitrarily close to the potential blowup time. 

\subsection{Derivation of the axisymmetric Navier--Stokes equations with swirl in high dimensions.}
To further explore potential blowup scenarios and address scaling instabilities, 
we introduce the generalized axisymmetric Navier--Stokes equations. One of our main contributions is a rigorous derivation of the axisymmetric Navier--Stokes in any integer dimension with $n>3$, which, to the best of our knowledge, has not appeared previously in the literature. The derivation is based on a novel representation of the velocity field in terms of cylindrical-like coordinates adapted to higher dimensions. 
An important feature of this formulation is that it allows for a natural extension to arbitrary positive real-valued dimensions while preserving many fundamental properties of the classical 3D axisymmetric Navier--Stokes equations. This generalization offers a powerful framework for studying potential singularity formation by leveraging space dimension as a dynamic parameter.

We first highlight some main ideas in deriving the axisymmetric Navier--Stokes equations in odd dimensions with $n=2m+1$ . Denote ${\bf x}=(x_1,x_2, \cdots, x_{2m},z)$. We define the three unit vectors below:
\begin{eqnarray*}
{\bf e}_r &=& r^{-1}(x_1,x_2,x_3,\cdots, x_{2m},0), \quad r =\sqrt{x_1^2+x_2^2+\cdots +x_{2m}^2} \;,\\
{\bf e}_{\theta} &=& r^{-1}(-x_2,x_1,-x_4,x_3,\cdots, -x_{2m},x_{2m-1},0),\\
{\bf e}_{z} &=& (0, 0,...,1).
\end{eqnarray*}
We call the velocity field ${\bf u}$ axisymmetric if it admits the following expression:
\begin{equation}
\label{vel-decomp}
{\bf u} = u^r(t,r,z) {\bf e}_r + u^{\theta}(t,r,z){\bf e}_{\theta} + u^z(t,r,z) {\bf e}_z .
\end{equation}
As we will see later, $u^r$ and $u^\theta$ are odd functions of $r$ as long as the solution remains smooth. Thus the above velocity field ${\bf u}$ is always well defined even at $r=0$ where $u^r {\bf e}_r =0$ and $u^\theta {\bf e}_\theta = 0$.

The choice of ${\bf e}_\theta$ plays a crucial role in our derivation of the axisymmetric Navier--Stokes equations in dimension $n > 3$. These three unit vectors 
 share the essential properties of the cylindrical coordinates:
\begin{equation}
({\bf e}_\theta \cdot \nabla ) {\bf e}_\theta = -\frac{1}{r} {\bf e}_r, \quad
({\bf e}_\theta \cdot \nabla ) {\bf e}_r = \frac{1}{r} {\bf e}_\theta, \quad
({\bf e}_r \cdot \nabla ) {\bf e}_\theta = 0, \quad
({\bf e}_r \cdot \nabla ) {\bf e}_r = 0 .
\end{equation}
As we will show in Section 2.4, it also satisfies a very important property 
\begin{equation}
\nabla \cdot (u^\theta {\bf e}_\theta) = 0 \;.
\end{equation}
As a consequence, the divergence free condition is reduced to
\begin{equation}
\nabla \cdot {\bf u} =
\frac{(r^{n-2} u^r)_r}{r^{n-2}} + \frac{(r^{n-2} u^z)_z}{r^{n-2}} = 0.
\end{equation} 
To derive the axisymmetric Navier--Stokes equations, we need to decompose the advection term:
\begin{eqnarray}
\label{conv-term-0}
({\bf u} \cdot \nabla) {\bf u} = 
\left (u^r \partial_r u^r + u^z \partial_z u^r - \frac{(u^\theta)^2}{r} \right ) {\bf e}_r 
+ \left (u^r \partial_r u^\theta + u^z \partial_z u^\theta + \frac{u^r}{r} u^\theta \right ) {\bf e}_\theta 
+ (u^r \partial_r u^z + u^z \partial_z u^z){\bf e}_z .
\end{eqnarray}
We will show that if the initial velocity field is of the form \eqref{vel-decomp}, it will remain to have the same form as long as the solution is smooth and we can derive the equivalent axisymmetric Navier--Stokes equations for $u^r$, $u^\theta$ and $u^z$ as functions of $r$ and $z$ only, see Section \ref{odd-dim-case} for the derivation. 

We have also derived the axisymmetric Navier--Stokes equations in even dimensions by separating ${\bf x}=(x_1,x_2, \cdots, x_{2m},z_1,z_2)$ into the first $2m$ dimensions plus the last two dimensions. The first $2m$ dimensions can be parametrized in the same way as in the odd-dimensional case and the last two dimensions are parametrized using the generalized polar coordinates, see Section \ref{even-dim-case}
 for the derivation and Section \ref{Illustration} for an illustration of the flow geometry in both four and five space dimensions. 

Denote by $u^\theta$, $\om^\theta$ and $\psi^\theta$ as the angular velocity, angular vorticity and angular stream function, respectively. The axisymmetric Navier--Stokes equations in odd dimensions are given below:
\begin{subequations}\label{eq:axisymmetric_NSE_0n}
\begin{align}
&\Gamma_{t}+u^r\Gamma_{r}+u^z \Gamma_{z} = \nu \left(\Gamma_{rr} + \frac{(n-4)}{r}\Gamma_{r} + \frac{(6-2n)}{r^2}\Gamma + \Gamma_{zz} \right),\label{eq:as_NSE_0_a}\\
&\om_{1,t}+u^r\om_{1,r}+u^z\om_{1,z}  = \left ( \frac{\Gamma^2}{r^4} \right )_z - (n-3) \psi_{1,z} \om_1 + \nu\left(\om_{1,rr} + \frac{n}{r}\om_{1,r} + \om_{1,zz} \right),\label{eq:as_NSE_0_b}\\
& -\left(\partial_r^2+\frac{n}{r}\partial_r+\partial_z^2\right)\psi_1 = \om_1,\label{eq:as_NSE_0_c}
\end{align}
\end{subequations}
where $\om_1 = \om^\theta/r$, $\psi_1=\psi^\theta/r$, $u^r = -(r^{n-2}\psi^\theta)_z/r^{n-2}, \; u^z  =(r^{n-2}\psi^\theta)_r/r^{n-2}$, and $\Gamma =r u^\theta$ is the total circulation. When $n=3$, we recover the $3$D axisymmetric Navier--Stokes equations. 


Building on this, we generalize the above axisymmetric Navier-Stokes equations \eqref{eq:axisymmetric_NSE_0n} to any dimension by taking $n$ to be any positive value. 
This generalized version of the axisymmetric Navier--Stokes equations enjoys many known properties of the $3$D axisymmetric Navier--Stokes equations,
 including the conservation of the total circulation, the incompressibility condition $(r^{n-2} u^r)_r + (r^{n-2}u^z)_z = 0$, and the conservation of the kinetic energy $ \int | {\bf u}|^2 r^{n-2} dr dz $ for smooth solutions and a given constant dimension $n$. Moreover, many of the known non-blowup criteria also apply to the generalized axisymmetric Navier-Stokes equations for a given constant dimension $n$.

\subsection{Scaling instability and a novel two-scale dynamic rescaling formulation}

The scaling instability is induced by the two scaling factors $R(t)$ and $Z(t)$ for a traveling wave singularity centered at $(R(t),Z(t))$. Since we impose an odd symmetry of $\Gamma$, $\omega_1$ and $\psi_1$ as a function of $z$, $Z(t)$ naturally characterizes the small scale of the solution. On the other hand, $R(t)$ is another scale measuring the distance between the singularity location and the symmetry axis $r=0$. These two scales can vary independently as $(R(t),Z(t)) \rightarrow (0,0)$.
If we use the rescaled variables $\xi=r/R(t)$ and $\eta = z/Z(t)$, 
we can see that $u^r = -r \psi_{1,z}$ is of order $R(t)/Z(t)$ larger than $u^z = (n-1)\psi_1 + r \psi_{1,r}$. This scaling imbalance seems to be responsible for the scaling instability for the traveling wave singularity.

A significant contribution of this work is the introduction of a novel two-scale dynamic rescaling formulation by rescaling the $r$ and $z$ directions independently.
This additional scaling parameter allows us to dynamically rescale $(R(t),Z(t))$ to a fixed location $(R_0,1)$. 
To balance the advection scaling along the $r$ and $z$ directions, we treat the dimension as an additional degree of freedom and dynamically vary the dimension by setting $n = 1+2 R(t)/Z(t)$. 
This seems to be very effective in eliminating the scaling instability. We observe a stable one-scale self-similar blowup solution with solution-dependent viscosity. Moreover, $R(t)/Z(t)$ seems to converge to a fixed constant as $t \rightarrow T$.

\vspace{-0.1in}
\subsection{Self-similar blowup of the generalized Navier-Stokes equations.} We solve the generalized axisymmetric Navier--Stokes equations with solution dependent viscosity of the form 
$\nu (t) = \nu_0 \|u_1\|_\infty Z(t)^2$, which is scaling invariant. Here $u_1 =u^\theta/r$. With $\nu_0=0.006$, we observe that the generalized axisymmetric Navier--Stokes equations with viscosity $\nu(t)$ develop a stable self-similar blowup with dimension $n(t) \rightarrow 3.188$ as $t \rightarrow T$. More specifically, the self-similar blowup solution is of the form:
\begin{eqnarray*}
{\psi}_{1} (t, r,z) &=& \frac{\lambda(t)}{(T-t)^{1/2}}  {\Psi} \left ( t, \frac{r}{\lambda(t)\sqrt{(T-t)}}, \frac{z}{\lambda(t)\sqrt{(T-t)}} \right ), \\
{u}_1 (t, r,z) &=& \frac{1}{(T-t)}  {V} \left ( t, \frac{r}{\lambda(t)\sqrt{(T-t)}}, \frac{z}{\lambda(t)\sqrt{(T-t)}} \right ), \\
{\omega}_1 (t, r,z) &=& \frac{1}{\lambda(t) (T-t)^{3/2}} {\Omega} \left ( t, \frac{r}{\lambda(t)\sqrt{(T-t)}}, \frac{z}{\lambda(t)\sqrt{(T-t)}} \right ), 
\end{eqnarray*}
where $u_1 =u^\theta/r$ and $\lambda (t) \approx (T-t)^{0.0233}$ and we normalize $V$ to satisfy 
$V(t,R(t),Z(t))=1$. The self-similar profile $({\Psi}, V, \Omega)$ seems to converge to a steady state as $t \rightarrow T$. If we denote $\lambda(t)\sqrt{(T-t)} \equiv (T-t)^{c_l}$, then we have $c_l \approx 0.5233$ and $\nu= \nu_0 \|u_1\|_\infty Z(t)^2 = \nu_0 (T-t)^{2c_l-1}$ since $Z(t) = (T-t)^{c_l}$. We observe that
the maximum vorticity scales like  $\|u_1(t)\|_\infty$. This implies that $\| {\vom}(t) \|_\infty$ scales like $O((T-t)^{-1})$, which is consistent with the Beale-Kato-Majda non-blowup criterion \cite{beale1984remarks}.

We say that the solution develops a nearly self-similar blowup if the rescaled profile $({\Psi}, V, \Omega)$ and the scaling parameters $c_l$ are sufficiently close to an approximate steady state as the rescaled time $\tau \rightarrow \infty$. We indeed observe that the rescaled profile seems to converge to a steady state as $\tau \rightarrow \infty$.

One surprising consequence of this solution dependent viscosity is that if the solution of the dynamic rescaling equations converges to a steady state, the self-similar profile satisfies the self-similar equations of the axisymmetric Navier--Stokes equations with constant viscosity $\nu_0$ after rescaling the $\xi$ variable:  
\begin{eqnarray}
\label{Dyn-rescale-eqn3}
&& (c_{l} \xi, c_{l}\eta) \cdot\nabla_{(\xi,\eta)} \widetilde{\Gamma}+\widetilde{\bf u} \cdot \nabla_{(\xi,\eta)} \widetilde{\Gamma} =  \nu_0  \widetilde{\Delta} \widetilde{\Gamma} ,\\
&&(c_{l} \xi, c_{l}\eta) \cdot\nabla_{(\xi,\eta)} \widetilde{\omega}_{1}+\widetilde{\bf u} \cdot \nabla_{(\xi,\eta)} \widetilde{\omega}_1 = c_\omega \widetilde{\omega}_1 + (\widetilde{u}_1^2)_\eta - (n-3)\widetilde{\psi}_{1,\eta} \widetilde{\omega}_1 +\nu_0 \Delta \widetilde{\omega}_1 ,\\
&& - \Delta \widetilde{\psi}_1 = \widetilde{\omega}_1 ,  \quad
\Delta = -\left(\partial_{\xi}^2+\frac{n}{\xi}\partial_{\xi} + \partial_{\eta}^2\right),
\end{eqnarray}
where $c_l$ and $c_\omega$ are scaling parameters,
$\widetilde{\Delta}$ is defined in the same way as the diffusion operator defined on the right hand side of \eqref{eq:as_NSE_0_a} except that we replace $r$ by $\xi$ and $z$ by $\eta$.

Our numerical experiments show that the effective space dimension approaches approximately 3.188 and converges toward 3 as the background viscosity $\nu_0$ diminishes. This observation suggests that the generalized Navier--Stokes equations can serve as a useful model for investigating the original 3D Navier--Stokes equations by gradually reducing the background viscosity $\nu_0$.


We observe that the solution dependent viscosity $\nu(t) \rightarrow 0$ as $t \rightarrow T$. However, the solution dependent viscosity still has an $O(1)$ effect on the self-similar blowup profile since the limiting rescaled profile does not satisfy the corresponding Euler equations as the vanishing viscosity limit. Instead, it satisfies the Navier-Stokes equations with constant viscosity $\nu_0$, which is quite surprising.


\vspace{-0.1in}
\subsection{Self-similar blowup of the rescaled Navier--Stokes model with constant viscosity}

We also investigate a rescaled Navier--Stokes model. Recall the 3D Navier-Stokes equations with $m=n=3$:  
\begin{subequations}\label{eq:axisymmetric_NSE_01}
\begin{align}
&\Gamma_{t}+u^r\Gamma_{r}+u^z \Gamma_{z} = \nu_1 \left(\Gamma_{rr} + \frac{({n}-4)}{r}\Gamma_{r} + \frac{(6-2{n})}{r^2}\Gamma + \Gamma_{zz} \right),\label{eq:as_NSE_01_a}\\
&\om_{1,t}+u^r\om_{1,r}+u^z\om_{1,z}  = \left ( \frac{\Gamma^2}{r^4} \right )_z + \nu_2\left(\om_{1,rr} + \frac{{n}}{r}\om_{1,r} + \om_{1,zz} \right),\label{eq:as_NSE_01_b}\\
& -\left(\partial_r^2+\frac{{n}}{r}\partial_r+\partial_z^2\right)\psi_1 = \om_1,\label{eq:as_NSE_01_c}
\end{align}
\end{subequations}
where $u^r= -r\psi_{1,z}, \; u^z =(m-1)\psi_1 + r \psi_{1,r}$. To balance the advection scaling along the $r$ and $z$ directions, we can either rescale $u^z$ by $R/Z$ or rescale $u^r$ by $Z/R$.
Since we observe that $Z/R < 1$ in our computation, in order not to penalize the stablizing effect of advection \cite{hou2008dynamic,lei2009stabilizing},
we rescale $u^z$ by $R/Z$ by choosing $m=1+2R/Z$ to balance the advection scaling along the $r$ and $z$ directions. However, rescaling $u^z$ by $R/Z$ with the choice $m = 1 + 2R/Z$ leads to an $m$-dimensional axisymmetric velocity field $(u^r, u^z)$, which is no longer three-dimensional. Since we omit the term $-(m - 3)\psi_{1,z} \omega_1$ on the right-hand side of \eqref{eq:as_NSE_0_b}, the resulting system does not correspond to a true $m$-dimensional Navier--Stokes equation, even when $m$ is an odd integer. For this reason, we refer to it as the {\it rescaled Navier--Stokes model}.

On the other hand, by choosing ${n} = 2m - 3$, we can show that the rescaled Navier--Stokes model conserves a generalized kinetic energy, provided that $R/Z$ remains constant. Although this model is not identical to the original 3D Navier--Stokes equations, it retains many of their key structural properties for constant $R/Z$, including the incompressibility condition $(r^{m-2} u^r)_r + (r^{m-2}u^z)_z = 0$, the conservation of total circulation $\Gamma$. Moreover, when $R=Z$ and $\nu_1=\nu_2$,  we can recover the original $3$D axisymmetric Navier--Stokes equations. More importantly, this rescaled Navier--Stokes model exhibits finite-time, nearly self-similar blowup. Thus, studying this rescaled Navier--Stokes model may offer valuable insights into the potential blowup mechanisms of the original $3$D Navier--Stokes equations.

In our study, we use a small viscosity coefficient ($\nu_1=0.0006$) for $\Gamma$ and a larger viscosity coefficient $\nu_2=10\nu_1$ for $\omega_1$.
We obtain the following scaling properties for the potential blowup solution:
\begin{eqnarray*} 
{\psi}_{1} (t, r,z) &=& \frac{\lambda(t)}{(T-t)^{1/2}}  {\Psi} \left ( t, \frac{r}{\lambda(t)\sqrt{(T-t)}}, \frac{z}{\lambda(t)\sqrt{(T-t)}} \right ), \\
{u}_1 (t, r,z) &=& \frac{1}{(T-t)}  {V} \left ( t, \frac{r}{\lambda(t)\sqrt{(T-t)}}, \frac{z}{\lambda(t)\sqrt{(T-t)}} \right ), \\
{\omega}_1 (t, r,z) &=& \frac{1}{\lambda(t)(T-t)^{3/2}} {\Omega} \left ( t, \frac{r}{\lambda(t)\sqrt{(T-t)}}, \frac{z}{\lambda(t)\sqrt{(T-t)}} \right ), 
\end{eqnarray*}
where
$\lambda (t) = (1+\epsilon |\log(T-t)|)^{-1/2}$ for a small $\epsilon$. 
We normalize $V$ to satisfy
$V (t,R(t),Z(t))=1$. Due to a smaller $\nu_1$, the total circulation
$\Gamma$ generates a sharp shock like traveling wave propagating toward the origin. 
This sharp shock front produces a regularized Delta function like profile for $\omega_1$ through the source term $(\Gamma^2/r^4)_z$ and the advection terms. 
The relatively large viscosity coefficient $\nu_2$ for $\omega_1$ stabilizes the shearing instability induced by the sharp front of $\Gamma$ and generates a stable nearly self-similar tornado like traveling wave singularity
at the origin with maximum vorticity increased by $O( 10^{30})$ and $n \approx 4.73$ by $\tau=155$ ($\tau$ is the rescaled time). 
We observe that the rescaled profile seems to converge to an approximate steady state after proper rescaling as the rescaled time $\tau \rightarrow \infty$.
Since the maximum vorticity $\| \vom(t) \|_\infty$ has the same scaling as $\|u_1\|_\infty$, 
we conclude that $\| \vom(t) \|_\infty$ scales like $O((T-t)^{-1})$, which is consistent with the Beale-Kato-Majda non-blowup criterion \cite{beale1984remarks}. 

\vspace{-0.1in}

%

\subsection{Confirming the blowup solution using various blowup criteria}



We have examined the generalized Ladyzhenskaya-Prodi-Serrin regularity criteria \cite{ladyzhenskaya1957,prodi1959,serrin1962} that are based on the estimate of the $L_t^q L_x^p$ norm of the velocity with ${n}/p + 2/q \leq 1$, which is formulated in $\R^n$. We study the cases of $(p,q) = (4{n}/3,8),\; (2{n},4),\; (3{n},3)$, and $(\infty,2)$ respectively. Denote by $\|{\bf u}(t)\|_{L^{p,q}} = \left (\int_0^t \|{\bf u(s)}\|_{L^p(\Omega)}^q ds \right )^{1/q}$.  Our numerical results show that $\|{\bf u}(t)\|_{L^{p,q}}^q$ blows up roughly with a logarithmic rate, $O(|\log(T-t)|)$ for $p$ large, e.g. $p=2{n}, 3{n}, \infty$. This provides strong evidence for the development of a potential finite time singularity of the rescaled Navier--Stokes model. 

We have further investigated the nonblowup criteria based on the $L^3$ norm of the $3$D velocity due to Escauriaza-Seregin-Sverak \cite{sverak2003}. Since the energy inequality is in the ${n}$-dimensional setting, we should consider the $L^{{n}}$ norm of the velocity \cite{palasek22}, which is scaling invariant.
We observe a mild logarithmic growth of 
$\| {\bf u}(t)\|_{L^{{n}}}$ for the rescaled Navier--Stokes model. Moreover, we examine the growth of the negative pressure and observe that $\|-p \|_\infty$ and 
$\| 0.5 | \nabla {\bf u}| + p \|_\infty$ blow up like $O(1/(T-t)$.
This provides further evidence for the potential finite time singularity of the rescaled Navier--Stokes model.

For the axisymmetric $3$D Navier--Stokes equations, there are two more non-blowup criteria. In the work by Yau et al \cite{chen2008lower,chen2009lower} and Sevrak et al \cite{sverak2009}, they exclude finite time blowup if the velocity field satisfies $\| {\bf u} \|_\infty \leq \frac{C}{\sqrt{T-t}}$ provided that $ \| r u^r \|_\infty $ and $\| r u^z \|_\infty $ remain bounded for $r \geq r_0 > 0$. Our numerical study shows that $ \| r u^r \|_\infty $ has a mild logarithmic growth in time. In the work by Wei \cite{wei2016criticality} (see also \cite{lei2017criticality}), finite time blowup of the $3$D axisymmetric Navier--Stokes is excluded if the condition $|\log (r)|^{3/2} | \Gamma (t, r,z)| \leq 1$ for $ r \leq \delta_0 < 1/2$. If we assume that their key estimate based on the Hardy inequality still holds, their result should also apply to the rescaled Navier--Stokes model. Our numerical result shows that $\max_{r \leq \delta_0}|\log (r)|^{3/2} | \Gamma (t, r,z)| $ can grow roughly like $O(|\log(T-t)|)$ as $t \rightarrow T$. This provides additional evidence for the finite time blowup of the rescaled Navier--Stokes model.

\vspace{-0.1in}
\subsection{Review of previous works}

For the $3$D Navier--Stokes equations, the partial regularity result due to Caffarelli--Kohn--Nirenberg \cite{caffarelli1982partial} is one of the best known results (see a simplified proof by Lin \cite{lin1998new}). This result implies that any potential singularity of the axisymmetric Navier--Stokes equations must occur on the symmetry axis. There have been some very interesting theoretical developments regarding the lower bound on the blow-up rate for axisymmetric Navier--Stokes equations \cite{chen2008lower,chen2009lower,sverak2009}. Another interesting development is a result due to Tao \cite{tao2016nse} who proposed an averaged three-dimensional Navier--Stokes equation that preserves the energy identity, but blows up in finite time. 

There have been a number of theoretical developments for the $3$D incompressible Euler equations, including the Beale--Kato--Majda blow-up criterion \cite{majda2002vorticity}, the geometric non-blow-up criterion due to Constantin--Fefferman--Majda \cite{cfm1996} and its Lagrangian analog due to Deng-Hou-Yu \cite{dhy2005}. Inspired by their work on the vortex sheet singularity \cite{caflisch89}, Caflisch and Siegel have studied complex singularity for $3$D Euler equation, see \cite{caflisch93,caflisch09}, and also \cite{frisch06} for the complex singularities for $2$D Euler equation.

In 2021, Elgindi \cite{elgindi2019finite} (see also \cite{Elg19}) proved an exciting result: the $3$D axisymmetric Euler equations develop a finite time singularity for a class of $C^{1,\alpha}$ initial velocity with no swirl and a very small $\alpha$. There have been a number of interesting theoretical results inspired by the Hou--Lou blowup scenario \cite{luo2014potentially,luo2014toward}, see e.g. \cite{kiselev2014small,choi2014on,choi2015,kryz2016,chen2019finite2,chen2020finite,chen2021finite3}
and the excellent survey article \cite{kiselev2018}. We remark that Huang-Qin-Wang-Wei recently proved the existence of exact self-similar blowup profiles for the gCLM model and the Hou-Luo model by using a purely analytic fixed point method in \cite{Huang2023a,Huang2023b}.

There has been substantial progress on
singularity formation of $3$D Euler equations in recent years. In \cite{ChenHou2023a,ChenHou2023b,ChenHou2025}, Chen and Hou have established a computer-assisted proof of finite time blowup for the $2$D Boussinesq and the $3$D axisymmetric Euler equations with smooth initial data and boundary by proving the nonlinear stability of an approximate self-similar profile.
In \cite{elgindi2023instability,elgindi2023invertibility}, Elgindi-Pasqualotto established blowup of $2$D Boussinesq and $3$D Euler equations (with large swirl) with $C^{1,\alpha}$ velocity and without boundary. In \cite{cordoba2023finite}, Cordoba-Martinez-Zoroa-Zheng developed a new method different from the above self-similar approach to establish blowup of axisymmetric Euler equations with no swirl and with ${\bf u} (t) \in C^{\infty}(R^3 \backslash O) \cap C^{1,\alpha} \cap L^2 $. In \cite{chen2023remarks}, Chen proved that such a blowup result can also be established by the self-similar approach. By adding an external force $f$ uniformly bounded in $C^{1,1/2-}$ up to the blowup time, the authors of \cite{cordoba2023blow} established blowup of $3$D Euler with smooth velocity. 

We remark that in two recent papers \cite{Tsai2022,Tsai2023} 
the authors studied  axisymmetric, swirl-free solutions of the Euler equations in dimensions $d \geq 4$. In \cite{Tsai2022}, Miller and Tsai investigated the potential finite time blowup in higher dimensions in $d \geq 4$ and proved a sufficient condition for finite time blowup. In \cite{Tsai2023}, Gustafson-Miller-Tsai established lower bounds on the growth of generalized anti-parallel vortex tube pair solutions in three and higher dimensions.

There have been relatively few papers on the numerical study regarding the potential blow-up of the $3$D Navier--Stokes equations. We refer to a recent survey paper \cite{protas2022} by Protas on a systematic search for potential singularities of the Navier--Stokes equations by solving PDE optimization problems. It concludes that ``No evidence
for singularity formation was found in extreme
Navier--Stokes flows constructed in this manner in three dimensions.'' There were a number of attempts to look for potential Euler singularities numerically, see \cite{gs1991,es1994,bp1994,kerr1993,hl2006,hl2008,luo2014potentially,luo2014toward,brenner2016euler,Hou-Huang-2022,Hou-Huang-2023,wang23}. 
We refer to a review article \cite{gibbon2008} for more discussions on potential Euler singularities.

Throughout this paper, we investigate finite-time blowup of two models for the Navier--Stokes equations. The first one is called the generalized Navier--Stokes equation defined in \eqref{eq:axisymmetric_NSE_0n} with $n=1+2R(t)/Z(t)$ and solution-dependent viscosity $\nu(t) = \nu_0 \|u_1\|_\infty Z(t)^2$, which is studied in Section~4 (see also the dynamic rescaling formulation \eqref{Dyn-rescale-eqn-Gamma-new}). The second one is called the rescaled Navier--Stokes model with two constant viscosity coefficients $(\nu_1,\nu_2)$ and $m=1+2R(t)/Z(t)$ ($n=2m-3$), which is defined in \eqref{eq:axisymmetric_NSE_01} and investigated in Section~5 (see also the dynamic rescaling formulation \eqref{Dyn-rescale-eqn-rescaled-NS}).

The rest of the paper is organized as follows. In Section \ref{sec:setup1}, we provide a rigorous derivation of the axisymmetric Navier--Stokes equations in both odd and even dimensions. We obtain the generalized Navier--Stokes equations by treating the dimension as a free parameter and perform the energy estimates for the generalized Navier--Stokes equations and the rescaled Navier--Stokes model. In Section \ref{sec:setup2}, we introduce our two-scale dynamic rescaling formulation. In Section \ref{sec:euler},  we investigate the self-similar blowup of the generalized Navier--Stokes equations with solution dependent viscosity. Section \ref{sec:nse} is devoted to the nearly self-similar blowup of the rescaled Navier--Stokes model with constant viscosity. 
Some concluding remarks are made in Section \ref{sec:conclude}. 

\vspace{-0.1in}
\section{Derivation of the axisymmetric Navier--Stokes equations with swirl in $\mathcal{\R}^n$}\label{sec:setup1}




In this section, we derive the axisymmetric Navier--Stokes equations with swirl in both odd and even dimensions. 
A natural way to extend the 3D axisymmetric Navier--Stokes equations to the $n$-dimensional axisymmetric Navier--Stokes equations is to use spherical coordinates in $\R^{n-1}$. 
In \cite{HouZhang2024}, we performed numerical study of the finite time self-similar blowup of the axisymmetric Euler equations with no swirl in three and higher space dimensions using $C^\alpha$ initial vorticity for a wide range of $\alpha$. Denote by ${\bf x} = (x_1,x_2,\cdots,x_{n-1},z)$.
We first define the $n$-dimensional cylindrical unit vectors
\begin{eqnarray*}
{\bf e}_r &=& (\cos \theta_1, \sin \theta_1 \cos \theta_2,...,\sin \theta_1, ...\cos \theta_{n-2},\sin \theta_1...\sin \theta_{n-2},0),\\
{\bf e}_{\theta_1} &=& (-\sin \theta_1, \cos \theta_1 \cos \theta_2,...,\cos \theta_1, ...\cos \theta_{n-2},\cos \theta_1...\sin \theta_{n-2},0),\\
{\bf e}_{\theta_{n-2}} &=& (0, 0,...,-\sin \theta_{n-2},\cos \theta_{n-2},0),\\
{\bf e}_{z} &=& (0, 0,...,1).
\end{eqnarray*}
We assume that the only nontrivial angular velocity is in $\theta_1$, denoted as $u^{\theta_1}$ and set $u^{\theta_j} \equiv 0$ for $j=2,3,...n-2$. 
We call the velocity field ${\bf u}$ axisymmetric if it admits the following expression:
\[
{\bf u} = u^r(t,r,z) {\bf e}_r + u^{\theta_1}(t,r,z){\bf e}_{\theta_1} + u^z(t,r,z) {\bf e}_z .
\]
Using the calculus on curvilinear coordinate, we obtain \cite{HouZhang2024}
\begin{equation}
\label{div-free-spherical-coord}
\nabla \cdot {\bf u} =
\frac{(r^{n-2} u^r)_r}{r^{n-2}} + \frac{(n-3)\cot (\theta_1)}{r} u^{\theta_1} + \frac{(r^{n-2} u^z)_z}{r^{n-2}} .
\end{equation}
Due to the appearance of $\cot(\theta_1)$ in front of  $u^{\theta_1}$ in the expression for $\nabla \cdot {\bf u}$, the divergence free condition cannot be satisfied if $u^r$, $u^{\theta_1}$ and $u^z$ are a function of $r$ and $z$ only. This seems to suggest that one cannot obtain high dimensional axisymmetric Euler or Navier--Stokes equations unless we impose $u^{\theta_1} \equiv 0$, which implies no swirl, if we adopt spherical coordinates to parameterize $\mathbb{S}^{n-2}$.

\vspace{-0.1in}
\subsection{Axisymmetric Navier--Stokes equations in odd dimensions}
\label{odd-dim-case}

To overcome this difficulty, we propose to use a different parametrization of $\mathbb{S}^{n-2}$ in the odd dimensional case with $n=2m+1$.
Denote ${\bf x}=(x_1,x_2, \cdots, x_{2m},z)$. We first define the three orthogonal unit vectors as follows:
\begin{eqnarray*}
{\bf e}_r &=& r^{-1}(x_1,x_2,x_3,\cdots, x_{2m},0), \quad r =\sqrt{x_1^2+x_2^2+\cdots +x_{2m}^2} \;,\\
{\bf e}_{\theta} &=& r^{-1}(-x_2,x_1,-x_4,x_3,\cdots, -x_{2m},x_{2m-1},0),\\
{\bf e}_{z} &=& (0, 0,...,1).
\end{eqnarray*}
The axisymmetric velocity field is expressed in terms of ${\bf e}_r$, ${\bf e}_{\theta}$ and ${\bf e}_{z}$ as follows:
\[
{\bf u} = u^r(t,r,z) {\bf e}_r + u^{\theta}(t,r,z){\bf e}_{\theta} + u^z(t,r,z) {\bf e}_z .
\]
This choice of ${\bf e}_\theta$ shares some essential properties similar to those of 3D cylindrical coordinates:
\begin{equation}
\label{basic-prop}
({\bf e}_\theta \cdot \nabla ) {\bf e}_\theta = -\frac{1}{r} {\bf e}_r, \quad
({\bf e}_\theta \cdot \nabla ) {\bf e}_r = \frac{1}{r} {\bf e}_\theta, \quad
({\bf e}_r \cdot \nabla ) {\bf e}_\theta = 0, \quad
({\bf e}_r \cdot \nabla ) {\bf e}_r = 0 .
\end{equation}

Another important property is that the angular component of the velocity field satisfies the divergence free condition exactly:
\begin{equation}
\label{div-free}
\nabla \cdot (u^\theta {\bf e}_\theta) = 0.
\end{equation}
As a consequence, the divergence free condition is reduced to
\begin{equation}
\nabla \cdot {\bf u} =
\frac{(r^{n-2} u^r)_r}{r^{n-2}} + \frac{(r^{n-2} u^z)_z}{r^{n-2}} = 0,
\end{equation} 
which is the same as the no-swirl case in the $n$-dimensional setting using spherical coordinates.

In order to derive the axisymmetric Navier--Stokes equations in $n$ dimensions, we need to decompose $({\bf u} \cdot \nabla ){\bf u}$ into the ${\bf e}_r$, ${\bf e}_\theta$ and ${\bf e}_z$ directions. By direct calculations using \eqref{basic-prop}, we can prove that
\begin{eqnarray}
\label{conv-term}
({\bf u} \cdot \nabla) {\bf u} = 
\left (u^r \partial_r u^r + u^z \partial_z u^r - \frac{(u^\theta)^2}{r} \right ) {\bf e}_r 
+ \left (u^r \partial_r u^\theta + u^z \partial_z u^\theta + \frac{u^r}{r} u^\theta \right ) {\bf e}_\theta 
+ (u^r \partial_r u^z + u^z \partial_z u^z){\bf e}_z ,
\end{eqnarray}
and 
\begin{subequations}
\label{diffus-term}
\begin{align}
\Delta {\bf u} &= 
\left (u^r_{rr} + \frac{n-2}{r} u^r_r - \frac{n-2}{r^2}u^r + u^r_{zz} \right ) {\bf e}_r + \left (u^\theta_{rr} + \frac{n-2}{r} u^\theta_r - \frac{n-2}{r^2}u^\theta + u^\theta_{zz} \right ) {\bf e}_\theta \\
&+\left (u^z_{rr} + \frac{n-2}{r} u^z_r  + u^z_{zz} \right ){\bf e}_z .
\end{align}
\end{subequations}
Using \eqref{conv-term} and \eqref{diffus-term}, we can derive the following equivalent system for $u^r$, $u^\theta$ and $u^z$:
\begin{subequations}
\label{axisym-primitive}
\begin{align}
&u^r_t + u^r \partial_r u^r + u^z \partial_z u^r - \frac{(u^\theta)^2}{r} = - p_r + \nu\left(u^r_{rr} + \frac{n-2}{r}u^r_{r} -\frac{n-2}{r^2}u^r + u^r_{zz}\right) ,\label{eq:as_NSE_ur}\\ 
&u^\theta_t + u^r \partial_r u^\theta + u^z \partial_z u^\theta + \frac{u^r u^\theta}{r}  = \nu\left(u^\theta_{rr} + \frac{n-2}{r}u^\theta_{r} - \frac{n-2}{r^2}u^\theta + u^\theta_{zz}\right) ,\label{eq:as_NSE_uth}\\ 
&u^z_t + u^r \partial_r u^z + u^z \partial_z u^z  = - p_z + \nu\left(u^z_{rr} + \frac{n-2}{r}u^z_{r} + u^z_{zz}\right) .\label{eq:as_NSE_uz}
\end{align}
\end{subequations}
The functions $u^r$, $u^\theta$ and $u^z$ can be extended as an odd or even function of $r$ and $z$. We will choose a smooth initial condition $u^r(0,r,z)=u^r_0(r,z)$, $u^\theta(0,r,z)=u^\theta_0(r,z)$ and
$u^z(0,r,z)=u^z_0(r,z)$ that satisfies the following symmetry properties: 

\begin{enumerate}
\item $u^r_0(r,z)$ is an odd function of $r$ and an even function of $z$; 

\item $u^\theta_0(r,z)$ is an odd function of both $r$ and $z$;

\item
$u^z_0(r,z)$ is an even function of $r$ and an odd function of $z$; 

\item
$p_0(r,z)$ is an even function of both $r$ and $z$. 
\end{enumerate}
It is easy to show that these symmetry properties are preserved dynamically by equations \eqref{eq:as_NSE_ur}--\eqref{eq:as_NSE_uz}. Since $u^r$ and $u^\theta$ are odd functions of $r$, we conclude that ${\bf u} = u^r(t,r,z) {\bf e}_r + u^{\theta}(t,r,z){\bf e}_{\theta} + u^z(t,r,z) {\bf e}_z$ is a smooth function of ${\bf x}$ as long as $(u^r,u^\theta,u^z)$ remains smooth as a function of $(r,z)$.

Similar to the $3$D axisymmetric Euler or Navier--Stokes equations, we introduce the angular vorticity $\omega^\theta$ and angular stream function $\psi^\theta$ by defining
$\omega^\theta = u_z^r  - u_r^z $ and 
$- \left(\partial_r^2+\frac{n-2}{r}\partial_r - \frac{n-2}{r^2}+\partial_z^2\right) \psi^\theta = \omega^\theta$. Using $\psi^\theta$, we can define $u^r$ and $u^z$ in terms of $\psi^\theta$ as follows:
\begin{equation}
\label{vel-def}
u^r = - \frac{(r^{n-2} \psi^\theta)_z}{r^{n-2}}, \quad 
u^z = \frac{(r^{n-2} \psi^\theta)_r}{r^{n-2}}.
\end{equation}
The divergence free condition 
$(r^{n-2} u^r )_r + (r^{n-2} u^z )_z = 0 $
is automatically satisfied. 
In terms of the vorticity and stream function formulation, we obtain the exact $n$-dimensional axisymmetric Navier--Stokes equations as follows:
\begin{subequations}\label{eq:axisymmetric_NSE_nD}
\begin{align}
&u^\theta_{t}+u^r u^\theta_{r}+u^z u^\theta_{z} = -\frac{u^r u^\theta}{r} + \nu \left(u^\theta_{rr} + \frac{n-2}{r}u^\theta_{r} - \frac{n-2}{r^2}u^\theta + u^\theta_{zz}\right),\label{eq:as_NSE_nD_a}\\
&\om^\theta_t+u^r\om^\theta_r+u^z\om^\theta_z =\frac{((u^\theta)^2)_z}{r} + (n-2) \frac{u^r}{r} \om^\theta + \nu\left(\om^\theta_{rr} + \frac{n-2}{r}\om^\theta_{r} - \frac{n-2}{r^2}\om^\theta + \om^\theta_{zz}\right) ,\label{eq:as_NSE_nD_b}\\
 & -\left(\partial_r^2+\frac{n-2}{r}\partial_r - \frac{n-2}{r^2}+\partial_z^2\right)\psi^\theta = \om^\theta . \label{eq:as_NSE_nD_c}
\end{align}
\end{subequations}
We can derive \eqref{eq:axisymmetric_NSE_01} from \eqref{eq:axisymmetric_NSE_nD} by using $\Gamma=r u^\theta$, $\omega_1=\omega^\theta/r$ and $\psi_1=\psi^\theta/r$. Note that the symmetry properties of $(u^r,u^\theta,u^z)$ imply that $\om^\theta$, $u^\theta$ and $\psi^\theta$ are odd functions of $r$ and $z$.

Building on our exact axisymmetric Navier--Stokes equations in odd dimensions, we define our generalized axisymmeytric Navier--Stokes equations by allowing $n$ to take any positive real value. We can show that the generalized axisymmetric Navier--Stokes equations still conserve the energy, see the derivation in the next subsection. When $n=3$, we recover the $3$D axisymmetric Navier--Stokes.

To derive the generalized rescaled Navier--Stokes model \eqref{eq:axisymmetric_NSE_01}, we restrict ourselves to the 3D axisymmetric Navier--Stokes equations but rescale the axial velocity $u^z$ by $R(t)/Z(t)$ to keep the scaling balance between $u^r$ and $u^z$. We define the velocity as  $u^r=-r\psi_{1,z}, \; u^z =(m-1)\psi_1 + r\psi_{1,r}$ with $m=1+2R(t)/Z(t)$ and ${n}=2m-3$. However, rescaling $u^z$ by $R/Z$ with the choice $m = 1 + 2R/Z$ leads to an $m$-dimensional axisymmetric velocity field $(u^r, u^z)$. Since we omit the term $-(m - 3)\psi_{1,z} \omega_1$ on the right-hand side of \eqref{eq:as_NSE_0_b} when we restrict ourself to $3$-dimension, the resulting system does not correspond to a true $m$-dimensional Navier--Stokes equation. For this reason, we call this model the {\it rescaled Navier--Stokes model}. When ${n}=3$, we recover the $3$D axisymmetric Navier--Stokes. We will show that this rescaled model satisfies
 a generalized energy equality in Section \ref{energy-conserv}. 

\subsection{Axisymmetric Navier--Stokes equations in even dimensions}
\label{even-dim-case}
Next, we derive the axisymmetric Navier--Stokes in even dimensions. Denote $n=2m+2$ and ${\bf x}=(x_1,x_2, \cdots, x_{2m},z_1,z_2)$. We use two approaches to derive the axisymmetric Navier--Stokes equations in even dimensions. In the first approach, we define ${\bf e}_r$ and ${\bf e}_{\theta}$ in the same way as in the odd-dimensional case and use ${\bf e}_{z_1}$ and ${\bf e}_{z_2}$ as the Cartesian unit vectors for the $z_1$, $z_2$ coordinates as follows:
\begin{eqnarray*}
{\bf e}_r &=& r^{-1}(x_1,x_2,x_3,\cdots, x_{2m},0,0)\;, \quad r =\sqrt{x_1^2+x_2^2+\cdots +x_{2m}^2} \;,\\
{\bf e}_{\theta} &=& r^{-1}(-x_2,x_1,-x_4,x_3,\cdots, -x_{2m},x_{2m-1},0,0)\;,\\
{\bf e}_{z_1} &=& (0, 0,...,1,0), \quad {\bf e}_{z_2} = (0, 0,...,0,1).
\end{eqnarray*}
The axisymmetric velocity field is expressed in terms of ${\bf e}_r$, ${\bf e}_{\theta}$, ${\bf e}_{z_1}$ and ${\bf e}_{z_2}$ as follows:
\[
{\bf u} = u^r(t,r,z_1,z_2) {\bf e}_r + u^{\theta}(t,r,z_1,z_2){\bf e}_{\theta} +
u^{z_1}(t,r,z_1,z_2){\bf e}_{z_1} + u^{z_2}(t,r,z_1,z_2) {\bf e}_{z_2} .
\]
The axisymmetry will correspond to the symmetry plane $(z_1,z_2)$ instead of a single symmetry $z$-axis in the odd dimensional case. Since ${\bf e}_r$ and ${\bf e}_{\theta}$ enjoy the same properties as in the odd-dimensional case, we can apply the same derivation as in the odd-dimensional case to derive the axisymmetric Navier--Stokes equations in the even-dimensional case and obtain the following system for $u^r$, $u^\theta$, $u^{z_1}$ and $u^{z_2}$ below
\begin{eqnarray*}
&&u^r_t + u^r \partial_r u^r + u^{z_1} \partial_{z_1} u^r + u^{z_2} \partial_{z_2} u^r - \frac{(u^\theta)^2}{r} = - p_r + \nu\left(u^r_{rr} + \frac{n-3}{r}u^r_{r} -\frac{n-3}{r^2}u^r + u^r_{z_1z_1}+ u^r_{z_2z_2}\right),\\ 
&&u^\theta_t + u^r \partial_r u^\theta + u^{z_1} \partial_{z_1} u^\theta + u^{z_2} \partial_{z_2} u^\theta + \frac{u^r u^\theta}{r}  = \nu\left(u^\theta_{rr} + \frac{n-3}{r}u^\theta_{r} - \frac{n-3}{r^2}u^\theta + u^\theta_{z_1z_1}+u^\theta_{z_2z_2}\right) ,\\ 
&&u^{z_1}_t + u^r \partial_r u^{z_1} + u^{z_1} \partial_{z_1} u^{z_1} + u^{z_2} \partial_{z_2} u^{z_1} = - p_{z_1} + \nu\left(u^{z_1}_{rr} + \frac{n-3}{r}u^{z_1}_{r} + u^{z_1}_{z_1 z_1}+u^{z_1}_{z_2 z_2}\right) ,\\
&&u^{z_2}_t + u^r \partial_r u^{z_2} + u^{z_1} \partial_{z_1} u^{z_2} + u^{z_2} \partial_{z_2} u^{z_2} = - p_{z_2} + \nu\left(u^{z_2}_{rr} + \frac{n-3}{r}u^{z_2}_{r} + u^{z_2}_{z_1 z_1}+u^{z_2}_{z_2 z_2}\right) .
\end{eqnarray*}
The pressure is determined by enforcing the following divergence free condition:
\begin{equation}
\frac{(r^{n-3} u^r)_r}{r^{n-3}} + \frac{(r^{n-3} u^{z_1})_{z_1}}{r^{n-3}}
+  \frac{(r^{n-3} u^{z_2})_{z_2}}{r^{n-3}} = 0.
\end{equation} 

We will choose smooth initial data $u^r(0,r,z_1,z_2)=u^r_0(r,z_1,z_2)$, $u^\theta(0,r,z_1,z_2)=u^\theta_0(r,z_1,z_2)$,
$u^{z_1}(0,r,z_1,z_2)=u^{z_1}_0(r,z_1,z_2)$ and
$u^{z_2}(0,r,z_1,z_2)=u^{z_2}_0(r,z_1,z_2)$ that satisfy the following symmetry properties: 

\begin{enumerate}
\item $u^r_0(r,z_1,z_2)$ is an odd function of $r$ and an even function of $z_1$ and $z_2$; 

\item $u^\theta_0(r,z_1,z_2)$ is an odd function of $r$, $z_1$ and $z_2$;

\item
$u^{z_1}_0(r,z_1,z_2)$ is an even function of $r$ and $z_2$ but an odd function of $z_1$; 

\item
$u^{z_2}_0(r,z_1,z_2)$ is an even function of $r$ and $z_1$ but an odd function of $z_2$; 

\item
$p_0(r,z_1,z_2)$ is an even function of $r$, $z_1$ and $z_2$. 
\end{enumerate}
It is easy to show that these symmetry properties are preserved dynamically by equations \eqref{eq:as_NSE_ur}--\eqref{eq:as_NSE_uz}. Since $u^r$ and $u^\theta$ are odd functions of $r$, we conclude that ${\bf u} = u^r(t,r,z_1,z_2) {\bf e}_r + u^{\theta}(t,r,z_1,z_2){\bf e}_{\theta} + u^{z_1}(t,r,z_1,z_2) {\bf e}_{z_1} + u^{z_2}(t,r,z_1,z_2) {\bf e}_{z_2}$ is a smooth function of ${\bf x}$ as long as $(u^r,u^\theta,u^{z_1},u^{z_2})$ remains smooth as a function of $(r,z_1,z_2)$.

One disadvantage of the first approach is that the solution depends on $(r,z_1,z_2)$, which is 3-dimensional. In the second approach, we apply the generalized polar coordinates for $(z_1,z_2)$ to further reduce it to a 2-dimensional problem. Specifically, we define the slightly modified four unit vectors
\begin{eqnarray*}
{\bf e}_r &=& r^{-1}(x_1,x_2,x_3,\cdots, x_{2m},0,0)\;, \quad r =\sqrt{x_1^2+x_2^2+\cdots +x_{2m}^2},\\ 
{\bf e}_{\theta} &=& r^{-1}(-x_2,x_1,-x_4,x_3,\cdots, -x_{2m},x_{2m-1},0,0)\;,\\
{\bf e}_{\phi}&=&\left(0,0,\cdots,0,-\frac{z_2}{z},\frac{z_1}{z}\right)\;,\quad
         {\bf e}_z=(0,0,\dots,0,\frac{z_1}{z},\frac{z_2}{z}), \quad z = \sqrt{x_{2m+1}^2+x_{2m+2}^2}.
\end{eqnarray*}
The axisymmetric velocity field is expressed in terms of ${\bf e}_r$, ${\bf e}_{\theta}$, ${\bf e}_{\phi}$ and ${\bf e}_{z}$ as follows:
\[
{\bf u} = u^r(t,r,z) {\bf e}_r + u^{\theta}(t,r,z){\bf e}_{\theta} +
u^{\phi}(t,r,z){\bf e}_{\phi} + u^z(t,r,z) {\bf e}_z .
\]
This choice of ${\bf e}_\theta$ and $ {\bf e}_\phi$ shares properties similar to those of the cylindrical and the polar coordinates:
\begin{eqnarray}
\label{basic-prop-2}
&&({\bf e}_\theta \cdot \nabla ) {\bf e}_\theta = -\frac{1}{r} {\bf e}_r, \quad
({\bf e}_\theta \cdot \nabla ) {\bf e}_r = \frac{1}{r} {\bf e}_\theta, \quad
({\bf e}_r \cdot \nabla ) {\bf e}_\theta = 0, \quad
({\bf e}_r \cdot \nabla ) {\bf e}_r = 0 ,\\
&& ({\bf e}_\phi \cdot \nabla ) {\bf e}_\phi = -\frac{1}{z} {\bf e}_z, \quad
({\bf e}_\phi \cdot \nabla ) {\bf e}_z = \frac{1}{z} {\bf e}_\phi, \quad
({\bf e}_z \cdot \nabla ) {\bf e}_\phi = 0, \quad
({\bf e}_z \cdot \nabla ) {\bf e}_z = 0 .
\end{eqnarray}
Another important property is that the angular components of the velocity field satisfy the divergence free condition exactly:
\begin{equation}
\label{div-free}
\nabla \cdot (u^\theta {\bf e}_\theta) = 0, \quad
\nabla \cdot (u^\phi {\bf e}_\phi) = 0 .
\end{equation}

We can extend the derivation for the odd dimensional axisymmetric Navier--Stokes to derive the following  equations for $u^r$, $u^\theta$, $u^z$, and $u^\phi$ for the even dimensional axisymmetric Navier-Stokes:
\begin{subequations}\label{eq:axisymmetric_even}
\begin{align}
    &u^r_t+u^r u^r_r+u^z u^r_z=\frac{({u^\theta})^2}{r}-p_r+\nu\Delta u^r-\nu\frac{(n-3)u^r}{r^2},\label{ur-eqn-even}\\
    &u^z_t+u^r u^z_r+u^z u^z_z=\frac{({u^\phi})^2}{z}-p_z+\nu\Delta u^z-\nu\frac{u^z}{z^2},\label{uz-eqn-even}\\
    &u^\theta_t+u^r u^\theta_r+u^z u^\theta_z=-\frac{u^ru^\theta}{r}+\nu\Delta u^\theta-\nu\frac{(n-3)u^\theta}{r^2},\label{uthe-eqn-even}\\
    &u^\phi_t+u^r u^\phi_r+u^z u^\phi_z=-\frac{u^zu^\phi}{z}+\nu\Delta u^\phi-\nu\frac{u^\phi}{z^2},\label{uphi-eqn-even}\\
    &\nabla \cdot {\bf u} = u^r_r+\frac{n-3}{r}u^r+u^z_z+\frac{1}{z}u^z=\frac{1}{r^{n-3}}(r^{n-3} u^r)_r+\frac{1}{z}(zu^z)_z=0 , \label{div-eqn-even}
\end{align}
\end{subequations}
where 
\[
\Delta f(r,z)=f_{rr}(r,z)+\frac{n-3}{r} f_r(r,z)+f_{zz}(r,z)+\frac{1}{z} f_z(r,z)\;.
\]

We will impose the following symmetry properties on the initial data: 

\begin{enumerate}
\item Both $u^r_0(r,z)$ and $u^\theta_0(r,z)$ are odd functions of $r$ and even functions of $z$; 

\item Both $u^\phi_0(r,z)$ and $u^z_0(r,z)$ are even functions of $r$ and odd functions of $z$;

\item
$p_0(r,z)$ is an even function of both $r$ and $z$. 
\end{enumerate}
These symmetry properties are preserved dynamically by equations \eqref{ur-eqn-even}--\eqref{div-eqn-even}. Since $(u^r,u^\theta)$ is an odd function of $r$ and an even function of $z$, while $(u^\phi,u^z)$ is an odd function of $z$ and an even function of $r$, we conclude that ${\bf u} = u^r {\bf e}_r + u^{\theta}{\bf e}_{\theta} + u^\phi {\bf e}_\phi +  u^z {\bf e}_z$ is a smooth function of ${\bf x}$ as long as $(u^r,u^\theta,u^\phi,u^z)$ remains smooth as a function of $(r,z)$.

The divergence free condition implies that there exists a stream function $\psi^\theta (r,z)$ that satisfies
\[
    -(r^{n-3}z\psi^\theta)_z=r^{n-3}zu^r,\quad
     (r^{n-3}z\psi^\theta)_r=r^{n-3}zu^z,
\]
which is equivalent to
\[
    u^r=-\frac{1}{z}(z\psi^\theta)_z,\quad
    u^z=\frac{1}{r^{n-3}}(r^{n-3}\psi^\theta)_r.
\]
Define $\omega^\theta=u^r_z-u^z_r$. Then $\om^\theta$ and $\psi^\theta$ satisfy
\[
     \om^\theta_t+u^r \om^\theta_r+u^z\om^\theta_z    
    =\frac{2 u^\theta u^\theta_z}{r}-\frac{2u^\phi u^\phi_r}{z}+\left(\frac{(n-3)u^r}{r}+\frac{u^z}{z}\right)\om^\theta+\nu\Delta \om^\theta-\nu\frac{(n-3)\om^\theta}{r^2}-\nu\frac{\om^\theta}{z^2},
\]
\[
-\left(\psi^\theta_{zz}+\psi^\theta_{rr}+\frac{n-3}{r}\psi^\theta_r-\frac{n-3}{r^2}\psi^\theta+\frac{1}{z}\psi^\theta_z-\frac{1}{z}\psi^\theta\right)=\om^\theta. \quad \quad \quad\quad\quad\quad\quad\quad
\quad\quad \quad\quad
\]
In summary, the axisymetric Navier--Stokes equations in even space dimensions are given below:
\begin{equation}
    \begin{aligned}       
        &u^\theta_t+u^r u^\theta_r+u^z u^\theta_z=-\frac{u^ru^\theta}{r}+\nu\left(\Delta -\frac{n-3}{r^2}\right) u^\theta,\\
        &u^\phi_t+u^r u^\phi_r+u^z u^\phi_z=-\frac{u^zu^\phi}{z}+\nu\left(\Delta -\frac{1}{z^2}\right)u^{\phi},\\
         & \om^\theta_t+u^r \om^\theta_r+u^z\om^\theta_z=\frac{2 u^\theta u^\theta_z}{r}-\frac{2u^\phi u^\phi_r}{z}+\left(\frac{(n-3)u^r}{r}+\frac{u^z}{z}\right)\om^\theta+\nu\left(\Delta -\frac{n-3}{r^2}-\frac{1}{z^2}\right)\om^\theta,\\
        &-\left(\Delta-\frac{n-3}{r^2}-\frac{1}{z^2}\right)\psi^\theta=\om^\theta,       
    \end{aligned}
\end{equation}
where $u^r=-\frac{1}{z}(z\psi^\theta)_z, \quad u^z=\frac{1}{r^{n-3}}(r^{n-3}\psi^\theta)_r$.
The symmetry properties of $(u^r,u^\theta,u^\phi,u^z)$ imply that $u^\theta$ is odd in $r$ and even in $z$ while $(\om^\theta,\psi^\theta)$ is an odd function of both $r$ and $z$. We can easily extend the local well-posedness analysis for the 3D axisymmetric Navier-Stokes equations to prove a similar local well-posedness result for the high dimensional axisymmetric Navier-Stokes equations. 

\subsection{Visualization of $ \mathbb{R}^4$ and $ \mathbb{R}^5$ .}
\label{Illustration}
\begin{figure}[!ht]
\centering
    \includegraphics[width=0.5\textwidth]{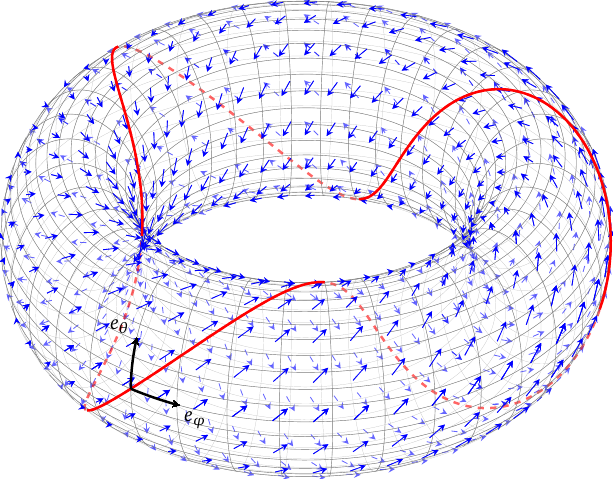}
    \caption[Torus-R4]{ Illustration of a velocity field on a Torus. Blue vectors: \( u^\theta {\bf e}_\theta + u^\varphi {\bf e}_\varphi \) at each point, which is always tangent to the surface of the torus. Here for a fixed \( (r, z) \), both \( u^\theta \) and \( u^\varphi \) are independent of the position on the torus surface. Red curve: one particular streamline of this velocity field. } 
    \label{fig:Torus}
\end{figure}

Taking $ \mathbb{R}^4 $ as an example, we actually view $ \mathbb{R}^4$ as $\mathbb{R}^2 \times \mathbb{R}^2$, and use polar coordinates in each $ \mathbb{R}^2 $ to describe the flow. By fixing $ (r, z) $ (where $ r > 0 $ and $ z > 0$), we obtain a cross-section of $ \mathbb{R}^4 $, which is a torus $ \mathbb{T}^2$. In Fig. \ref{fig:Torus}, we may assume that the radius of the torus $\mathbb{T}^2$ is $z$ while the radius of the cross section of $\mathbb{T}^2$ is $r$. The flow on this torus is described by the components $ u^\theta$ and $ u^\phi$, which represent the tangential velocity fields along the angular directions $\theta$ and $\phi$ within the torus. These components capture the motion of the fluid on the same torus at a given $(r, z)$ point. In addition to this, there is also fluid motion between different tori, which is described by the radial components $ u^r$ and $ u^z$, corresponding to the flow in the radial directions $ r$ and $ z$. These components govern the movement of the fluid from one torus to another, thus capturing the motion both within each torus and across different tori in $ \mathbb{R}^4$.

In the case of $ \mathbb{R}^5$, we view $\mathbb{R}^5$ as $ \mathbb{R}^4 \times \mathbb{R}$. The vector $\mathbf{e}_\theta$ can be considered as the unit tangent vector to some circles around the $ z$-axis. Given $ (r\cos \alpha_1 \cos \alpha_2, r\cos \alpha_1 \sin \alpha_2, r\sin \alpha_1 \cos \alpha_3, r\sin \alpha_1 \sin \alpha_3, z_0) \in\mathbb{R}^5$, the circle passing through this point is described by:
\[
\begin{aligned}
    x_1(\theta) &= r \cos \alpha_1 \cos(\alpha_2 + \theta), & x_2(\theta) &= r \cos \alpha_1 \sin(\alpha_2 + \theta), \\
    x_3(\theta) &= r \sin \alpha_1 \cos(\alpha_3 + \theta), & x_4(\theta) &= r \sin \alpha_1 \sin(\alpha_3 + \theta), & z(\theta) &= z_0.
\end{aligned}
\]
This is a parametric description of a circle. We can rewrite it as:
\[
\begin{aligned}
    &(x_1(\theta), x_2(\theta), x_3(\theta), x_4(\theta), z(\theta)) \\
    &= (0, 0, 0, 0, z_0) + \big( r \cos \alpha_1 \cos \alpha_2, r \cos \alpha_1 \sin \alpha_2, r \sin \alpha_1 \cos \alpha_3, r \sin \alpha_1 \sin \alpha_3,0 \big) \cos \theta \\
    &\quad + \big( -r \cos \alpha_1 \sin \alpha_2, r \cos \alpha_1 \cos \alpha_2, -r \sin \alpha_1 \sin \alpha_3, r \sin \alpha_1 \cos \alpha_3 ,0\big) \sin \theta.
\end{aligned}
\]
Hence, it is indeed a circle of radius $r$ on the plane passing through the point $(0,0,0,0,z_0)$ and spanned by ${\bf e}_r$ and ${\bf e}_\theta$ at 
$(x_1(0), x_2(0), x_3(0), x_4(0), z(0)) $ . 
The derivative of this parametrization is:
\[
(x_1'(\theta), x_2'(\theta), x_3'(\theta), x_4'(\theta), z'(\theta)) = (-x_2(\theta), x_1(\theta), -x_4(\theta), -x_3(\theta), 0) \varparallel \mathbf{e}_\theta.
\]
Therefore, $\mathbf{e}_\theta$ is the unit tangent vector of the circles described above. As $ \theta$ changes from 0 to $ 2\pi$, the point moves along the circle, rotating once around the $ z$-axis in $\mathbb{R}^5$. It is easy to show that every point in $ \mathbb{R}^5$, except those on the $ z$-axis, lies on one of these circles, and that different circles do not intersect. Therefore, fixing $(r, z)$ with $r >0$, these circles form a 3-dimensional sphere $\mathbb{S}^3$. To visualize this, we can plot the following figure, showing that the circles form a torus $\mathbb{T}^2$, which is a cross-section of $ \mathbb{S}^3$. The component $u^\theta$ represents the flow field on $\mathbb{S}^3$ corresponding to a specific $(r, z)$, where the fluid moves around the aforementioned circles. Additionally, the components $ u^r$ and $ u^z$ describe the fluid flow between different $ \mathbb{S}^3$'s corresponding to different $(r,z)$. These components capture both the motion within each sphere and the movement of the fluid between different spheres in $ \mathbb{R}^5$.

\begin{figure}[!ht]
\centering
    \includegraphics[width=0.5\textwidth]{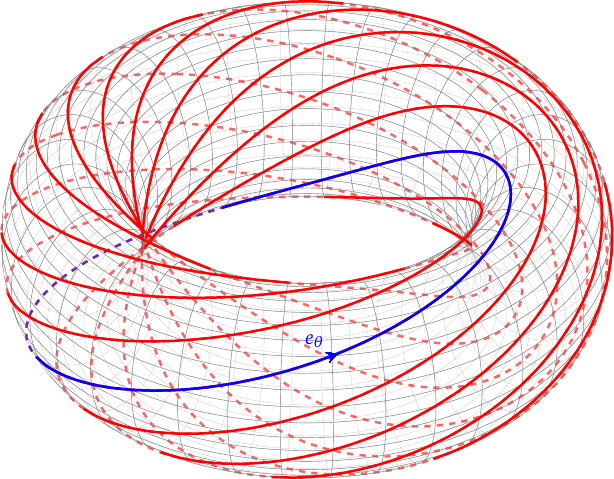}
    \caption[Circle-R5]{ Illustration of circles forming a Torus. } 
    \label{fig:Torus2}
\end{figure}

\subsection{Proofs of some key properties}
In this subsection, we will prove the key properties \eqref{basic-prop},\eqref{div-free}, \eqref{conv-term} and
\eqref{diffus-term} in the odd dimensional case with $n=2m+1$.

\vspace{0.05in}
\noindent
{\bf Proof of \eqref{basic-prop}}:  We will only prove ${\bf e}_\theta \cdot \nabla {\bf e}_\theta = - r^{-1} {\bf e}_r$. The other properties can be proved similarly. 
By definition of ${\bf e}_\theta$, we have for any $1 \leq  j \leq m$
\begin{eqnarray*}
({\bf e}_\theta \cdot \nabla) (-x_{2j}/r) &=&
\sum_{i=1}^m \left ( (-x_{2i}/r)\partial_{x_{2i-1}} (-x_{2j}/r) + (x_{2i-1}/r)\partial_{x_{2i}} (-x_{2j}/r) \right )\\
& = &  -x_{2j-1}/r^2 +
\sum_{i=1}^m \left ( -\frac{x_{2j} x_{2i-1}x_{2i}}{r^4} + \frac{x_{2j} x_{2i-1}x_{2i}}{r^4} \right ) = -x_{2j-1}/r^2 \;,
\end{eqnarray*}
\begin{eqnarray*}
({\bf e}_\theta \cdot \nabla) (x_{2j-1}/r) &=&
\sum_{i=1}^m \left ( (-x_{2i}/r)\partial_{x_{2i-1}} (x_{2j-1}/r) + (x_{2i-1}/r)\partial_{x_{2i}} (x_{2j-1}/r) \right )\\
& = &  -x_{2j}/r^2 +
\sum_{i=1}^m \left ( \frac{x_{2j-1} x_{2i-1}x_{2i}}{r^4} - \frac{x_{2j-1} x_{2i-1}x_{2i}}{r^4} \right ) = -x_{2j}/r^2 \;.
\end{eqnarray*}
This proves ${\bf e}_\theta \cdot \nabla {\bf e}_\theta = - r^{-1} {\bf e}_r$.

\vspace{0.05in}
\noindent
{\bf Proof of \eqref{div-free}.}
By definition, we have
\begin{eqnarray*}
\nabla \cdot (u^\theta {\bf e}_\theta ) & = &
\sum_{i=1}^m \left (\partial_{x_{2i-1}} ( u^\theta (-x_{2i}/r)) + \partial_{x_{2i}} ( u^\theta (x_{2i-1}/r)) \right ) \\
& = & \sum_{i=1}^m \partial_r (u^\theta)\left ( \frac{\partial r}{\partial_{x_{2i-1}}} (-x_{2i}/r) +  \frac{\partial r}{\partial_{x_{2i}}} (x_{2i-1}/r)\right ) + \sum_{i=1}^m u^\theta 
\left ( \partial_{x_{2i-1}} (-x_{2i}/r) +  \partial_{x_{2i}} (x_{2i-1}/r) \right ) \\
& = & \sum_{i=1}^m \partial_r (u^\theta)\left ( \frac{-x_{2i-1}x_{2i}} {r^2} + \frac{x_{2i}x_{2i-1}} {r^2} \right ) + \sum_{i=1}^m u^\theta 
\left ( \frac{x_{2i}x_{2i-1}}{r^3} - \frac{x_{2i-1}x_{2i}}{r^3} \right ) = 0\; .
\end{eqnarray*}

\vspace{0.05in}
\noindent
{\bf Proof of \eqref{conv-term}}. Direct calculations give 
\begin{eqnarray*}
({\bf u} \cdot \nabla ) {\bf u} & = &
(u^r {\bf e}_r + u^\theta {\bf e}_\theta + u^z {\bf e}_z) \cdot \nabla (u^r {\bf e}_r + u^\theta {\bf e}_\theta + u^z {\bf e}_z) \\
&=& u^r ({\bf e}_r\cdot \nabla)(u^r {\bf e}_r) 
+ u^r ({\bf e}_r\cdot \nabla)(u^\theta {\bf e}_\theta)
+ u^r ({\bf e}_r\cdot \nabla)(u^z {\bf e}_z)\\
&+& u^\theta ({\bf e}_\theta \cdot \nabla)(u^r {\bf e}_r) 
+ u^\theta ({\bf e}_\theta \cdot \nabla)(u^\theta {\bf e}_\theta)
+ u^\theta ({\bf e}_\theta \cdot \nabla)(u^z {\bf e}_z)\\
&+&u^z ({\bf e}_z\cdot \nabla)(u^r {\bf e}_r) 
+ u^z ({\bf e}_z\cdot \nabla)(u^\theta {\bf e}_\theta)
+ u^z ({\bf e}_z\cdot \nabla)(u^z {\bf e}_z) .
\end{eqnarray*}
Using \eqref{basic-prop} and the fact that $u^r$, $u^\theta$ and $u^z$ are functions of $r$ and $z$ only, we obtain
\begin{eqnarray*}
({\bf u} \cdot \nabla ) {\bf u} & = &
(u^r (\partial_r u^r) {\bf e}_r
+ u^r (\partial_r u^\theta) {\bf e}_\theta
+ u^r (\partial_r u^z) {\bf e}_z ) + ((u^\theta u^r) r^{-1} {\bf e}_\theta -
(u^\theta)^2 r^{-1} {\bf e}_r) \\
&+& u^z (\partial_z u^r) {\bf e}_r
+ u^z (\partial_z u^\theta) {\bf e}_\theta
+u^z (\partial_z u^z) {\bf e}_z \\
&=& \left (u^r \partial_r u^r + u^z \partial_z u^r - 
\frac{(u^\theta)^2}{r} \right ) {\bf e}_r + \left (u^r \partial_r u^\theta + u^z \partial_z u^\theta  +
\frac{u^\theta u^r}{r} \right ) {\bf e}_\theta
+ (u^r \partial_r u^z + u^z \partial_z u^z) {\bf e}_z .
\end{eqnarray*}

\vspace{0.05in}
\noindent
{\bf Proof of \eqref{diffus-term}}. Define $\hat{\bf e}_r = r{\bf e}_r$ and $\hat{\bf e}_\theta = r{\bf e}_\theta$. Note that $\hat{\bf e}_r$ and $\hat{\bf e}_\theta$ are linear in $x_i$ ($i=1,2,\cdots, n-1$). Thus we have $\Delta \hat{\bf e}_r =0$ and $\Delta \hat{\bf e}_\theta = 0$. Direct calculations give
\begin{eqnarray*}
\Delta {\bf u} &=& \Delta \left (\frac{u^r}{r} \hat{\bf e}_r \right ) +  \Delta \left (\frac{u^\theta}{r} \hat{\bf e}_\theta \right ) + \Delta u^z {\bf e}_z\\
&=& \Delta \left (\frac{u^r}{r}\right )\hat{\bf e}_r + 2 \partial_r \left ( \frac{u^r}{r} \right ) {\bf e}_r
+\Delta \left (\frac{u^\theta}{r}\right )\hat{\bf e}_\theta + 2 \partial_r \left ( \frac{u^\theta}{r} \right ) {\bf e}_\theta + \Delta u^z {\bf e}_z \\
&=& \left ( r \Delta \left (\frac{u^r}{r}\right ) + 2 \partial_r \left ( \frac{u^r}{r} \right )\right ) {\bf e}_r
+\left ( r \Delta \left (\frac{u^\theta}{r}\right ) + 2 \partial_r \left ( \frac{u^\theta}{r} \right )\right )  {\bf e}_\theta + \Delta u^z {\bf e}_z \; .
\end{eqnarray*}
Using the well-known property of the $n$-dimensional Laplacian for a scalar function $f(r,z)$:
\[
\Delta f(r,z) = f_{rr} + \frac{n-2}{r} f_r + f_{zz},
\]
we can simplify the above formula to prove \eqref{diffus-term}.

\subsection{Energy estimates for the generalized Navier--Stokes equations}
\label{energy-conserv}

In this subsection, we will derive the energy estimates for both the generalized Navier--Stokes equations and the rescaled Navier--Stokes model for a given constant dimension $n$. We first rewrite the $\Gamma$-equation in terms of $u_1 = u^\theta/r$ as follows:
\begin{equation}
\label{eqn-nse-u1}
u_{1,t}+u^ru_{1,r}+u^zu_{1,z} =2u_1\psi_{1,z} + \nu_1\left(u_{1,rr} + \frac{n}{r}u_{1,r} + u_{1,zz} \right).
\end{equation}
To perform energy estimates for the generalized Navier--Stokes equations \eqref{eq:axisymmetric_NSE_nD}, we multiply the $u_1$-equation \eqref{eqn-nse-u1} by $u_1 r^n dr dz$ and integrate over the whole domain. Similarly, we multiply equation \eqref{eq:as_NSE_nD_b} for $\omega_1$ and equation \eqref{eq:as_NSE_nD_c} for $\psi_1$ by $\psi_1r^n dr dz$, and integrate over the whole domain. Using the divergence form of the diffusion operator, 
\[
u_{1,rr} + \frac{n}{r}u_{1,r} + u_{1,zz} = \left (\frac{ (r^n u_{1,r})_r}{r^n} + u_{1,zz} \right ),
\]
we perform integration by parts to obtain
\begin{eqnarray*}
&&\int u_1 (u_{1,rr} + \frac{n}{r}u_{1,r} + u_{1,zz}) \;r^n dr dz = - \int |\nabla u_1 |^2 \;r^n dr dz,\\
&&-\int \psi_1 (\psi_{1,rr} + \frac{n}{r}\psi_{1,r} + \psi_{1,zz}) \;r^n dr dz =  \int |\nabla \psi_1 |^2 \;r^n dr dz,\\
&&\int \psi_1 (\omega_{1,rr} + \frac{n}{r}\omega_{1,r} + \omega_{1,zz}) \;r^n dr dz =  -\int |\Delta \psi_1 |^2 \;r^n dr dz,
\end{eqnarray*}
where we have used $-\Delta \psi_1 = \omega_1$. There are no contributions from the boundary when we perform integration by parts since $\psi_1$, $u_1$ and $\omega_1$ are odd functions of $z$ and even functions of $r$.  Using $u^r = - (r^{n-1}\psi_1)_z/r^{n-2}$, $u^z =  (r^{n-1}\psi_1)_r/r^{n-2}$ and the incompressibility condition, we get
\begin{eqnarray*}
&&\int u_1 (-u^r u_{1,r}-u^z u_{1,z} + 2\psi_{1,z}u_1) r^n dr dz \\
&& = 
\int \left (\frac{1}{2}(r^{n-1}\psi_1)_z (u_1^2)_r - \frac{1}{2}(r^{n-1}\psi_1)_r (u_1^2)_z \right ) r^2 drdz + \int 2 \psi_{1,z}u_1^2 r^n dr dz =  \int \psi_{1,z}u_1^2 \;r^n dr dz.
\end{eqnarray*}

Similarly, we obtain by integration by parts that

\begin{eqnarray*}
&&\int \psi_1 (-u^r \omega_{1,r}-u^z \omega_{1,z} + (u_1^2)_z - (n-3) \psi_{1,z}\omega_1) r^n dr dz \\
&& = 
\int \left (\psi_1 (r^{n-1}\psi_1)_z \omega_{1,r} - \psi_1(r^{n-1}\psi_1)_r \omega_{1,z} \right ) r^2 drdz + \int (-\psi_{1,z}u_1^2 - (n-3)\psi_1\psi_{1,z}\omega_1) r^n dr dz \\
&& =  (n-3) \int \psi_1 \psi_{1,z}\omega_1 \;r^n dr dz - \int (\psi_{1,z}u_1^2 + (n-3)\psi_1\psi_{1,z}\omega_1) r^n dr = - \int \psi_{1,z}u_1^2 \;r^n dr dz.
\end{eqnarray*}

By adding the above two estimates, we observe that the right hand side terms cancel each other. Thus, we have

\begin{eqnarray*}
\frac{1}{2}\frac{d}{dt} \int (u_1^2 + |\nabla \psi_1|^2) r^n dr dz = -\nu \int (  |\nabla u_1 |^2 + |\Delta \psi_1|^2) \;r^n dr dz.
\end{eqnarray*}
In terms of the original physical variables, we get by a direct computation that 
\[
\int (u_1^2 + |\nabla \psi_1|^2) r^n dr dz = \int (|u^\theta|^2 + |u^r|^2 + u^z|^2) r^{n-2} dr dz,
\]
which is the familiar kinetic energy $\int |{\bf u}|^2 r^{n-2} dr dz$ in $n$ dimensions. The damping term from the viscous term can be shown to be equivalent to $-\nu \int (|\nabla {\bf u}|^2) r^{n-2} dr dz$.

We can perform a similar energy estimate for the rescaled Navier--Stokes model. As we can see from our above energy estimate, the term $-(n-3)\psi_{1,z}\omega_1$ from the $\omega_1$ equation plays an important role in canceling the contribution from the advection terms. Since the term $-(n-3)\psi_{1,z}\omega_1$ is absent from  the $\omega_1$ equation in the rescaled Navier--Stokes model, we need to change the weight in our energy estimate to cancel the contribution from the advection terms. We proceed as follows:
\begin{eqnarray*}
&&\int u_1 (-u^r u_{1,r}-u^z u_{1,z} + 2\psi_{1,z}u_1) r^n dr dz \\
&& = 
\int \left (\frac{1}{2}(r^{m-1}\psi_1)_z (u_1^2)_r - \frac{1}{2}(r^{m-1}\psi_1)_r (u_1^2)_z \right ) r^{n-m+2} drdz + \int 2 \psi_{1,z}u_1^2 r^n dr dz\\
&& =  \frac{(5-m)}{2}\int \psi_{1,z}u_1^2 \;r^n dr dz =  \frac{(7-n)}{4}\int \psi_{1,z}u_1^2 \;r^n dr dz .
\end{eqnarray*}

Similarly, we obtain by integration by parts that
\begin{eqnarray*}
&&\int \psi_1 (-u^r \omega_{1,r}-u^z \omega_{1,z} + (u_1^2)_z ) r^n dr dz \\
&& = 
\int \left (\psi_1 (r^{m-1}\psi_1)_z \omega_{1,r} - \psi_1(r^{m-1}\psi_1)_r \omega_{1,z} \right ) r^{n-m+2} drdz - \int \psi_{1,z}u_1^2  r^n dr dz \\
&& =  (2m-3-n) \int \psi_1 \psi_{1,z}\omega_1 \;r^n dr dz - \int \psi_{1,z}u_1^2 \; r^n dr
 = - \int \psi_{1,z}u_1^2 \;r^n dr dz,
\end{eqnarray*}
where we have used $m=(n+3)/2$ to cancel the first term on the right hand side in the second to the last step. Since we assume $n < 7$, we obtain the following generalized energy estimate by forming a linear combination of the two energy terms to cancel the right hand sides
\begin{eqnarray*}
\frac{1}{2}\frac{d}{dt} \int (u_1^2 + \frac{(7-n)}{4} |\nabla \psi_1|^2) r^n dr dz = -\nu_1 \int  |\nabla u_1 |^2 \;r^n dr dz - \nu_2 \frac{(7-n)}{4} \int  |\Delta \psi_1|^2 \;r^n dr dz.
\end{eqnarray*}
We can further express $\int (u_1^2 + \frac{(7-n)}{4} |\nabla \psi_1|^2) r^n dr dz$ as
$\int ((u^\theta)^2 + \frac{(7-n)}{4} (|u^r|^2+|u^z|^2) r^{n-2} dr dz$.

\section{A novel two-scale dynamic rescaling formulation}
\label{sec:setup2}

In this subsection, we introduce a novel two-scale dynamic rescaling formulation to study potential nearly self-similar blowup solution of the generalized Navier--Stokes equations. 
The dynamic rescaling formulation was introduced in \cite{mclaughlin1986focusing,  landman1988rate} to study the self-similar blowup of the nonlinear Schr\"odinger equations. This formulation is also called the modulation technique in the literature and has been developed by Merle, Raphael, Martel, Zaag, Soffer, Weinstein, and others. It has been a very effective tool to analyze the formation of singularities for many problems like the nonlinear Schr\"odinger equation \cite{soffer90,kenig2006global,merle2005blow}, compressible Euler equations \cite{buckmaster2019formation,buckmaster2019formation2}, the nonlinear wave equation \cite{merle2015stability}, the nonlinear heat equation \cite{merle1997stability}, the generalized KdV equation \cite{martel2014blow}, and other dispersive problems. 
Recently, this method has been applied to study singularity formation in incompressible fluids \cite{chen2019finite2,elgindi2019finite} .

For other nonlinear PDEs such as the nonlinear Schr\"odinger or Keller-Segel system, one can eliminate these unstable modes by using the symmetry properties of the solution and studying the spectral properties of the compact linearized operator around an explicit ground state, see e.g. \cite{merle2005blow,Collot2019SpectralAF}. In our case, we do not have an explicit ground state and the linearized operator is not compact. We need to enlarge the solution space by lifting the space dimension above $3$, and varying the space dimension dynamically to enforce the scaling balance between the advection along the $r$ and $z$ directions. This effectively eliminates the scaling instability and leads to one-scale blowup. 

We first define the following dynamically rescaled profiles and the rescaled time variable $\tau$.
	\begin{eqnarray*}
  \widetilde{u}_1 (\tau, \xi,\eta) = C_{u}(\tau)u_1\left(t(\tau), C_{lr}(\tau)\xi, C_{lz}(\tau)\eta\right), \\
  \widetilde{\omega}_1 (\tau, \xi,\eta) = C_{\omega}(\tau)\omega_1\left(t(\tau), C_{lr}(\tau)\xi, C_{lz}(\tau)\eta \right), \\
 \widetilde{\psi}_1 (\tau, \xi,\eta) = C_{\psi}(\tau)\psi_1\left(t(\tau), C_{lr}(\tau)\xi, C_{lz}(\tau)\eta \right),
  \end{eqnarray*}
where $ \xi = r/C_{lr} $, $\eta = z/C_{lz} $, and $u_1 =u^\theta/r$,
\begin{equation}
\label{scaling-parameter-Cu}
 C_{u}(\tau)=e^{\int_0^\tau c_{u}(s)\mathrm{d}s},\; C_\omega(\tau)=e^{\int_0^\tau c_\omega(s)\mathrm{d}s}, \;
 C_\psi(\tau)=e^{\int_0^\tau c_\psi(s)\mathrm{d}s},
 \end{equation}
 and
 \begin{equation}
\label{scaling-parameter-Cl}
 C_{lr}(\tau)=e^{-\int_0^\tau c_{lr}(s)\mathrm{d}s}, \;  C_{lz}(\tau)=e^{-\int_0^\tau c_{lz}(s)\mathrm{d}s},
\; t(\tau)=\int_0^\tau C_{\psi}(s) C_{lz} (s)\mathrm{d}s.
\end{equation}
Here $\tau$ is the rescaled time variable satisfying $d\tau/dt = ( C_{\psi} C_{lz} )^{-1}.$
Then, the generalized $n$-dimensional axisymmetric Navier--Stokes equations in the $\widetilde{\psi}_1, \; \widetilde{u}_1,\; \widetilde{\omega}_1$ variables can be described by the following two-scale dynamic rescaling equations
\begin{subequations}
\label{Dyn-rescale-eqn-nD}
\begin{align}
& \widetilde{u}_{1,\tau}+c_{lr} \xi \widetilde{u}_{1,\xi} + c_{lz}\eta \widetilde{u}_{1,\eta}+\widetilde{\bf u} \cdot \nabla_{(\xi,\eta)} \widetilde{u}_1 = c_u \widetilde{u}_1 + 2\widetilde{u}_1\widetilde{\psi}_{1,\eta} + \frac{\nu_1 C_\psi}{C_{lz}} \Delta \widetilde{u} ,\\
&\widetilde{\omega}_{1,\tau}+c_{lr} \xi \widetilde{\omega}_{1,\xi} + c_{lz}\eta \widetilde{\omega}_{1,\eta} + \widetilde{\bf u} \cdot \nabla_{(\xi,\eta)} \widetilde{\omega}_1 = c_\omega \widetilde{\omega}_1 + (\widetilde{u}_1^2)_\eta - (n-3)\widetilde{\psi}_{1,\eta} \widetilde{\omega}_1 + \frac{\nu_2 C_\psi}{C_{lz}}\Delta \widetilde{\omega} ,\\
& - \Delta \widetilde{\psi}_1 = \widetilde{\omega}_1 ,  \quad
\Delta = \left(\delta^2\partial_{\xi}^2+\delta^2 \frac{n}{\xi}\partial_{\xi} + \partial_{\eta}^2\right),
\end{align}
\end{subequations}
where $\delta = C_{lz}(\tau)/C_{lr}(\tau)$ and the rescaled velocity field is given by $\widetilde{\bf u} = (\tilde{u}^\xi, \tilde{u}^\eta)$, $\tilde{u}^\xi = -\xi \widetilde{\psi}_{1,\eta}$, $\tilde{u}^\eta = (n-1)\widetilde{\psi}_1 + \xi \widetilde{\psi}_{1,\xi}$, and the scaling parameters ($c_{lz}, c_{lr},c_\psi, c_u, c_\omega)$ satisfy the rescaling relationship
\begin{equation}
\label{scaling-relation-cu}
c_\psi = c_u+c_{lz}, \quad c_\omega = c_u - c_{lz}.
\end{equation}

In the case of the rescaled Navier--Stokes model, we drop $(n-3)\widetilde{\psi}_{1,\eta} \widetilde{\omega}_1$ from the $\widetilde{\omega}_1$ equation, and change the velocity field to $u^r= -(r^{m-2}\psi^\theta)_z/r^{m-2}, \; u^z =(r^{m-2}\psi^\theta)_r/r^{m-2}$ with $m=(n+3)/2$.

We have {\bf three free scaling parameters $c_{lr}$, $c_{lz}$ and $c_u$ to choose to enforce the normalization conditions. Using the two-scale dynamic rescaling formulation and introducing space dimension as a new degree of freedom are two key ingredients that enable us to eliminate the scaling instability}. If we use the traditional one-scale dynamic rescaling formulation with $c_{lr}=c_{lz}$, we could not eliminate this scaling instability.  We enforce the normalization conditions that $\widetilde{u}_1$ achieves its maximum at $(\xi,\eta) = (R_0,1)$ with $R_0 = 3.6927$ and $\|\widetilde{u}_1\|_\infty$ being fixed to be $1$. We remark that $(R_\omega,Z_\omega)$ converges to a fixed position $(R_\omega,Z_\omega)$ as $\tau $ increases. 

If the scaling parameters converge to a constant value as $\tau \rightarrow \infty$, 
\[(c_{lz}(\tau), c_{lr}(\tau),c_\psi(\tau), c_u(\tau), c_\omega(\tau)) \rightarrow (c_{lz}, c_{lr},c_\psi, c_u, c_\omega),  \
\] 
we can obtain the actual blowup rate in the physical time variable by inverting the mapping from $\tau$ back to $t$. To simplify the derivation, we assume that $(c_{lz}(\tau), c_{lr}(\tau),c_\psi(\tau), c_u(\tau), c_\omega(\tau))$ are time independent. Then we obtain the following asymptotic scaling results:
\[
C_\psi(\tau)=e^{c_\psi \tau }, \quad
 C_{lz}(\tau)=e^{-c_{lz} \tau},
 \]
 which implies that
 \[
 t(\tau) = \int_0^\tau C_\psi(s)C_{lz}(s) ds = \frac{1}{c_\psi - c_{lz}} (e^{(c_\psi - c_{lz})\tau} - 1).
 \]
By inverting this relationship, we obtain
\[
\tau = \frac{1}{c_{lz}-c_\psi}\log \left ( \frac{1}{T-t} \right ),
\]
where $T = \frac{1}{c_{lz} - c_\psi}$. By substituting $\tau = \frac{|\log(T-t)|}{c_{lz}-c_\psi}$ back to $C_u$, $C_\omega$, $C_{lz}$ etc and using $c_\psi = c_u+c_{lz}$ and $ c_\omega = c_u - c_{lz}$, we obtain the following scaling formula:
\[
C_{lz} (\tau) = (T-t)^{\widehat{c}_{lz}}, \quad
C_{lr} (\tau) = (T-t)^{\widehat{c}_{lr}},
\]
where $\widehat{c}_{lz} = c_{lz}/(c_{lz}-c_\psi)$ and $\widehat{c}_{lr} = c_{lr}/(c_{lz}-c_\psi)$. The blowup rates are given by
\[
\frac{1}{C_u (\tau)} = \frac{1}{(T-t)}, \quad
\frac{1}{C_\omega (\tau)} = \frac{1}{(T-t)^{1+\widehat{c}_{lz}}}, \quad
\frac{1}{C_\psi (\tau)} = \frac{1}{(T-t)^{1-\widehat{c}_{lz}}},
\]
where $\widehat{c}_\omega = c_\omega/(c_{lz}-c_\psi) = -1 - \widehat{c}_{lz}$ and $\widehat{c}_\psi = c_\psi/(c_{lz}-c_\psi) = -1 + \widehat{c}_{lz}$.
In our computation, we monitor closely these normalized scaling exponents $\widehat{c}_\omega,\; \widehat{c}_\psi, \; \widehat{c}_{lz}$ and $\widehat{c}_{lr}$ to study their scaling properties.

As we mentioned before, we will use the conservative $(\Gamma, \omega_1, \psi_1)$ formulation in our computation. Using the relationship $\Gamma = r^2 u_1$, we can obtain an equivalent dynamic rescaling formulation for $(\Gamma,\omega_1, \psi_1)$ as follows:
\begin{subequations}
\label{Dyn-rescale-eqn-Gamma}
\begin{align}
& \widetilde{\Gamma}_{\tau}+c_{lr} \xi \widetilde{\Gamma}_{\xi} + c_{lz}\eta \widetilde{\Gamma}_{\eta}+\widetilde{\bf u} \cdot \nabla_{(\xi,\eta)} \widetilde{\Gamma} = c_\Gamma \widetilde{\Gamma} + \frac{\nu_1 C_\psi}{C_{lz}} \widetilde{\Delta} \widetilde{\Gamma} ,\\
&\widetilde{\omega}_{1,\tau}+c_{lr} \xi \widetilde{\omega}_{1,\xi} + c_{lz}\eta \widetilde{\omega}_{1,\eta}+\widetilde{\bf u} \cdot \nabla_{(\xi,\eta)} \widetilde{\omega}_1 = c_\omega \widetilde{\omega}_1 + (\widetilde{u}_1^2)_\eta + (3-n)\widetilde{\psi}_{1,\eta} \widetilde{\omega}_1 + \frac{\nu_2 C_\psi}{C_{lz}}\Delta \widetilde{\omega} ,\\
& - \Delta \widetilde{\psi}_1 = \widetilde{\omega}_1 , \quad \Delta =\left(\delta^2\partial_{\xi}^2+\delta^2 \frac{n}{\xi}\partial_{\xi} + \partial_{\eta}^2\right), 
\end{align}
\end{subequations}
where $c_\Gamma = c_u + 2c_{lr} = 2(c_{lr}-c_{lz})$ and
\begin{equation}
\label{Gamma-Delta}
\widetilde{\Delta}\widetilde{\Gamma} = \left(\delta^2\partial_{\xi}^2\widetilde{\Gamma} +\delta^2 \frac{n-4}{\xi}\partial_{\xi}\widetilde{\Gamma} + \delta^2 \frac{6-2n}{\xi^2} \widetilde{\Gamma}+ \partial_{\eta}^2 \widetilde{\Gamma}\right).
\end{equation}
We still choose $c_{lr}$ and $c_{lz}$ to enforce that $\widetilde{u}_1$ achieves its maximum at $(\xi,\eta) = (R_0,1)$ and choose $c_u$ to fix the maximum value of $\widetilde{u}_1$ to be $1$.

In the case of the generalized  Navier--Stokes equations with solution dependent viscosity $\nu (\tau)= \nu_0 \|u_1\|_\infty Z(t)^2$, we note that $ \|u_1\|_\infty Z(t)^2 =  \| \widetilde{u}_1 \|_\infty = C_{lz}/C_\psi$. Thus, we have 
$\frac{\nu C_{\psi}}{C_{lz}} = \nu_0$.
As a result, the dynamic rescaling formulation is the same as that for the generalized axisymmetric Navier--Stokes equations with a constant viscosity $\nu_0$. If $c_{lz} \rightarrow c_l$ and $c_{lr} \rightarrow c_l$, then we have $\delta(\tau) \rightarrow \lambda_0$. In our computation, we obtain $\lambda_0  \approx 0.914$. If the dynamic rescaled solution converges to a steady state, we can rescale the $\xi$ variable to $\xi/\lambda_0$ to obtain a one-scale solution with $\delta = 1$ and $c_{lr}=c_{lz} \equiv c_l$.  

If the solution of the dynamic rescaling equations converges to a steady state, the steady state satisfies the following self-similar equations: 
\begin{eqnarray}
\label{Dyn-rescale-eqn3new}
&& (c_{l} \xi, c_{l}\eta) \cdot\nabla_{(\xi,\eta)} \widetilde{\Gamma}+\widetilde{\bf u} \cdot \nabla_{(\xi,\eta)} \widetilde{\Gamma} =  \nu_0  \widetilde{\Delta} \widetilde{\Gamma} ,\\
&&(c_{l} \xi, c_{l}\eta) \cdot\nabla_{(\xi,\eta)} \widetilde{\omega}_{1,\eta}+\widetilde{\bf u} \cdot \nabla_{(\xi,\eta)} \widetilde{\omega}_1 = c_\omega \widetilde{\omega}_1 + (\widetilde{u}_1^2)_\eta - (n-3)\widetilde{\psi}_{1,\eta} \widetilde{\omega}_1 +\nu_0 \Delta \widetilde{\omega} ,\\
&& - \Delta \widetilde{\psi}_1 = \widetilde{\omega}_1 ,  \quad
\Delta = -\left(\partial_{\xi}^2+\frac{n}{\xi}\partial_{\xi} + \partial_{\eta}^2\right),
\end{eqnarray}
where we have used $c_\Gamma = 2(c_{lr}-c_{lz}) = 0$. 
In our study, we observe that $c_l$ and $n$ decrease as $\nu_0$ decreases.
It would be interesting to solve the self-similar equations directly with $\nu_0$ as a continuation parameter. Due to the total circulation conservation of the generalized Euler equations, we expect to have $c_l \rightarrow 1/2$ as $\nu_0\rightarrow 0$. We also observe that the dimension $n$ decreases toward $n=3$ as we decrease $\nu_0$, but the solution suffers from the Kelvin-Helmholtz instability when $\nu_0$ is below certain threshold.

\subsection{The importance of resolving the far field solution}

We remark that capturing the correct far field decay rate is essential to capture the correct scaling properties of the solution. Resolving the far field solution also plays an important role in capturing the slow growth rate of the $L^n$ norm of the velocity and the dynamic growth of $\| r u^r\|_\infty$.  On the other hand, resolving the near field is also important since we need to compute accurately the location of the maximum of $V_1$ and its amplitude in order to enforce our normalization conditions. 
It is also worth emphasizing the importance of using the conservative $(\Gamma,\omega_1,\psi_1)$ formulation in our two-scale dynamic rescaling formulation. This enables us to capture the nearly self-similar blowup. 

In our computation, we expand the domain size by $R(\tau)^{-1/5}$ and $Z(\tau)^{-1/5}$ in the $r$ and $z$ directions, respectively and apply the homogeneous Neumann boundary conditions for the stream function $\widetilde{\psi}_1$, $\widetilde{\Gamma}$ and $\widetilde{\omega}_1$. In the case of the rescaled Navier--Stokes model with two constant viscosity coefficients, we have $(R(\tau),Z(\tau)) = ( 3.3\cdot 10^{-15},  6.2\cdot 10^{-16})$ by the end of our computation.
The domain size has increased to $[0,3 \cdot 10^7]\times [0, 1.65 \cdot 10^7]$ from the initial domain size of $O(10^{4})$.

\vspace{-0.1in}
\subsection{The operator splitting strategy}

We will adopt an operator splitting strategy developed in \cite{HouZhang2024} to enforce the normalization conditions. To enforce the normalization conditions accurately at every time step, we utilize the operator splitting method.  We denote by ${\bf v} = (\widetilde{\Gamma}, \widetilde{\omega}_1)$. We will split the time evolution of ${\bf v}$ into two parts:
\[
{\bf v}_{\tau} =F({\bf v})+G({\bf v}),
\]
where $F({\bf v})$ contains the original terms in the generalized Navier--Stokes equations and $G({\bf v})$ contains the linear terms that control the rescaling, i.e. $G({\bf v}) = -c_{lr} \xi {\bf v}_\xi - c_{lz} \eta {\bf v}_\eta + c_{\bf v} {\bf v}$. We view $\widetilde{\psi}_1$ as a function of $\widetilde{\omega}_1$ through the Poisson equation. The operator splitting method allows us to solve the dynamic rescaling formulation by solving ${\bf v}_{\tau} =F({\bf v})$ and
${\bf v}_{\tau} =G({\bf v})$ alternatively.
We can use the forward Euler method to solve ${\bf v}_{\tau} =F({\bf v})$. In the second step, we can obtain a closed form solution to ${\bf v}_{\tau} =G({\bf v})$ as follows:
\[
 {\bf v}(\tau, \xi, \eta) = C_{\bf v}(\tau) {\bf v}(0, C_{lr} \xi, C_{lz} \eta),
 \] 
where $C_{\bf v} = \exp \left ( \int_0^\tau c_{\bf v}(s) ds \right )$,
$C_{lr} = \exp \left ( -\int_0^\tau c_{lr}(s) ds \right )$ and 
$C_{lz} = \exp \left ( -\int_0^\tau c_{lz}(s) ds \right )$.
In the first step, solving ${\bf v}_{\tau} =F({\bf v})$ will violate the normalization conditions. But we will correct this error in the second step by solving ${\bf v}_{\tau} =G({\bf v})$ with a smart choice of $C_{\bf v},C_{lr}$ and $C_{lz}$. In other words, at every time step when we solve ${\bf v}_{\tau} =G({\bf v})$, we can exactly enforce the normalization conditions of fixing the location of the maximum of $\widetilde{u}_1$ to be at $(R_0,1)$ by rescaling the $\xi$ and $\eta$ coordinates. We could also adopt Strang's splitting method to improve the splitting scheme to second order accuracy.


\section{Blowup of the generalized Navier-Stokes equations with solution dependent viscosity}
\label{sec:euler}


We first review our previous adaptive mesh computation for the potential singular behavior of the 3D axisymmetric Euler and Navier--Stokes equations in \cite{Hou-euler-2022,Hou-nse-2022}. In this study, we used the following initial condition (recall $u_1 =u^\theta/r$ and $\omega_1 = \omega^\theta/r$):
\begin{equation}
\label{eq:initial-data}
u_1 (0,r,z) =\frac{12000(1-r^2)^{18}
\sin(2 \pi z)}{1+12.5(\sin(\pi z))^2}, \quad \om_1(0,r,z)=0, \quad r \leq 1.
\end{equation}
The flow is completely driven by large swirl initially.
The other two velocity components are set to zero initially. 
Note that $u_1$ is an odd and periodic function of $z$ with period $1$. The oddness of $u_1$ induces the oddness of $\omega_1$ dynamically through the vortex stretching term in the $\omega_1$-equation. 
It is worth emphasizing that the specific power $18$ is important, which determines the ratio of the scales along the $r$ and $z$ directions and enforces a rapid decay near the boundary $r=1$. 
The specific form of the denominator is also important. It breaks the even symmetry of $\sin(2 \pi z)$ with respect to $z=1/4$ along the $z$ direction with a bias toward $z=0$.
This initial condition generates a solution that has comparable scales along the $r$ and $z$ directions, leading to a one-scale traveling solution moving toward the origin. 

We imposed a periodic boundary condition in $z$ with period $1$ and no-slip no-flow boundary condition at $r=1$. Since $u^\theta,\om^\theta,\psi^\theta$ are  odd functions of $r$ \cite{liuwang2006},  
$u_1,\om_1,\psi_1$ are even functions of $r$. Thus, we imposed the following pole conditions:
$u_{1,r}(t,0,z) = \om_{1,r}(t,0,z) = \psi_{1,r}(t,0,z) = 0.$
To numerically compute the potential singularity formation of the axisymmetric Navier--Stokes equations 
with initial condition \eqref{eq:initial-data}, we adopted the numerical methods developed in \cite{Hou-Huang-2022,Hou-Huang-2023}. In particular, we design an adaptive mesh by constructing two adaptive mesh maps for $r$ and $z$ explicitly. The computation was performed in the transformed domain using a uniform mesh. When we map back to the physical domain, we obtain a highly adaptive mesh. We refer to Appendix A in \cite{Hou-Huang-2023} for more detailed discussions.

\subsection{Numerical setup and implementation details}
In this section, we will investigate the asymptotically self-similar blowup of the generalized Navier-Stokes equations with solution dependent viscosity using the two-scale dynamic rescaling formulation. We will use the late stage solution obtained by the adaptive mesh computation in \cite{Hou-nse-2022} using resolution $1536\times 1536$ at the time $T_1=0.0022868$ by which the maximum vortcity has increased by $1.070090e+06$. The maximum of $u_1$ is located at $R(T_1)=1.2311815e-04$ and $Z(T_1) = 3.3302426e-05$. Denote by $C_l = Z(T_1)$ and $V_{M} = \max (u_1(T_1))$. We rescale the solution at $T_1$ using the parabolic scaling invariant property. More precisely, we define $\xi = r/C_l$ and $\eta = z/C_l$, and
\begin{eqnarray*}
\widetilde{u}_1 (\xi, \eta) = u_1(r,z)/V_{M}, \;
\widetilde{\omega}_1 (\xi, \eta) = \omega_1(r,z)/(V_{M}/C_l),\;
\widetilde{\psi}_1 (\xi, \eta) = \psi_1(r,z)/(V_{M} C_l).
\end{eqnarray*}
Through the above rescaling, the maximum of $\widetilde{u}_1$ is equal to 1 and the location at which $\widetilde{u}_1$ achieves its global maximum is rescaled to $\xi = R_0 = 3.692705$ and $\eta =1$. The computational domain is mapped from the original physical domain $[0, 1]\times [0, 0.5]$ to the rescaled domain $[0, 1/C_l]\times [0, 0.5/C_l]$. In our two-scale dynamic rescaling computation, we always rescale the maximum of $\widetilde{u}_1$ to be equal to 1 at a fixed location $(R_0,1)$ for all times. 
We also apply a soft cut-off to the rescaled initial data to suppress the far field tail and use the modified rescaled initial data for the dynamic rescaling formulation. The rescaled initial data can be made available upon request.

We apply an adaptive mesh similar to what we used in \cite{Hou-euler-2022,Hou-nse-2022} to cover the rescaled domain $[0, 1/C_l]\times [0, 0.5/C_l]$. Moreover we use a second order finite difference scheme for space discretization of the $\Gamma$ and $\omega_1$ equations and a second order Runge-Kutta scheme for time discretization with adaptive time step that satisfies the standard CFL condition for the advection terms and the time step constraint for the diffusion terms. A second order finite element method is used to solve the stream function $\psi_1$ to obtain the velocity $u^r$ and $u^z$. 

Due to the use of a highly adaptive mesh and variable coefficients, the finite element discretization of the elliptic equation for the stream function $\psi_1$ produces a stiffness matrix that has a large condition number. We use a direct solver (backslash in Matlab) to invert the stiffness matrix that has a built-in backward error analysis functionality to control the effect of round-off errors. In \cite{ChenHou2023a,ChenHou2023b}, we have used the same solver to solve the dynamical rescaling formulation with a similar adaptive mesh.  We have performed careful resolution study and demonstrated that we can achieve a very small residual error of order $10^{-10}$ in a bounded domain and $10^{-7}$ in the whole domain using a weighted norm. A similar study has been performed in \cite{luo2014toward} where a quadruple precision computation was used to validate the double precision computation when solving the Poisson equation by the finite element method with a similar adaptive mesh.

It turns out that the choice of the viscosity coefficient plays a crucial role in generating a stable and self-similar blowup of the generalized Naver--Stokes equations. 
In our study, we choose $\nu = \nu_0 \|u_1\|_\infty Z(t)^2$ with $\nu_0=0.006$. Here $(R(t),Z(t))$ is the position where $u_1$ achieves its maximum. Note that $\|u_1\|_\infty Z(t)^2$ is scaling invariant. This choice of viscosity is to enforce the balance between the vortex stretching terms and the diffusion terms for both the $u_1$ and $\omega_1$ equations. Another way to interpret this solution dependent viscosity is that it is chosen such that the cell Reynolds number is finite and independent of the small scales of the physical solution.

Using the scaling relationship $c_\psi = c_u + c_{lz}$ from \eqref{scaling-relation-cu}, we can show that 
\[
\|u_1\|_\infty Z(t)^2 = C_{lz}/C_\psi ,
\]
which implies that 
\[
\nu(t) \frac{C_\psi}{C_{lz}} = \nu_0.
\]
Thus, the dynamic rescaling formulation for the generalized Navier--Stokes equations with solution dependent viscosity $\nu(t)=\nu_0 \|u_1\|_\infty Z(t)^2$ reduce to the following dynamic rescaling formulation:
\begin{subequations}
\label{Dyn-rescale-eqn-Gamma-new}
\begin{align}
& \widetilde{\Gamma}_{\tau}+c_{lr} \xi \widetilde{\Gamma}_{\xi} + c_{lz}\eta \widetilde{\Gamma}_{\eta}+\widetilde{\bf u} \cdot \nabla_{(\xi,\eta)} \widetilde{\Gamma} = c_\Gamma \widetilde{\Gamma} + \nu_0 \widetilde{\Delta} \widetilde{\Gamma} ,\\
&\widetilde{\omega}_{1,\tau}+c_{lr} \xi \widetilde{\omega}_{1,\xi} + c_{lz}\eta \widetilde{\omega}_{1,\eta}+\widetilde{\bf u} \cdot \nabla_{(\xi,\eta)} \widetilde{\omega}_1 = c_\omega \widetilde{\omega}_1 + (\widetilde{u}_1^2)_\eta + (3-n)\widetilde{\psi}_{1,\eta} \widetilde{\omega}_1 + \nu_0 \Delta \widetilde{\omega}_1 ,\\
& - \Delta \widetilde{\psi}_1 = \widetilde{\omega}_1 , \quad \Delta =\left(\delta^2\partial_{\xi}^2+\delta^2 \frac{n}{\xi}\partial_{\xi} + \partial_{\eta}^2\right), 
\end{align}
\end{subequations}
where $\widetilde{\bf u} = (\widetilde{u}^\xi,\widetilde{u}^\eta)$, $\widetilde{u}^\xi = -\xi \widetilde{\psi}_{1,\eta}$, $\widetilde{u}^\eta = (n-1) \widetilde{\psi}_1+\xi\widetilde{\psi}_{1,\xi} $, $n=1+2R(\tau)/Z(\tau)$, where $(R(\tau),Z(\tau))$ is where $\widetilde{u}_1=\widetilde{\Gamma}/\xi^2$ achieves its global maximum, $c_\Gamma = c_u + 2c_{lr} = 2(c_{lr}-c_{lz})$ and $\widetilde{\Delta}$ is defined in \eqref{Gamma-Delta}.
We choose $c_{lr}$ and $c_{lz}$ to enforce that $\widetilde{u}_1$ achieves its maximum at $(\xi,\eta) = (R_0,1)$ and choose $c_u$ to fix the maximum value of $\widetilde{u}_1$ to be $1$, and update $c_\psi$ and $c_\omega$ using $c_\psi = c_u +c_{lz}$ and $c_\omega = c_u - c_{lz}$.  

An important consequence of this choice of viscosity is that the self-similar profile of the blowup solution satisfies the generalized self-similar Navier--Stokes equations with {\it constant viscosity coefficient} $\nu_0$. This explains why we can maintain the balance between the vortex stretching term and the diffusion term in the self-similar space variables $\xi=r/C_{lr}$ and $\eta = z/C_{lz}$. 

In the rest of this section, we will present numerical evidences of finite time self-similar blowup of the generalized Navier-Stokes equations 
\eqref{eq:axisymmetric_NSE_0n} with solution dependent viscosity $\nu(t) = \nu_0 \|u_1\|_\infty Z(t)^2$.

\subsection{Rapid growth of maximum vorticity}

In this subsection, we investigate how the profiles of the solution evolve in time. We will use the numerical results computed on the adaptive mesh using resolution $(n_1,n_2) = (1024,1024)$. We have computed the numerical solution up to time $\tau=185$ when the solution is still well resolved. 

\begin{figure}[!ht]
\centering
    \includegraphics[width=0.4\textwidth]{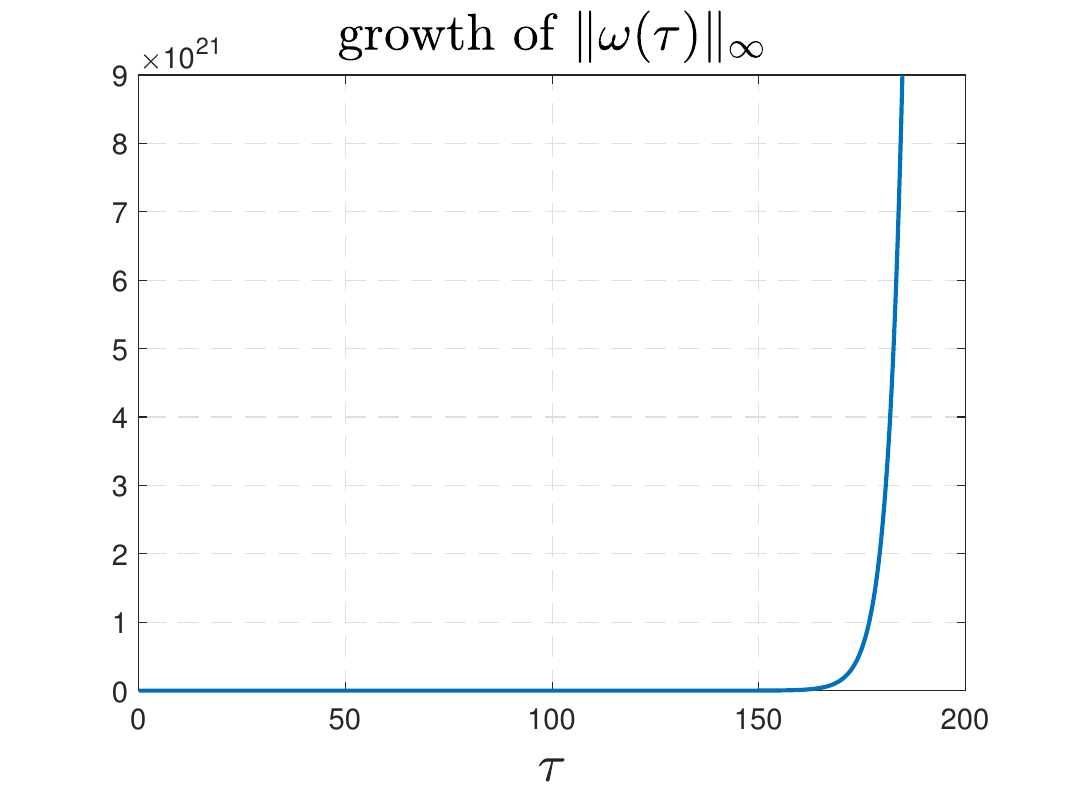}
 \includegraphics[width=0.4\textwidth]{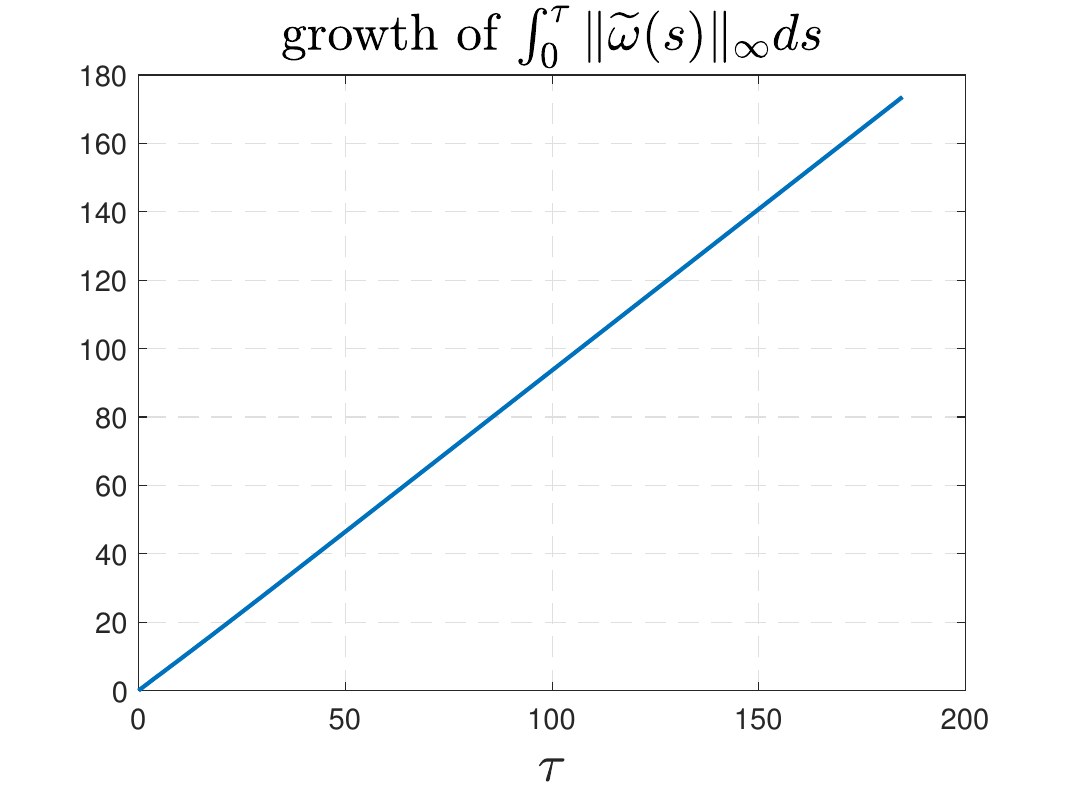} 
    \caption[Lpq-time]{ Left plot: Dynamic growth of the maximum vorticity of the generalized Navier--Stokes equations as a function of $\tau$. Right plot: The growth of $\int_0^\tau\|\widetilde{\vom}\|_{L^\infty} ds$ as a function of $\tau$. The perfect linear fitting implies that $\|\vom(t)\|_{L^\infty} =O\left (\frac{1}{T-t} \right )$.} 
    \label{fig:max-vorticity-decay}
\end{figure}

In Figure \ref{fig:max-vorticity-decay}(a), we plot the dynamic growth of the maximum vorticity as a function of $\tau$. We observe a rapid growth of $\| \vom \|_\infty$. By the end of the computation, the maximum vorticity has increased by a factor of $9\times 10^{21}$. To best of our knowledge, such a large growth rate of the maximum vorticity has not been reported for the $3$D Euler or Navier--Stokes equations in the literature. In Figure \ref{fig:max-vorticity-decay}(b), we plot the time integral of the maximum vorticity $\int_0^t  \| \vom (s)\|_\infty ds = \int_0^\tau  \| \widetilde{ \vom} (s')\|_\infty ds' $. We observe a perfect linear growth of the time integral of the maximum vorticity with respect to $\tau$. Since $\tau = c_0 |\log (T-t)|$, this implies that $\| \vom \|_\infty \sim \frac{1}{T-t}$. This violates the Beale-Kato-Majda blowup criterion \cite{beale1984remarks}.

  \begin{figure}[!ht]
\centering
    \includegraphics[width=0.4\textwidth]{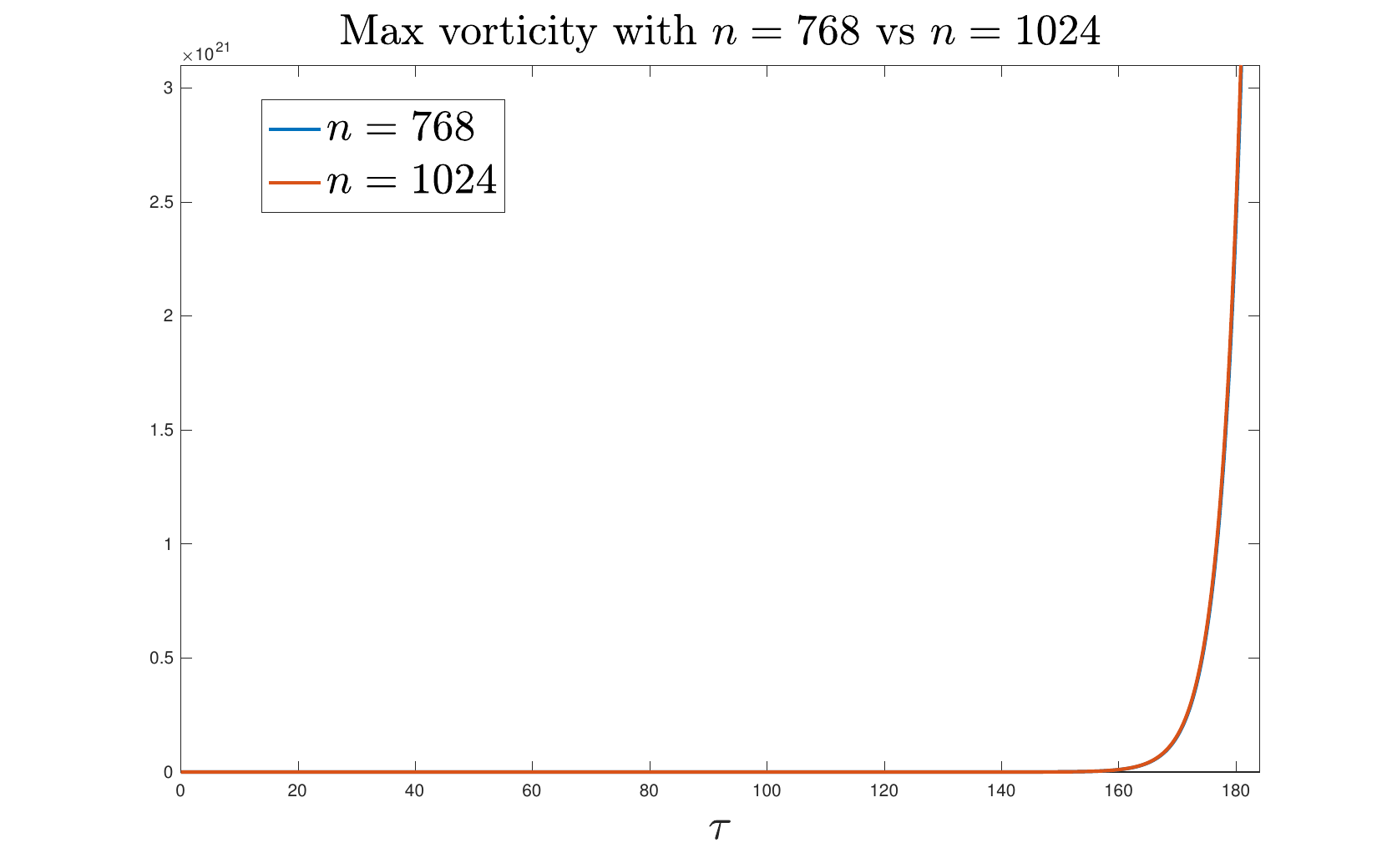}
    \includegraphics[width=0.4\textwidth]{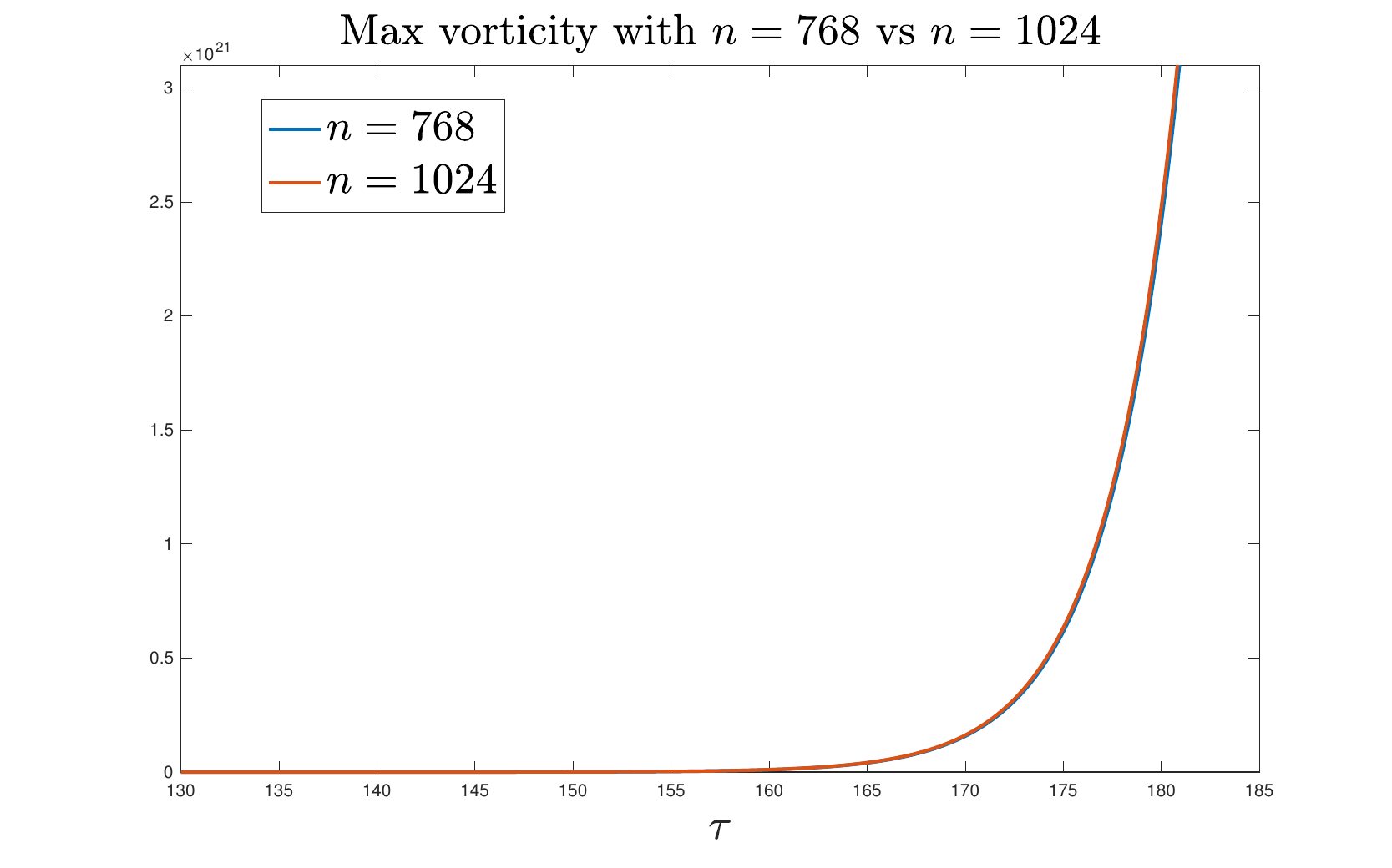} 
    \caption[Max vorticity compare]{ Growth of maximum vorticity for the generalized Navier--Stokes equations. Left plot: Comparison of $\|\vom (\tau)\|_{L^\infty}/\|\vom (0)\|_{L^\infty}$ in time, $n=768$ (blue) vs $n=1024$ (red). Right plot:  Zoomed-in version. Here $n$ stands for the numerical resolution of using an $n\times n$ grid.} 
    \label{fig:max_vort-comparison_nse-decay}
\end{figure}

In Figure \ref{fig:max_vort-comparison_nse-decay} (a), we compare the growth rate of the maximum vorticity using two different resolutions, $768\times 768$ vs $1024\times 1024$. A zoomed-in version is provided in Figure \ref{fig:max_vort-comparison_nse-decay} (b). We can see that the maximum vorticity computed by the resolution $1024\times 1024$ grows slightly faster than that computed by the resolution $768\times 768$. This indicates that the higher resolution captures the growth of the maximum vorticity more accurately, but their difference is very small, indicating that the solution is well resolved by the $1024\times 1024$ grid.

\begin{figure}[!ht]
\centering
    \begin{subfigure}[b]{0.40\textwidth}
        \centering
        \includegraphics[width=1\textwidth]{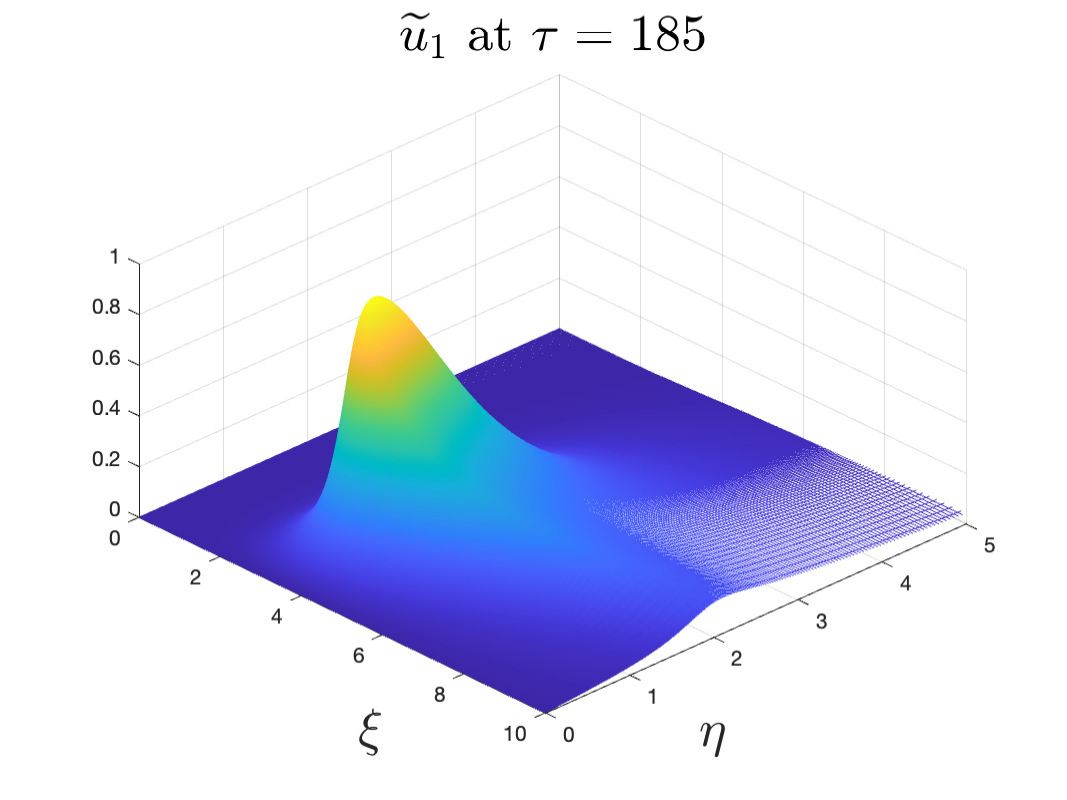}
        \caption{Rescaled profile for $\widetilde{u}_1$}
    \end{subfigure}
    \begin{subfigure}[b]{0.40\textwidth}
        \centering
        \includegraphics[width=1\textwidth]{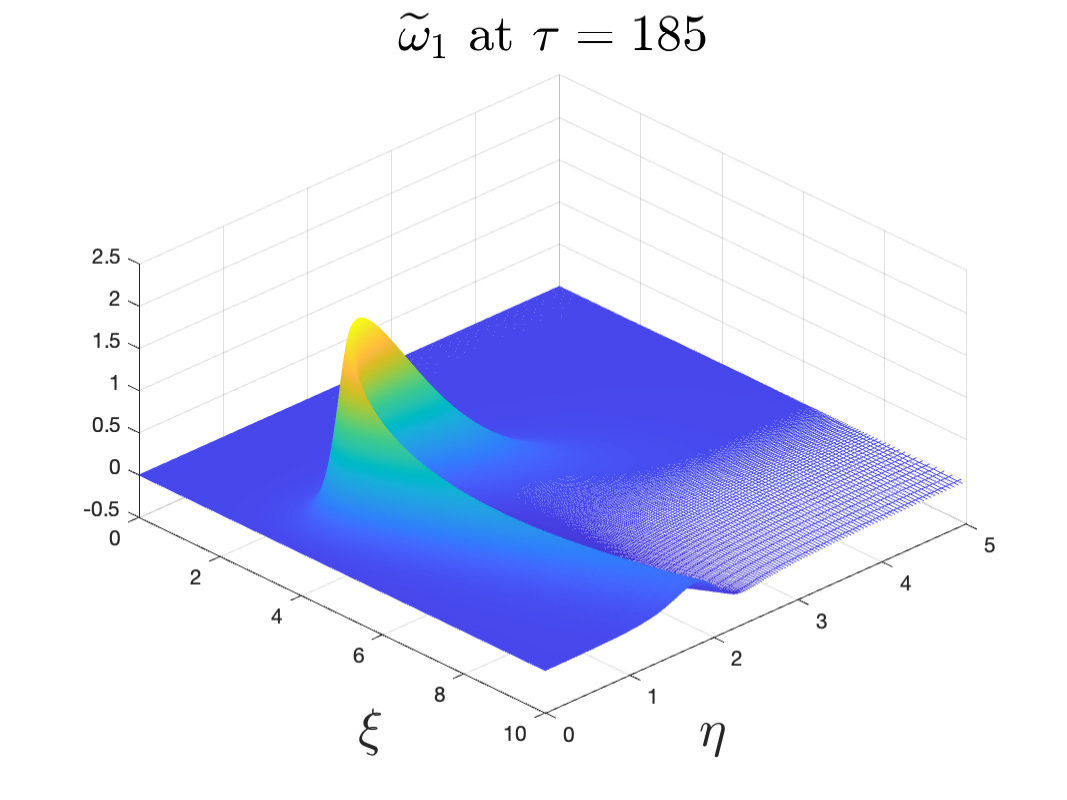}
        \caption{Rescaled profile for $\widetilde{\omega}_1$}
    \end{subfigure}
    \begin{subfigure}[c]{0.40\textwidth}
        \centering
        \includegraphics[width=1\textwidth]{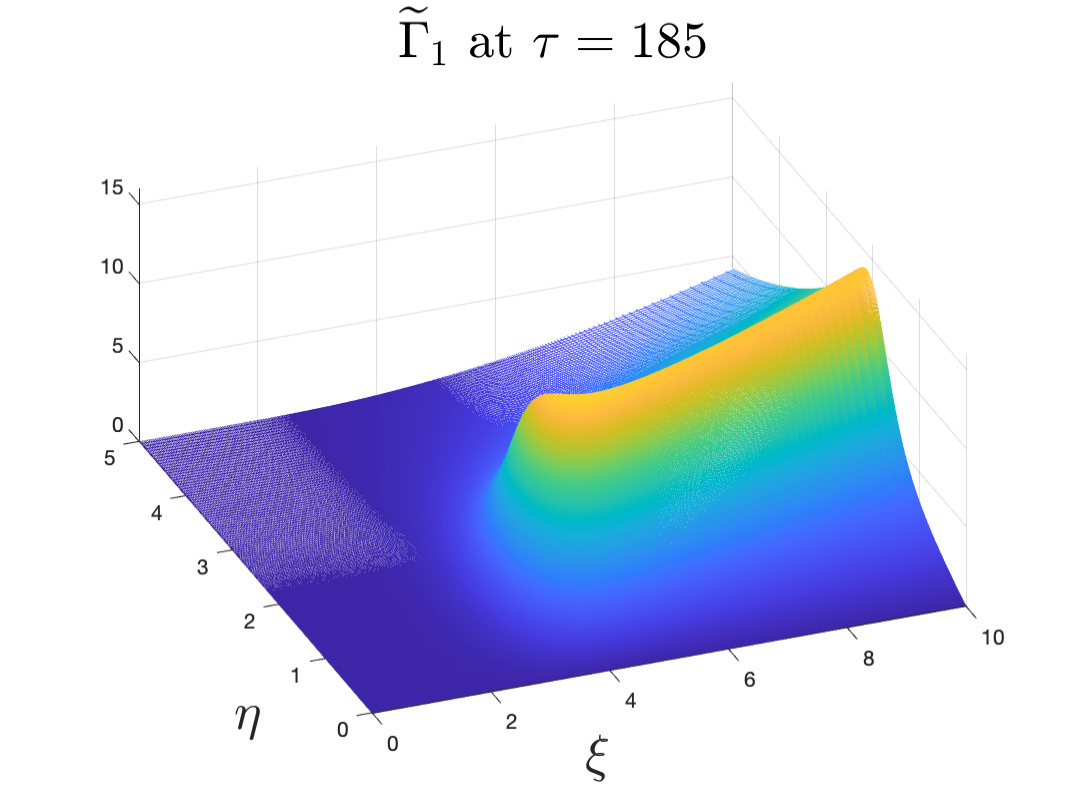}
        \caption{Rescaled profile for $\widetilde{\Gamma}$}
    \end{subfigure}
    \begin{subfigure}[d]{0.40\textwidth}
        \centering
        \includegraphics[width=1\textwidth]{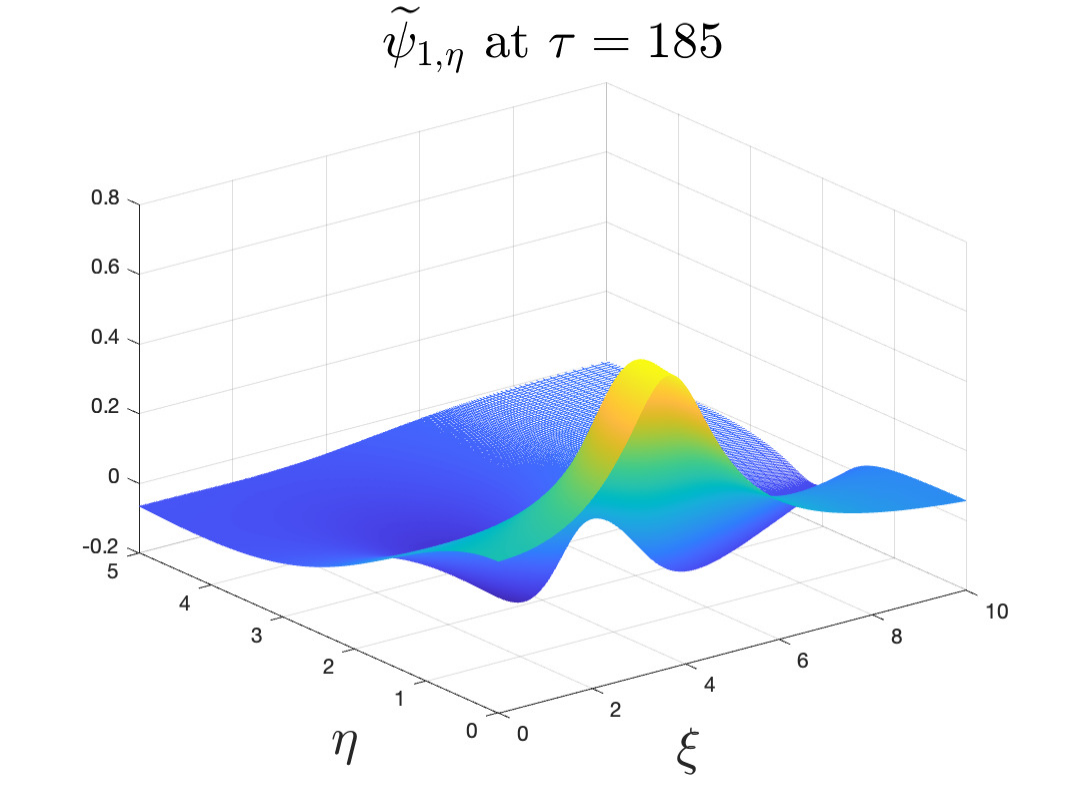}
        \caption{Rescaled profile for $\widetilde{\psi}_{1\eta}$}
    \end{subfigure}
    \caption[Profile]{The local view of rescaled profiles at time $\tau = 185$ for the generalized Navier-Stokes equations. (a) $\widetilde{u}_1$; (b) $\widetilde{\omega}_1$; (c) $\widetilde{\Gamma}$; (d) $\widetilde{\psi}_{1,\eta}$.}  
     \label{fig:profiles-decay}
\end{figure}

In Figure \ref{fig:profiles-decay}, we present the $3$D solution profiles of $(\widetilde{u}_1$, $\widetilde{\om}_1$, $\widetilde{\Gamma}$, $\widetilde{\psi}_{1,\eta}$)
 at time $ \tau =185$. By this time, the maximum vorticity has increased by a factor of $9 \times 10^{21}$ as we can see from Figure \ref{fig:max-vorticity-decay}. We observe that the singular support of the profiles travels toward the origin with distance of order $O(10^{-15})$. Note that the position $(R,Z)$ where $u_1$ achieves its maximum  has been fixed to be at $(R_0,1)$. Due to the viscous regularization, the profile of $\widetilde{u}_1$ remains relatively smooth near $(R,Z)$. Moreover, the thin structure for $\om_1$ that we observed for the $3$D Euler equations in \cite{Hou-euler-2022} becomes much smoother. The tail part of $u_1$ and $\om_1$ is quite smooth and decays rapidly into the far field.
 
Due to the relatively small viscosity $\nu$ for $\widetilde{\Gamma}$, we observe that the total circulation $\widetilde{\Gamma}$ develops a traveling wave with a relatively sharp front, propagating toward the origin.  The diffusion term in the $\widetilde{\omega}_1$ equation regularizes the nearly singular source term due to the sharp traveling wave profile of $\widetilde{\Gamma}$ and generates a regularized Delta function like profile for $\widetilde{\omega}_1$. We observe that $\widetilde{\psi}_{1,\eta}$ achieves its maximum value at $\eta = 0$ near $\xi = R_0$. This property is crucial in generating a traveling wave that propagates toward $\eta = 0$, overcoming the destabilizing effect of the transport along the $\eta$ direction.

We observe that the rescaled profile $\widetilde{\omega}_1$ decays rapidly in the far field with boundary values of $O(10^{-18})$. Similar observation also applies to $u_1$ whose boundary values are of order $O(10^{-12})$. The rescaled profile $\widetilde{\psi}_1$ also has a fast decay in the far field with boundary values of order $O(10^{-6})$. On the other hand, a portion of $\Gamma$ in the near field is transported to the far field, resulting in  $O(1)$ values of $\widetilde{\Gamma}$ in the far field.  We note that $\widetilde{\Gamma}$ contributes to the generalized Navier--Stokes equations only through the vortex stretching term $(\widetilde{\Gamma}^2/\xi^4)_\eta$, which is extremely small with order $O(10^{-34})$ at the far field boundary. Therefore, we do not need to enforce the decay of $\widetilde{\Gamma}$ in the far field. Moreover,  the boundary values of the transport terms for $\widetilde{\Gamma}$ and $\widetilde{\omega}_1$ are of order $10^{-9}$ and $10^{-29}$, respectively. 


\subsection{The streamlines} 

\begin{figure}[!ht]
\centering
    \begin{subfigure}[b]{0.40\textwidth}
        \centering
        \includegraphics[width=1\textwidth]{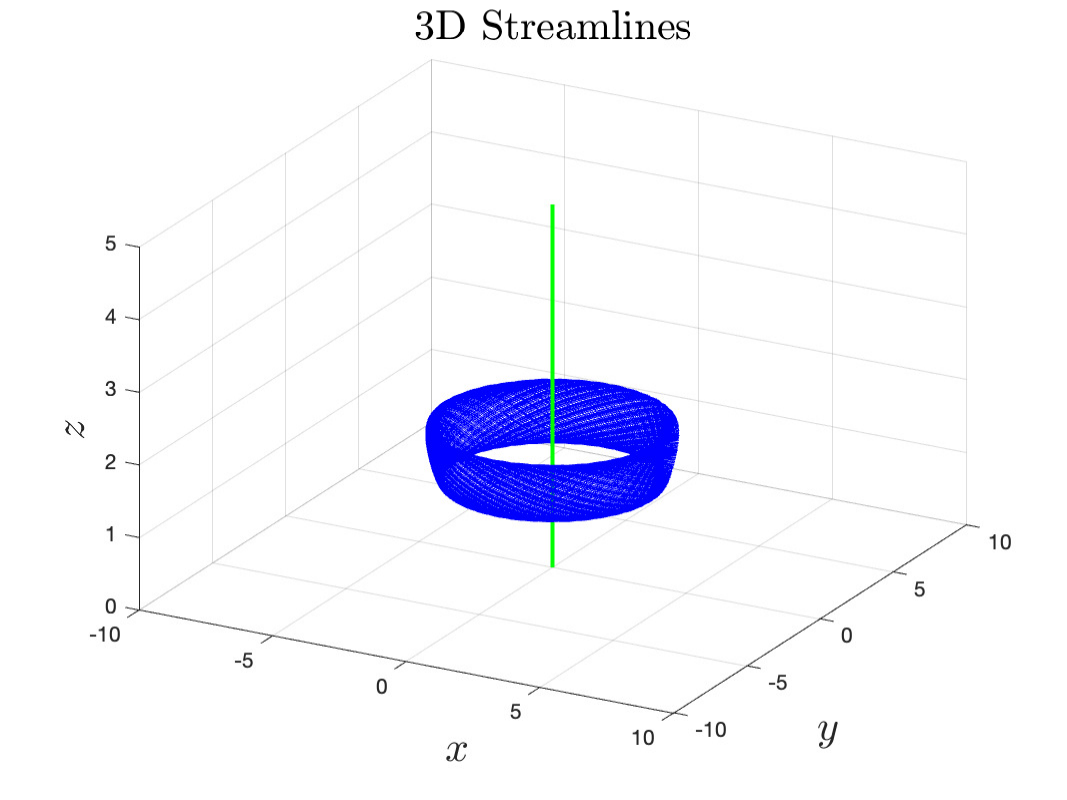}
        \caption{$r_0 = 4$, $z_0 = 2$}
    \end{subfigure}
    \begin{subfigure}[b]{0.40\textwidth}
        \centering
        \includegraphics[width=1\textwidth]{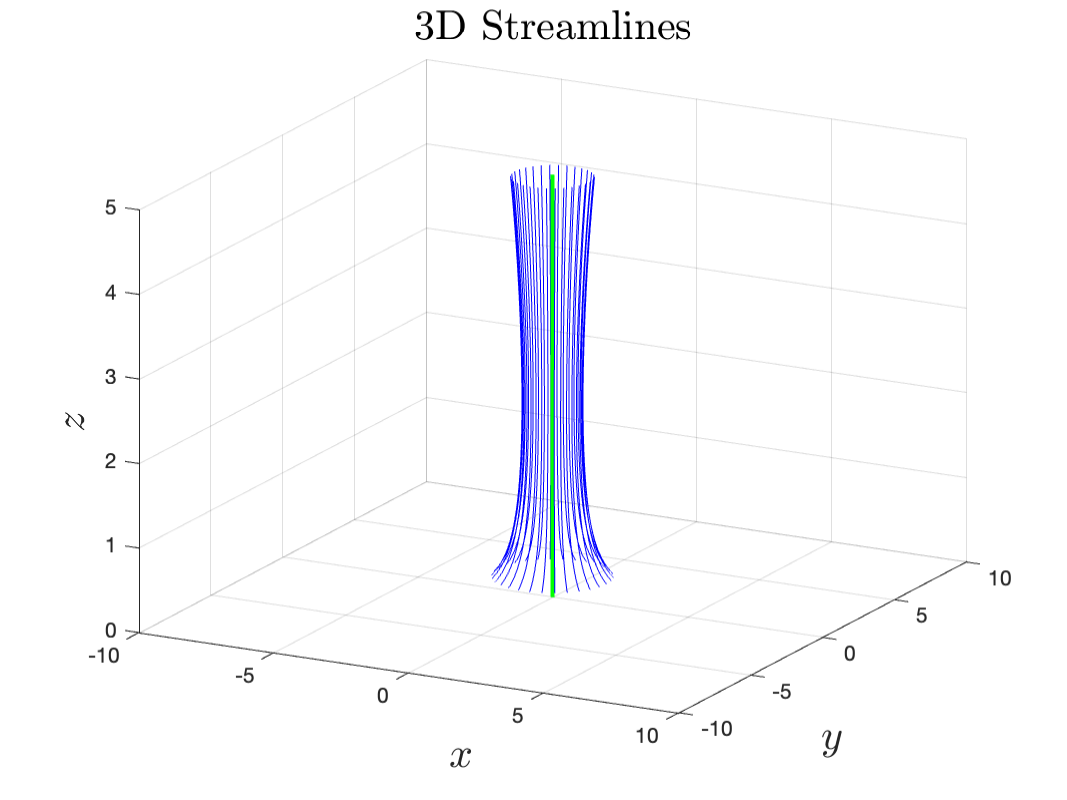}
        \caption{$r_0 = 2$, $z_0 = 0.25$}
    \end{subfigure}
    \begin{subfigure}[c]{0.40\textwidth}
        \centering
        \includegraphics[width=1\textwidth]{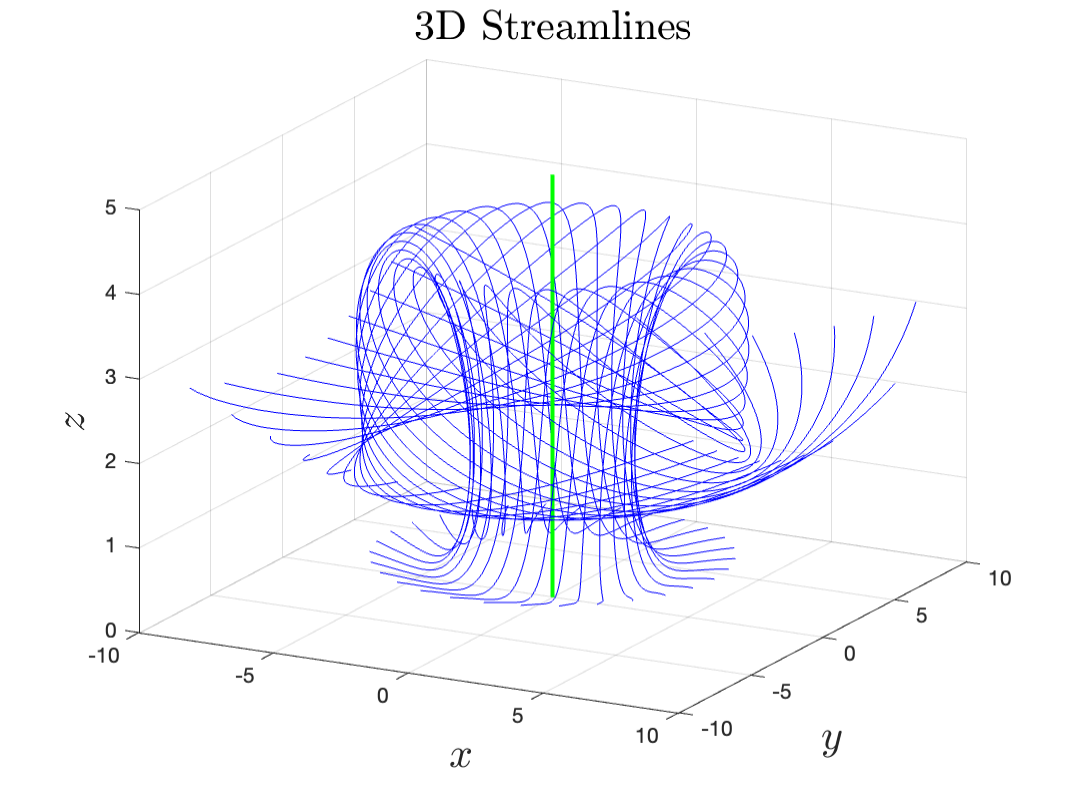}
        \caption{$r_0 = 6$, $z_0 = 0.5$}
    \end{subfigure}
    \begin{subfigure}[d]{0.40\textwidth}
        \centering
        \includegraphics[width=1\textwidth]{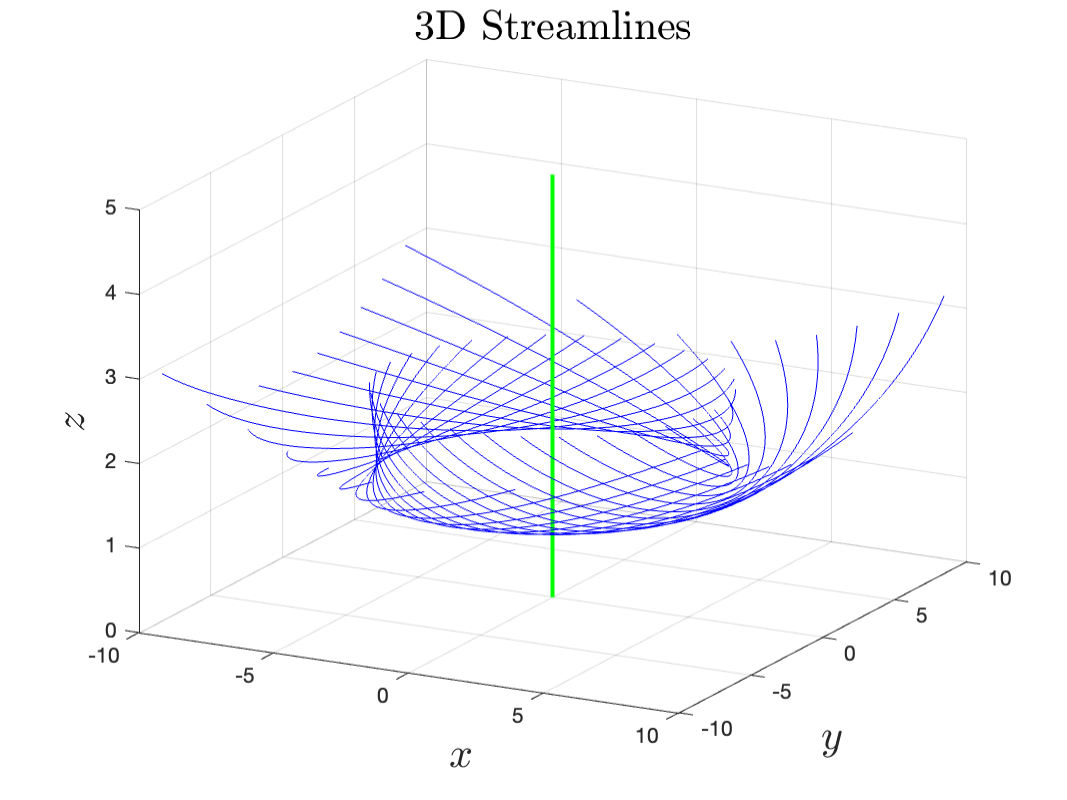}
        \caption{$r_0 = 6$, $z_0 = 2.5$}
    \end{subfigure}
    \caption[Global streamline]{The streamlines of $(u^r(t),u^\theta(t),u^z(t))$ for the generalized Navier-Stokes equations at time $\tau = 185$ with initial points given by (a) $(r_0,z_0) = (4,2)$, streamlines form a torus; (b) $(r_0,z_0) = (2,0.25)$, streamlines go straight upward; (c) $(r_0,z_0) = (6,0.5)$, streamlines first go downward, then travel inward, finally go upward; (d) $(r_0,z_0) = (6,2.5)$, streamlines spin downward and outward. The green pole is the symmetry axis $r=0$.}  
     \label{fig:streamline-nse-decay}
\end{figure}

In this subsection, we investigate the features of the velocity field. We first study the velocity field by looking at the induced streamlines. 
Interestingly the induced streamlines look qualitatively the same as those obtained for the $3$D Navier--Stokes equations in a periodic cylinder \cite{Hou-nse-2022}. In Figure~\ref{fig:streamline-nse-decay}, we plot the streamlines induced by the velocity field ${\bf \widetilde{u}}$ at $\tau = 185$ for different initial points. By this time, the ratio between the maximum vorticity and the initial maximum vorticity, i.e. $\|\vom (t)\|_{L^\infty}/\|\vom (0)\|_{L^\infty}$, has increased by a factor of $9\times 10^{21}$. 

The velocity field resembles that of a tornado spinning around the symmetry axis (the green pole). In  Figure~\ref{fig:streamline-nse-decay}(a) with $(r_0,z_0) = (4,1.5)$, we observe that the streamlines form a torus spinning around the symmetry axis. In  Figure~\ref{fig:streamline-nse-decay}(b) with $(r_0,z_0) = (2,0.25)$, the streamlines go straight upward without any spinning. In  Figure~\ref{fig:streamline-nse-decay}(c) with $(r_0,z_0) = (6,2.5)$, the streamlines first go downward, then travel inward and upward, finally travel downward and spin outward. In Figure~\ref{fig:streamline-nse-decay}(d) with $(r_0,z_0) = (6,2)$, the streamlines first spin downward and then outward. The solution behaves qualitatively the same as what we observed for the axisymmetric Navier--Stokes equations in a periodic cylinder \cite{Hou-nse-2022}. 

\subsection{The 2D flow} 
\begin{figure}[!ht]
\centering
    \includegraphics[width=0.40\textwidth]{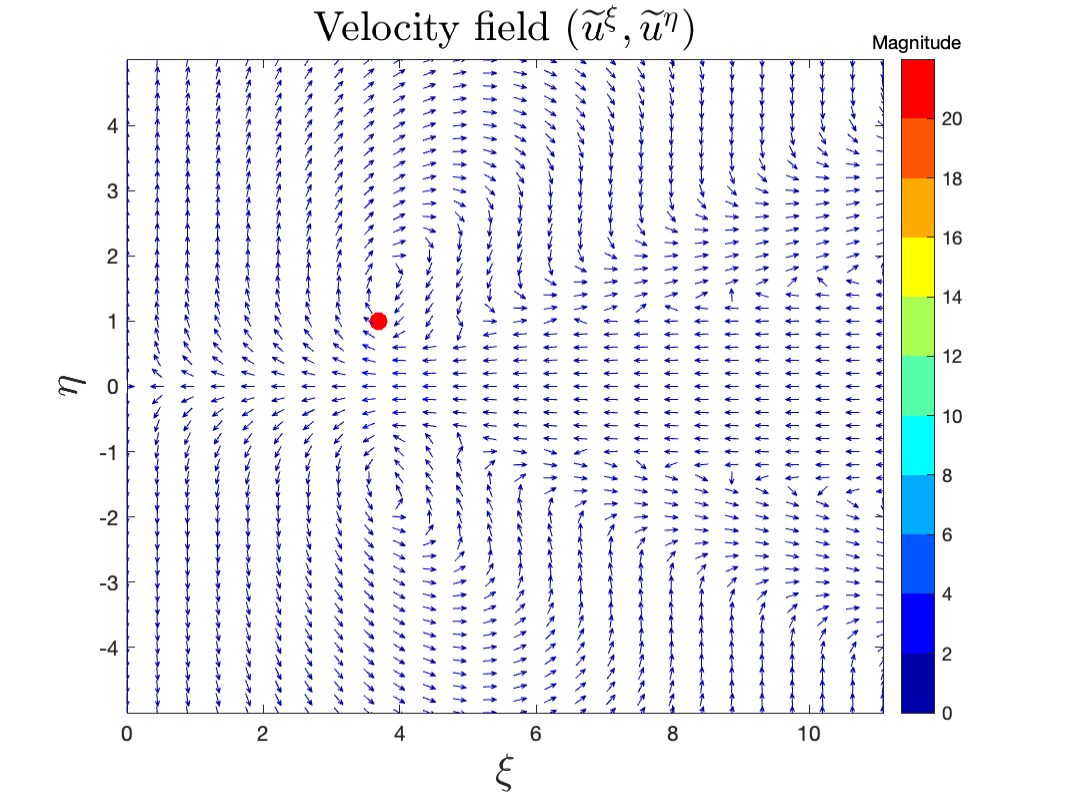}
    \includegraphics[width=0.40\textwidth]{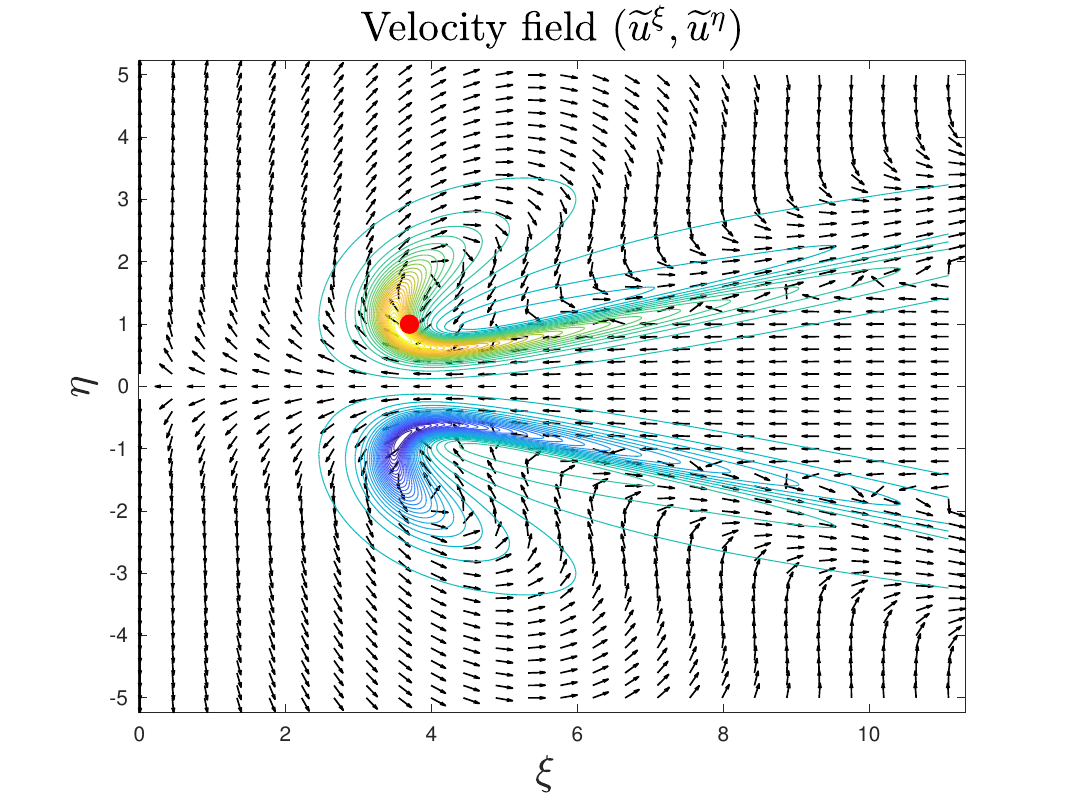} 
    \caption[Dipole]{The dipole structure of $\widetilde{\omega}_1$  and the induced dynamically rescaled velocity field $(\widetilde{u}^\xi,\widetilde{u}^\eta)$ for the generalized Navier-Stokes equations at $\tau=185$. Recall $\xi = r/C_{lr}$ and $\eta = z/C_{lz}$. The arrows in both plots are used to illustrate the direction of the velocity field. They are normalized to have the same length for better plotting quality and do not have any physical meaning here. Left plot: the velocity vector. Right plot: the velocity vector with the $\omega_1$ contour as background. The red dot is the position $(R,Z)$ where $\widetilde{u}_1$ achieves its maximum.}  
     \label{fig:dipole_nse-decay}
\end{figure}

To understand the phenomena in the most singular region as shown in Figure~\ref{fig:streamline-nse-decay}, we study the $2$D velocity field $(u^r,u^z)$. In Figure~\ref{fig:dipole_nse-decay}(a)-(b), we plot the dipole structure of $\widetilde{\omega}_1$ in a local symmetric region and the hyperbolic velocity field induced by the dipole structure in a local microscopic domain $[0,R_b]\times [0,Z_b]$ at  $\tau=185$.
The dipole structure for the generalized Navier--Stokes equations look qualitatively similar to that of the $3$D Navier--Stokes equations in a periodic cylinder \cite{Hou-nse-2022}. As in the case of the $3$D Navier--Stokes equations in a periodic cylinder, the negative radial velocity near $\eta=0$ induced by the antisymmetric vortex dipoles pushes the solution toward $\xi=0$, then move upward away from $\eta=0$. This is one of the driving mechanisms for a potential singularity on the symmetry axis. Since the value of $\widetilde{u}_1$ becomes very small near the symmetry axis $\xi=0$, the streamlines almost do not spin around the symmetry axis, as illustrated in Figure~\ref{fig:streamline-nse-decay}(b).  

We also observe that the velocity field $(u^r(t),u^z(t))$ forms a closed circle right above $(R,Z)$. The corresponding streamlines are trapped in the circle region in the $\xi \eta$-plane, which is responsible for the formation of the spinning torus that we observed earlier in Figure~\ref{fig:streamline-nse-decay}(a).

\begin{figure}[!ht]
\centering
    \includegraphics[width=0.40\textwidth]{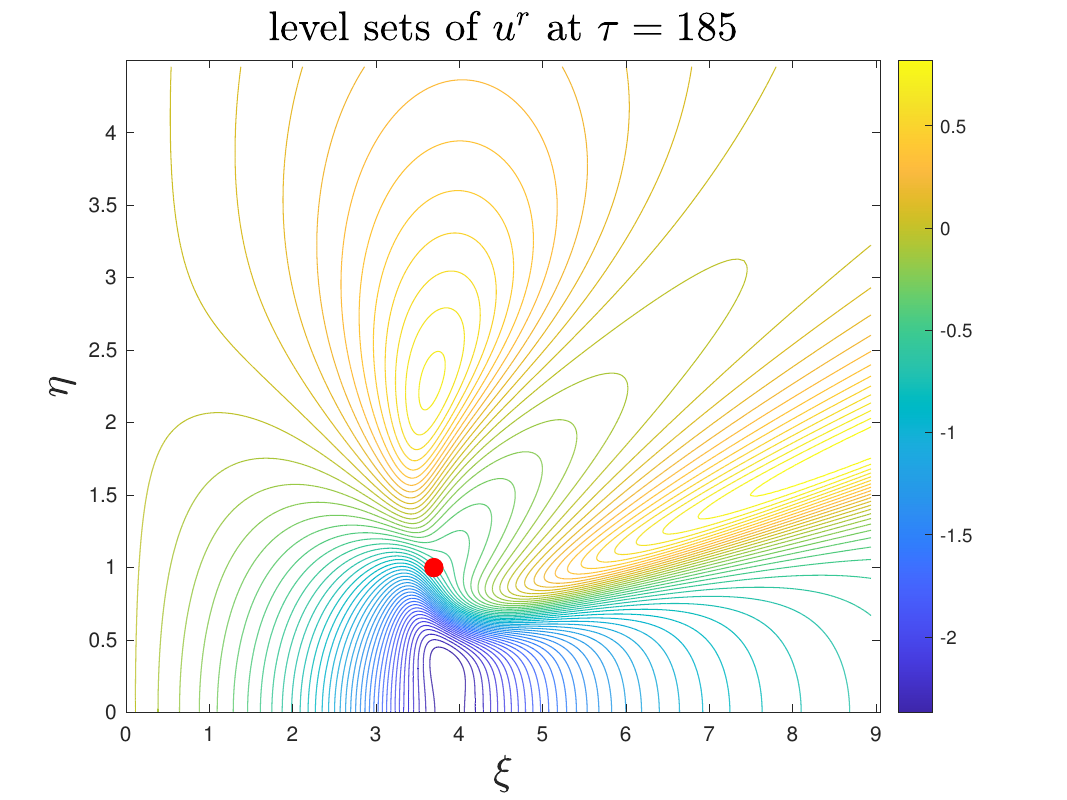}
    \includegraphics[width=0.40\textwidth]{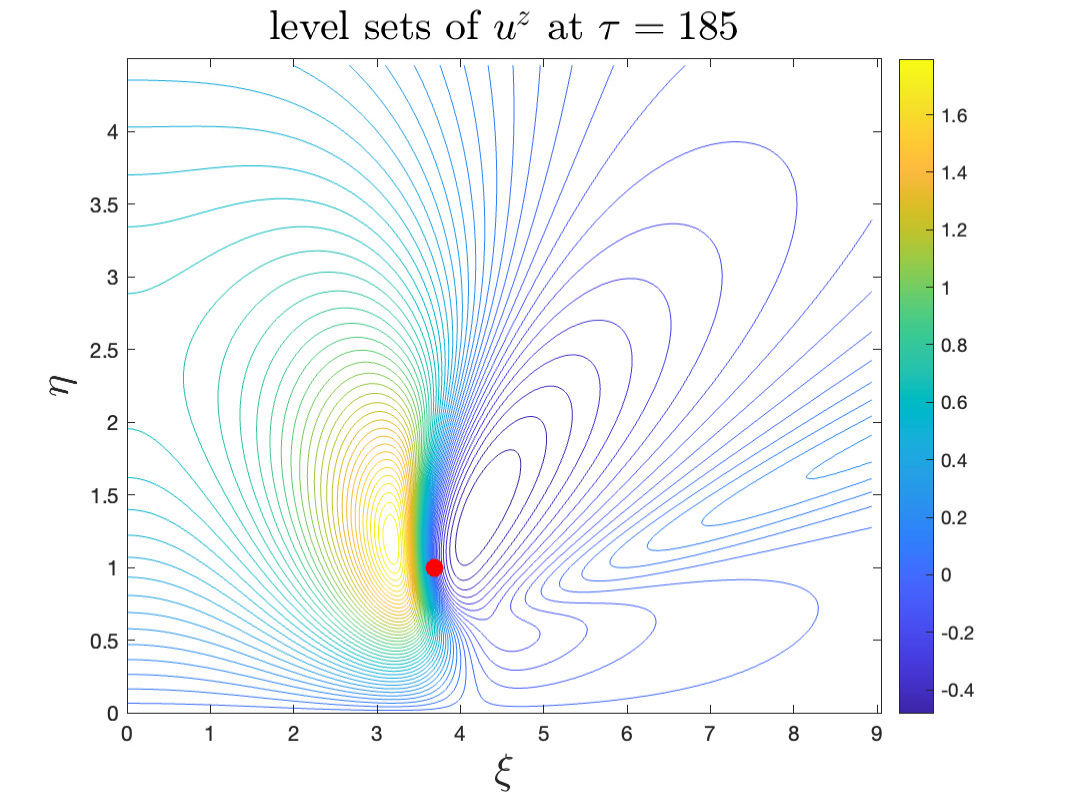} 
    \caption[Velocity level sets]{The level sets of $\tilde{u}^\xi$ (left) and $\tilde{u}^\eta$ (right) for the generalized Navier-Stokes equations at $\tau = 185$. The red point is the maximum location $(R,Z)$ of $\widetilde{u}_1$.}  
       \label{fig:velocity_levelset_euler}
\end{figure}

We can also understand this hyperbolic flow structure from the velocity contours in Figure~\ref{fig:velocity_levelset_euler} (a)-(b). As we can see from Figure~\ref{fig:velocity_levelset_euler}(a), the radial velocity $u^r$ is negative and large in amplitude below the red dot $(R,Z)$ where 
$\widetilde{u}_1$ achieves its maximum, pushing the flow toward the symmetry axis $\xi=0$. But it becomes large and positive above $(R,Z)$, pushing the flow outward. Similarly, we can see from Figure \ref{fig:velocity_levelset_euler}(b) that the axial velocity $u^z$ is negative and large in amplitude to the right hand side of $(R,Z)$, pushing the flow downward toward $\eta = 0$. But it becomes large and positive on the left hand side of $(R,Z)$, pushing the flow upward away from $\eta = 0$. This is the driving mechanism for forming the hyperbolic flow structure near the origin.

\begin{figure}[!ht]
\centering
    \includegraphics[width=0.4\textwidth]{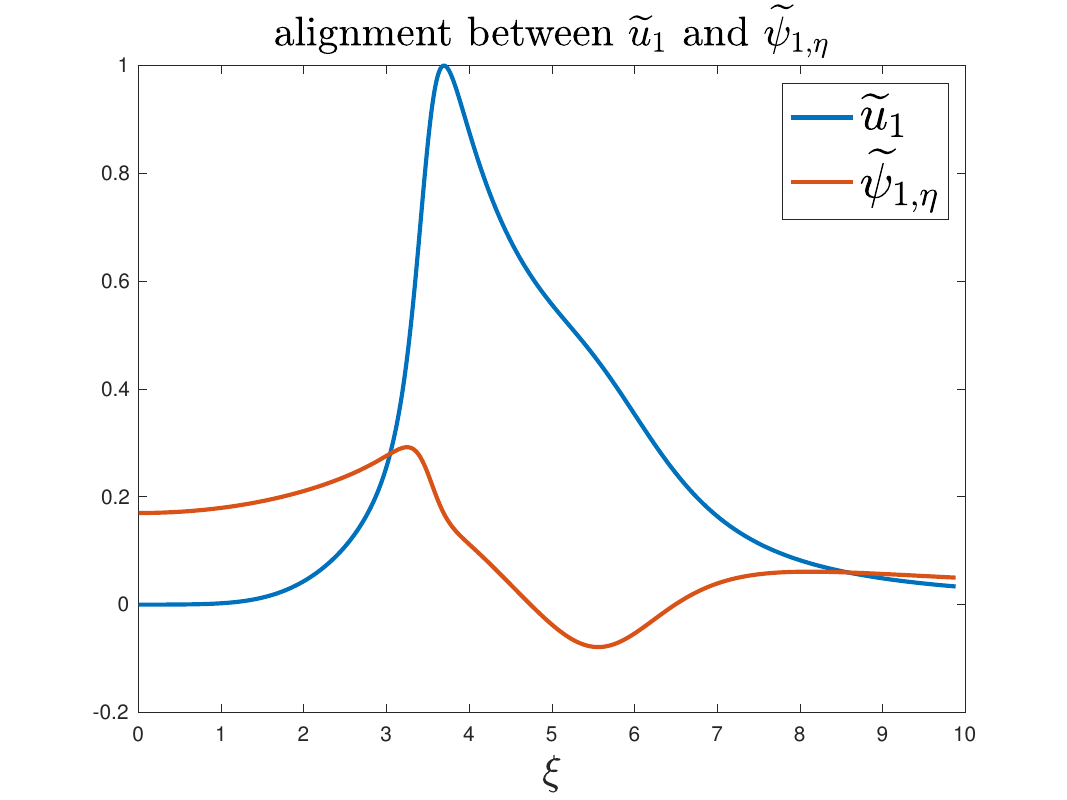}
    \includegraphics[width=0.4\textwidth]{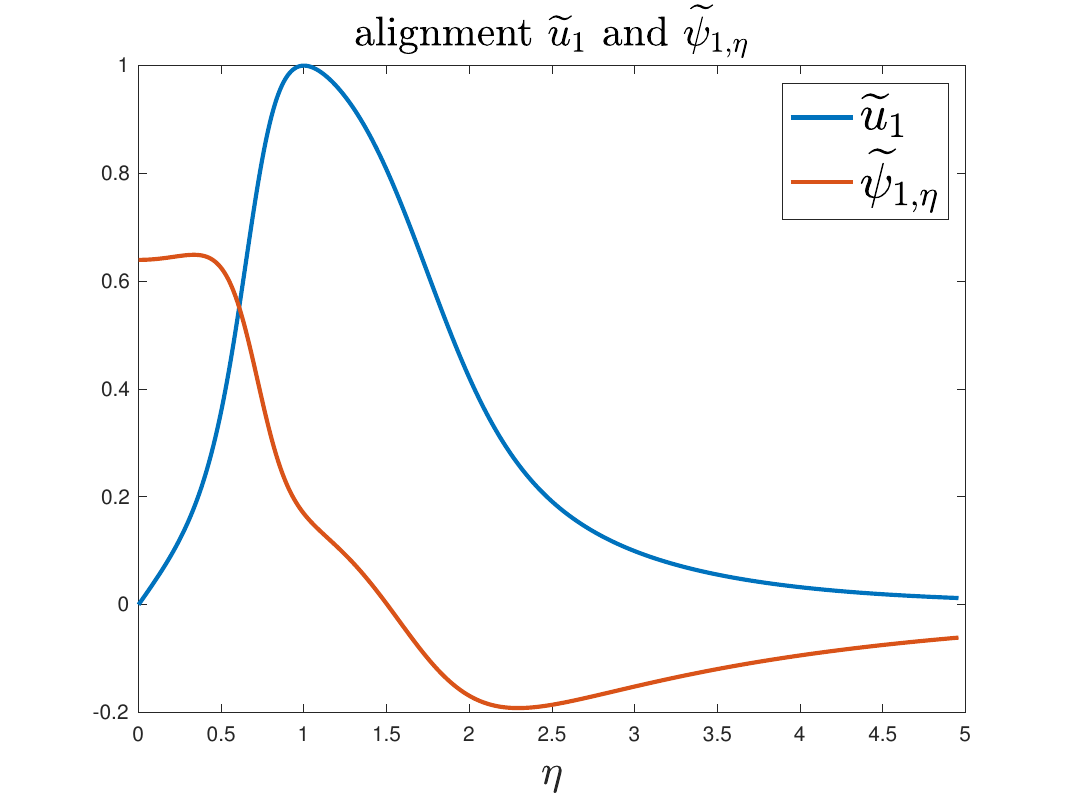} 
    \caption[Alignment-space]{ Left plot: Alignment between alignment $\widetilde{u}_1$ and $\widetilde{\psi}_{1,\eta}$ at $\eta = Z$ as a function of $\xi$ for the generalized Navier-Stokes equations at $\tau = 185$. Right plot:  Alignment between alignment $\widetilde{u}_1$ and $\widetilde{\psi}_{1,\eta}$ at $\xi = R$ as a function of $\eta$ at $\tau = 185$.} 
    \label{fig:alignment-space-decay}
\end{figure}

%
%


In Figure \ref{fig:alignment-space-decay}(a)-(b), we demonstrate the alignment between $\widetilde{\psi}_{1\eta}$ and $\widetilde{u}_1$ at $\tau =185$. Although the maximum vorticity has grown a lot by this time, the local solution structures have remained qualitatively the same in the late stage of the computation.  This shows that the viscous effect has a strong stabilizing effect that enhances the nonlinear alignment of vortex stretching. We also observe that $\widetilde{\psi}_{1\eta}$ is relatively flat in the region $\{(\xi,\eta) | 0 \leq \xi \leq 0.9R, \; 0 \leq 
\eta \leq 0.5Z\}$. This property is critical for $\widetilde{u}_1$ to remain large between the sharp front and $\xi=0$, thus avoiding developing a two-scale structure.


\subsection{Alignment of vortex stretching}

\begin{figure}[!ht]
\centering
    \includegraphics[width=0.4\textwidth]{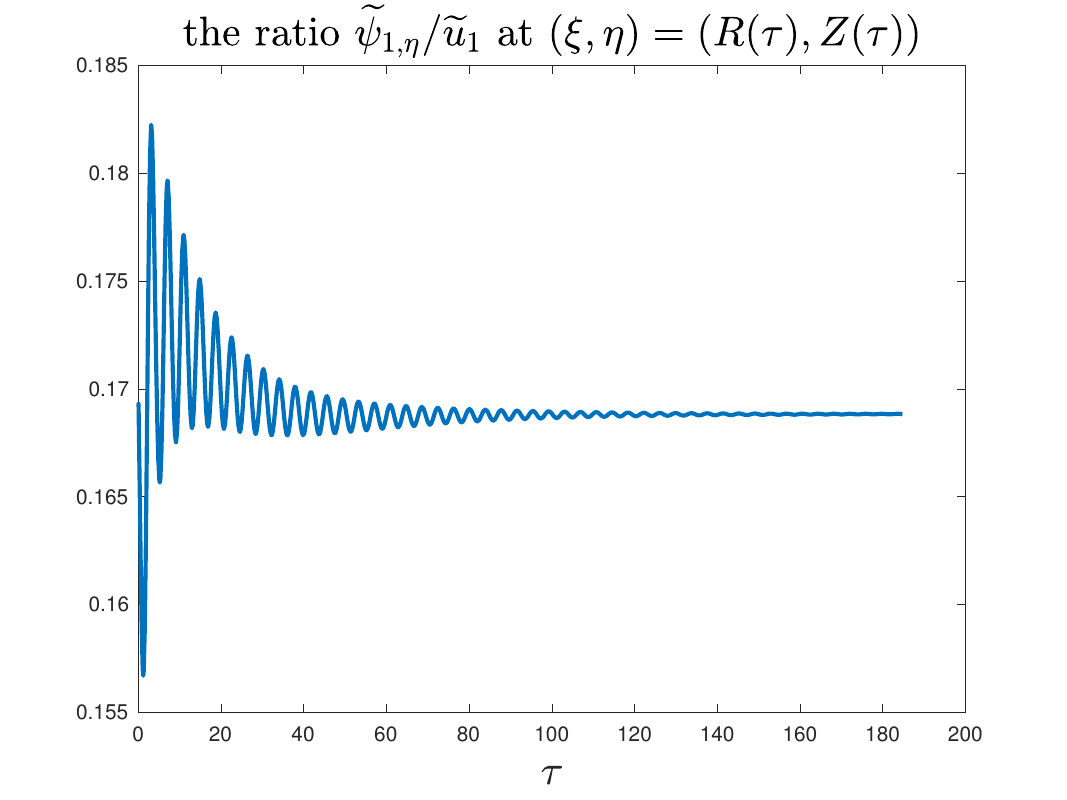}
 \includegraphics[width=0.4\textwidth]{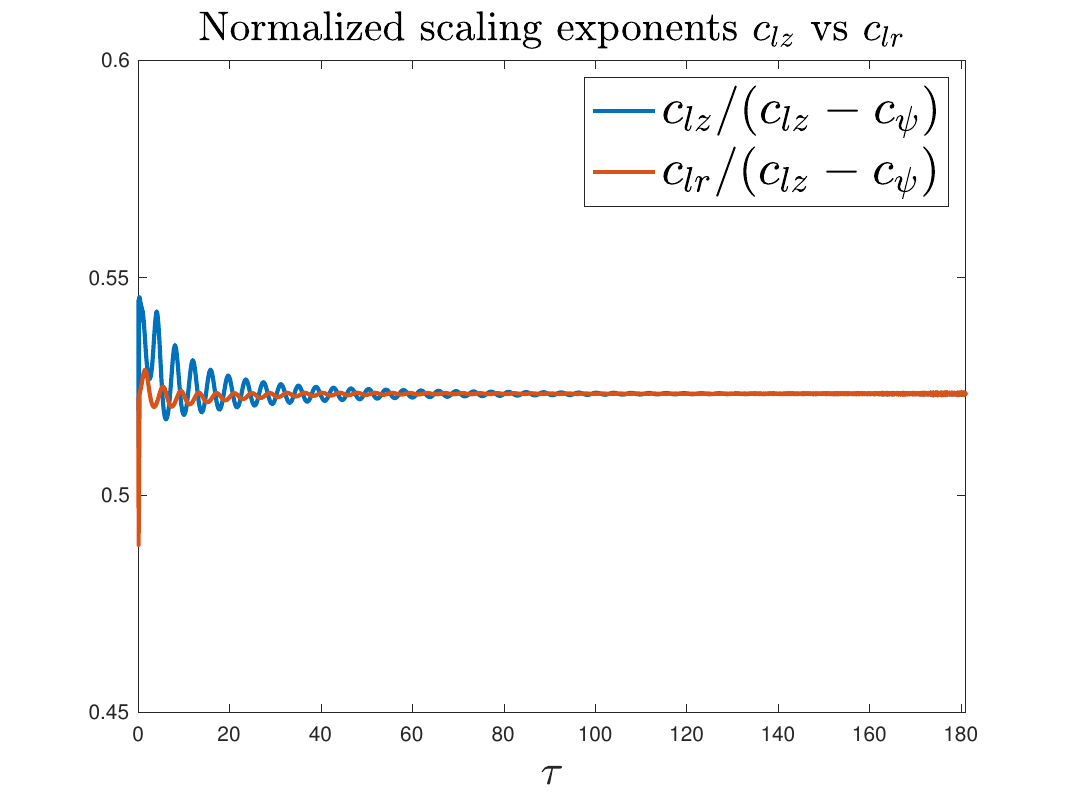} 
    \caption[Alignment-time]{ Left plot: The ratio between $\widetilde{\psi}_{1,\eta}$ and $\widetilde{u}_1$ for the generalized Navier-Stokes equations at $(\xi,\eta)=(R,Z)$ as a function of $\tau$. Right plot: The normalized scaling exponent $c_{lz}/(c_{lz}-c_\psi)$ and $c_{lr}/(c_{lz}-c_\psi)$ as a function of $\tau$.} 
    \label{fig:alignment-time-2-decay}
\end{figure}

Due to the viscous regularization, the solution becomes smoother and is more stable. We are able to compute up to a time when $(R(t),Z(t))$ is very close to the origin. This is something we could not achieve for the $3$D Navier--Stokes equations in a periodic cylinder \cite{Hou-nse-2022}. In Figure \ref{fig:alignment-time-2-decay}(a),  we observe that the positive alignment between $\widetilde{u}_1$ and $\widetilde{\psi}_{1
\eta}$ converges to a constant as $\tau$ increases. This indicates that the generalized Navier--Stokes equations with solution dependent viscosity achieves a self-similar scaling relationship. We observe some mild oscillations in time for the alignment and the normalized coefficients in Figure \ref{fig:alignment-time-2-decay}. This is due to the rapid decay of the solution dependent viscosity $\nu(t)$ in time (see Figure \ref{fig:ratio-vort-diff-w1-decay}(b)), which is not strong enough to stabilize the shearing instability induced by the sharp front of $\Gamma$ in the generalized Euler equations.

In the Figure \ref{fig:alignment-time-2-decay}(b), we plot the normalized scaling parameters 
$\widehat{c}_{lz} = c_{lz}/(c_{lz}-c_\psi)$ and $\widehat{c}_{lr} = c_{lr}/(c_{lz}-c_\psi)$. We observe that they converge to the same constant $c_l = 0.523$ as $\tau$ increases.  This provides some evidence that the dynamically rescaled solution seems to converge to an approximate steady state as $\tau \rightarrow \infty$.

If we kept the space dimension $n=3$ for all time, we observed that $c_{lz}$ and $c_{lr}$ did not converge to the same value as $\tau$ increases. It would generate a potential two-scale blowup, which is not compatible with the scaling properties of a Navier--Stokes blowup. However, when we vary the space dimension as $n(\tau) = 1+2R(\tau)/Z(\tau)$, we observe that $c_{lz}$ and $c_{lr}$ magically converge to the same value as $\tau$ increases. This shows that using the space dimension as an extra degree of freedom can eliminate the scaling instability and generate a one-scale self-similar blowup. Since $\|u_1\|_\infty= 1/(T-t)$ and $Z(t) = (T-t)^{c_l}$, the solution dependent viscosity is given by $\nu= \nu_0 \|u_1(\tau)\|_\infty Z(t)^2 = \nu_0 (T-t)^{2c_l-1} = \nu_0 (T-t)^{0.1272}$.

\subsection{Balance between vortex stretching and diffusion and self-similar profile}

In this subsection, we study the balance between the vortex stretching terms and the diffusion terms for both $\widetilde{u}_1$ and $\widetilde{\omega}_1$ equations. 
Using the scaling relationship $c_u = c_\psi - c_{lz}$ given by \eqref{scaling-relation-cu}, we can easily show that the solution dependent viscosity satisfies
 \[
 \nu (\tau) = \nu_0 \|u_1\|_{\infty} Z(t)^2 = \nu_0 C_{lz}/C_\psi ,
 \] 
which exactly cancels the scaling factor $C_\psi/C_{lz}$ in front of the rescaled diffusion term. This is why we obtain constant viscosity $\nu_0$ in the dynamic rescaling formulation for both $\widetilde{u}_1$ and $\widetilde{\omega}_1$ equations.

\begin{figure}[!ht]
\centering
    \includegraphics[width=0.4\textwidth]{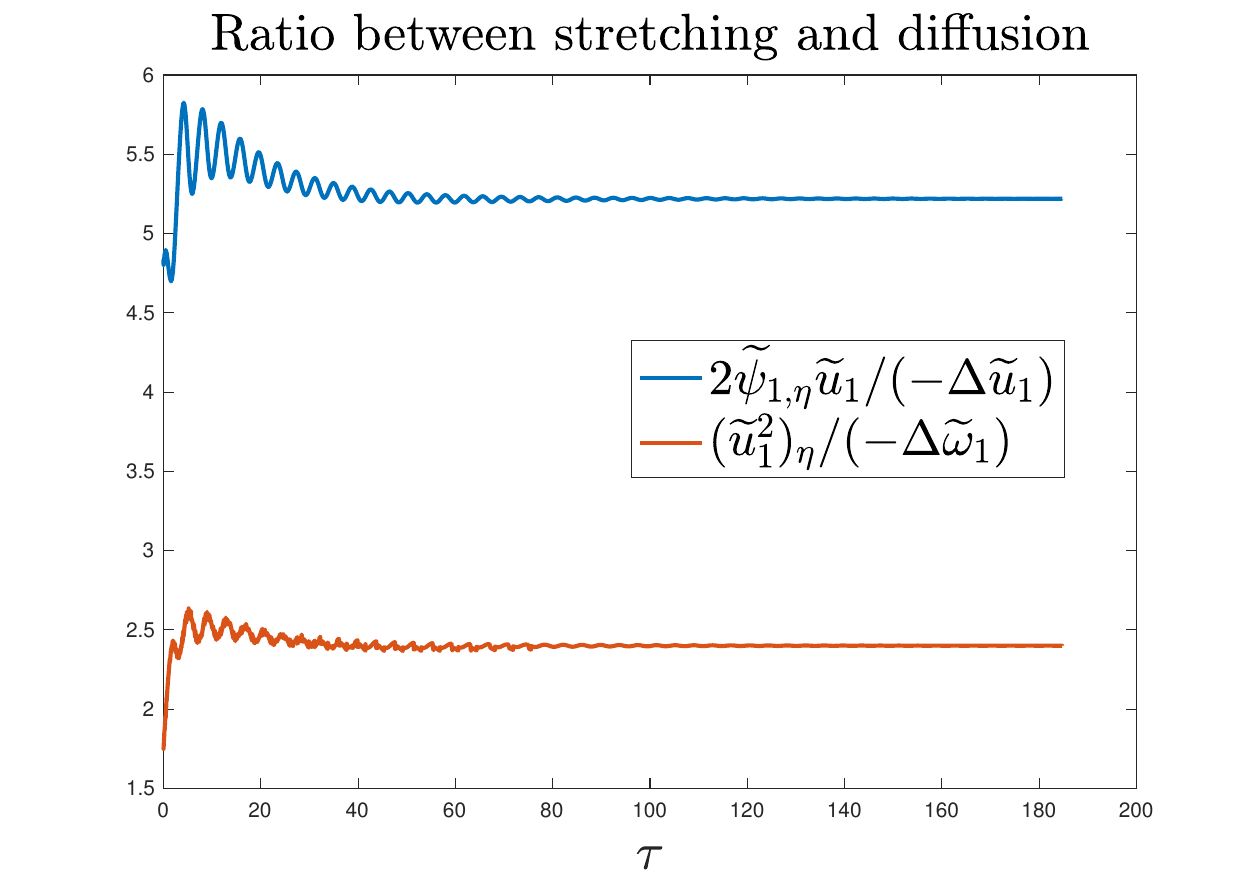}
 \includegraphics[width=0.4\textwidth]{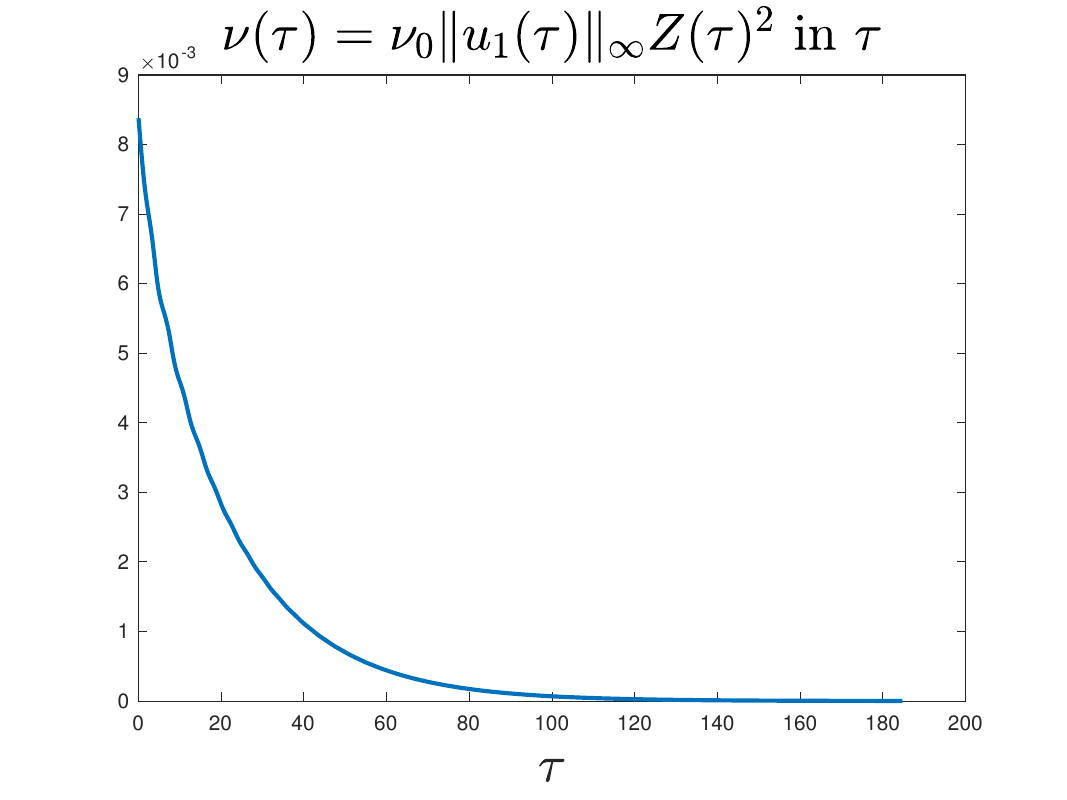} 
    \caption[VelL3-time]{ Left plot: Ratio between the vortex stretching term $2\widetilde{\psi}_{1,\eta}\widetilde{u}_1$ and the diffusion term $-\nu_0\Delta \widetilde{u}_1$ at $(R,Z)$ (blue) and ratio between the vortex stretching term $(\widetilde{u})^2_{\eta}$ and the diffusion term $-\nu_0\Delta \widetilde{\omega}_1$ at $(R_\omega,Z_\omega)$ (red) in $\tau$. Right plot: Viscosity $\nu(\tau) = \nu_0 \|u_1(\tau)\|_\infty Z(\tau)^2$ in $\tau$ with $\nu(185) = 1.3 \cdot 10^{-6}$.} 
    \label{fig:ratio-vort-diff-w1-decay}
\end{figure}

In Figure \ref{fig:ratio-vort-diff-w1-decay}(a), we plot the ratio between the vortex stretching term $2 \widetilde{\psi}_{1,\eta} \widetilde{u}_1$ 
and the diffusion term $-\nu_0\Delta \widetilde{u}_1$
 at $(R,Z)$ where $\widetilde{u}_1$ achieves its maximum. We observe that this ratio converges to a constant value $5.22$ as $\tau$ increases. This shows that the vortex stretching term dominates the diffusion term. 
In the same figure, we also plot the ratio between $(\widetilde{u}_1^2)_\eta$ and 
$-\nu_0\Delta \widetilde{\omega}_1$
 at $(R_\omega,Z_\omega)$ where $\widetilde{\omega}_1$ achieves its maximum. We do not include the contribution from the term $(n-3)\widetilde{\psi}_{1,\eta}\widetilde{\omega}_1$ since $(n-3)\widetilde{\psi}_{1,\eta}\widetilde{\omega}_1$ is only $7\%$ of $(\widetilde{u}_1^2)_\eta$ at $(R_\omega,Z_\omega)$. We observe that the ratio of these two terms converges to a constant value $2.4$ as $\tau$ increases, which again shows that the vortex stretching term dominates the diffusion, but the diffusion term has a nontrivial contribution as the solution develops a self-similar blowup. The balance between the vortex stretching terms and the diffusion terms is crucial in maintaining the robust nonlinear growth of the maximum vorticity in time. 
 
In Figure \ref{fig:ratio-vort-diff-w1-decay}(b), we plot the solution dependent viscosity $\nu(\tau) = \nu_0 \|u_1\|_\infty Z(\tau)^2$ as a function of $\tau$.  We observe that this solution dependent viscosity decays to zero as $\tau$ increases with $\nu(185) = 1.3 \cdot 10^{-6}$.  However, the solution dependent viscosity still has an $O(1)$ effect on the self-similar blowup profile since the limiting rescaled profile does not satisfy the corresponding Euler equations as the vanishing viscosity limit. Instead, it satisfies the Navier-Stokes equations with constant viscosity $\nu_0$.

\begin{figure}[!ht]
\centering
    \includegraphics[width=0.4\textwidth]{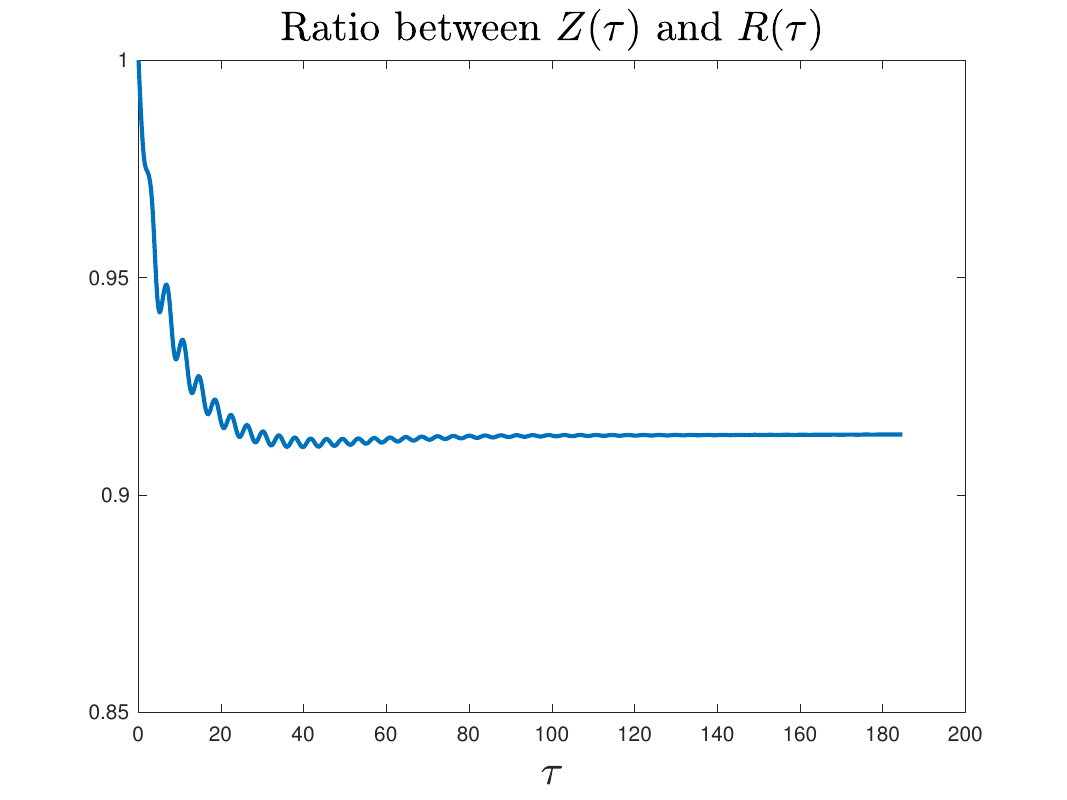}
 \includegraphics[width=0.4\textwidth]{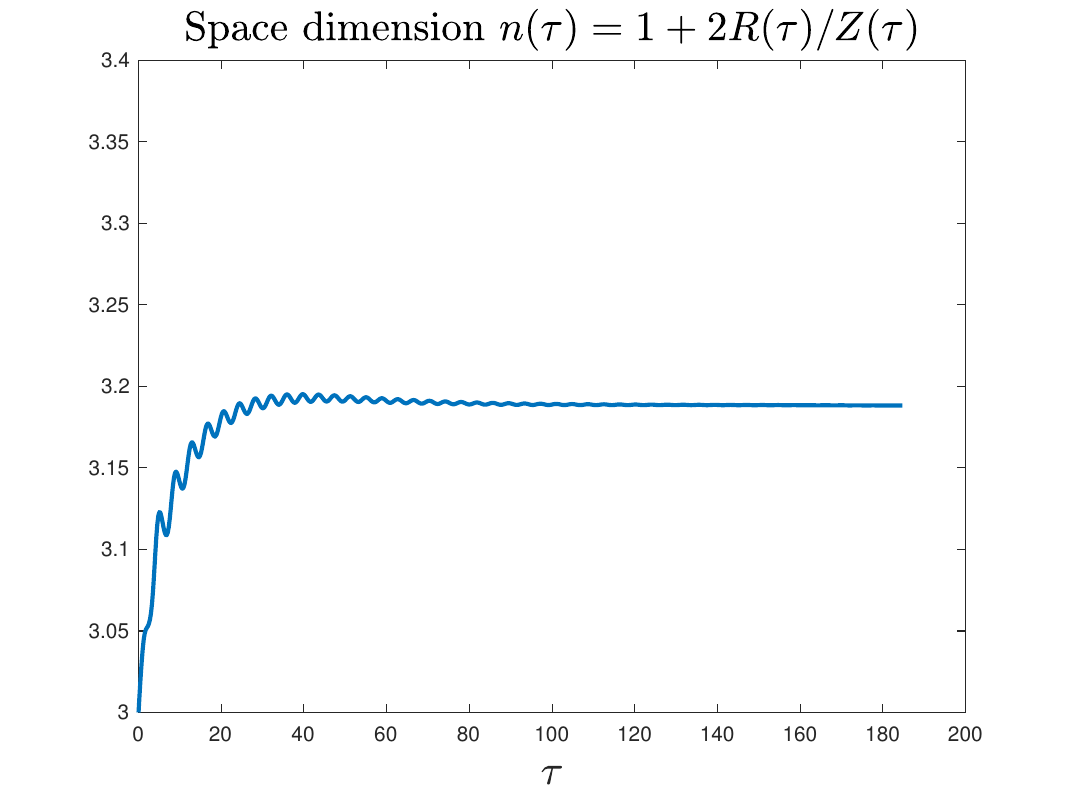} 
    \caption[VelL3-time]{ Left plot: Ratio of $Z(\tau)/R(\tau)$  as a function of $\tau$ for the generalized Navier-Stokes equations. Right plot: The space dimension $n(\tau) = 1 +2R(\tau)/Z(\tau)$ as a function of $\tau$ with $n(185) = 3.188$.} 
    \label{fig:ratio-Z-R-decay}
\end{figure}

%

In Figure \ref{fig:ratio-Z-R-decay}(a),
we plot the ratio between $Z(\tau)$ and $R(\tau)$. We observe that the ratio $Z(\tau)/R(\tau)$ converges to a constant value of $0.914$. This shows that we have a one-scale blowup. In Figure \ref{fig:ratio-Z-R-decay}(b), we plot the space dimension $n(\tau) = 1+2R(\tau)/Z(\tau)$ as a function of $\tau$. We observe that the space dimension is relatively flat and seems to converge to a constant value $n=3.188$ by $\tau =185$.  This provides further evidence that the dynamically rescaled solution seems to converge to an approximate steady state as $\tau \rightarrow \infty$.

\begin{figure}[!ht]
\centering
    \includegraphics[width=0.4\textwidth]{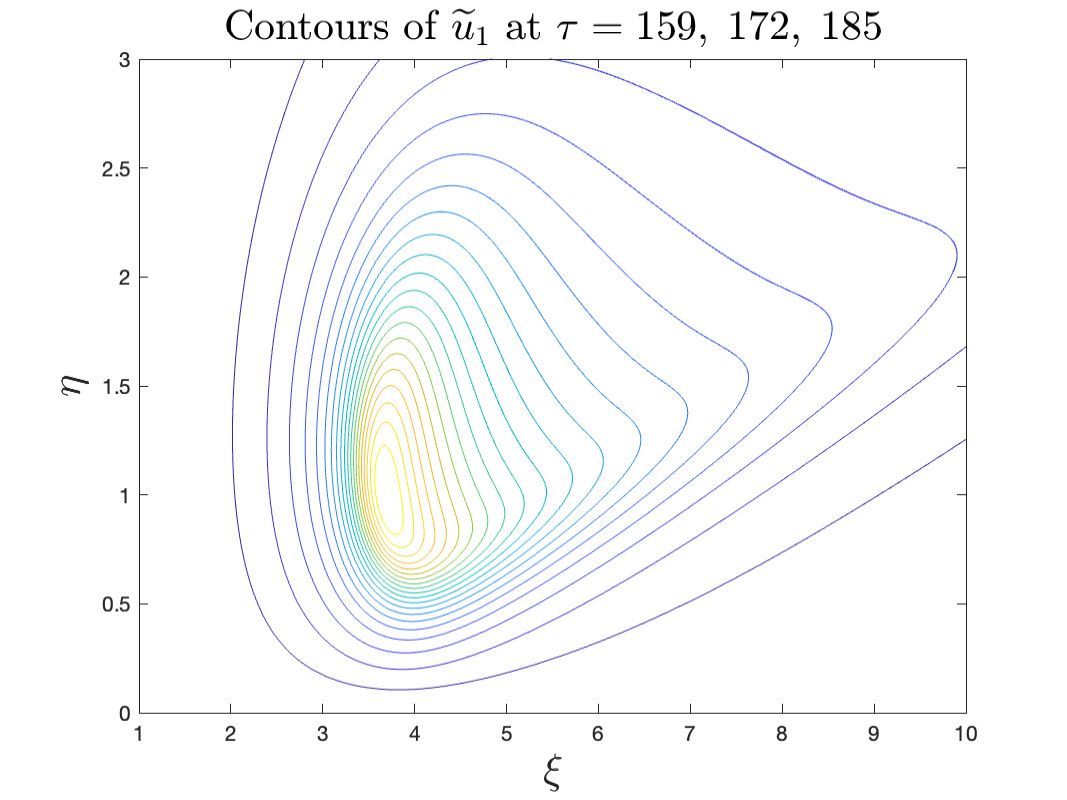}
 \includegraphics[width=0.4\textwidth]{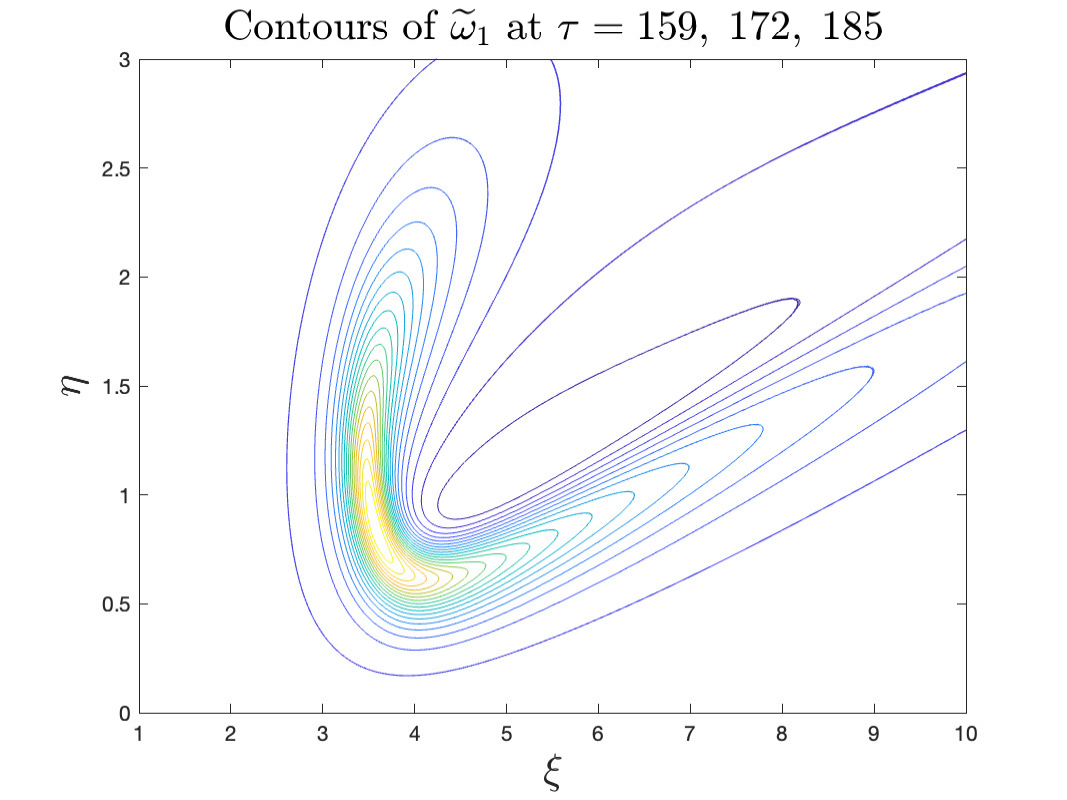} 
    \caption[VelL3-time]{ Left plot: Contours of $\widetilde{u}_1$ with respect to $(\xi,\eta)$ for the generalized Navier-Stokes equations. Right plot: Contours of $\widetilde{\omega}_1$ with respect to $(\xi,\eta)$ at $\tau = 159, \;172,\; 185$ during which the maximum vorticity has grown by a factor of $1029$.} 
    \label{fig:Contours-decay}
\end{figure}

 In Figure \ref{fig:Contours-decay}, we plot the contours of $\widetilde{u}_1$  and $\widetilde{\omega}_1$ 
as a function of $(\xi, \eta)$ for three different time instants, $\tau=159, \; 172, \; 185$ using resolution $1024\times1024$. 
During this time interval, the maximum vorticity has increased  by a factor of $1029$. We observe that these contours are almost indistinguishable from each other. We conjecture that the dynamically rescaled solution converges to a steady state as $\tau \rightarrow \infty$, which implies that ${u}_1$  and ${\omega}_1$ develop a self-similar blowup with scaling proportional to 
$C_{lz} \sim (T-t)^{c_l}$ and $C_{lr} \sim (T-t)^{c_l}$ ($c_l = 0.523$). 

\section{Blowup of the rescaled Navier--Stokes model with constant viscosity}
\label{sec:nse}
In this section, we will investigate the nearly self-similar blowup of the rescaled Navier--Stokes model with two constant viscosity coefficients. If we choose the two viscosity coefficients to be the same constant viscosity $\nu_0$, we find that the solution of the rescaled Navier--Stokes model either develops a turbulent flow if $\nu_0$ is too small or becomes a laminar flow if $\nu_0$ is too large.  After performing many experiments, we find that $\nu_1=6\cdot 10^{-4}$ and $\nu_2 = 6\cdot 10^{-3}$ seem to give robust one-scale nearly self-similar blowup. This choice of viscosity coefficients produces a stable nonlinear alignment of vortex stretching and nearly self-similar scaling properties.  

We solve the following dynamic rescaling formulation in our numerical study:
\begin{subequations}
\label{Dyn-rescale-eqn-rescaled-NS}
\begin{align}
& \widetilde{\Gamma}_{\tau}+c_{lr} \xi \widetilde{\Gamma}_{\xi} + c_{lz}\eta \widetilde{\Gamma}_{\eta}+\widetilde{\bf u} \cdot \nabla_{(\xi,\eta)} \widetilde{\Gamma} = c_\Gamma \widetilde{\Gamma} + \frac{\nu_1 C_\psi}{C_{lz}} \widetilde{\Delta} \widetilde{\Gamma} ,\\
&\widetilde{\omega}_{1,\tau}+c_{lr} \xi \widetilde{\omega}_{1,\xi} + c_{lz}\eta \widetilde{\omega}_{1,\eta}+\widetilde{\bf u} \cdot \nabla_{(\xi,\eta)} \widetilde{\omega}_1 = c_\omega \widetilde{\omega}_1 + (\widetilde{u}_1^2)_\eta + \frac{\nu_2 C_\psi}{C_{lz}}\Delta \widetilde{\omega}_1 ,\\
& - \Delta \widetilde{\psi}_1 = \widetilde{\omega}_1 , \quad \Delta =\left(\delta^2\partial_{\xi}^2+\delta^2 \frac{n}{\xi}\partial_{\xi} + \partial_{\eta}^2\right), 
\end{align}
\end{subequations}
where $\widetilde{\bf u} = (\widetilde{u}^\xi,\widetilde{u}^\eta)$, $\widetilde{u}^\xi = -\xi \widetilde{\psi}_{1,\eta}$, $\widetilde{u}^\eta = (m-1) \widetilde{\psi}_1+\xi\widetilde{\psi}_{1,\xi} $, $m=1+2R(\tau)/Z(\tau)$, $n=2m-3$, $c_\Gamma = c_u + 2c_{lr} = 2(c_{lr}-c_{lz})$ and $\widetilde{\Delta}$ is defined in \eqref{Gamma-Delta}.
We still choose $c_{lr}$ and $c_{lz}$ to enforce that $\widetilde{u}_1$ achieves its maximum at $(\xi,\eta) = (R_0,1)$ and choose $c_u$ to fix the maximum value of $\widetilde{u}_1$ to be $1$. We update $c_\psi$ and $c_\omega$ using $c_\psi = c_u +c_{lz}$ and $c_\omega = c_u - c_{lz}$. 

In the rest of this section, we will present numerical evidences of finite time self-similar blowup of the rescaled Navier-Stokes model  
\eqref{eq:axisymmetric_NSE_01} with two constant viscosity coefficients.

\subsection{Rapid growth of maximum vorticity}

In this subsection, we investigate how the profiles of the solution evolve in time. We will use the numerical results computed on the adaptive mesh of resolution $(n_1,n_2) = (1024,1024)$. We have computed the numerical solution up to time $\tau=155$. 

\begin{figure}[!ht]
\centering
    \begin{subfigure}[b]{0.40\textwidth}
        \centering
        \includegraphics[width=1\textwidth]{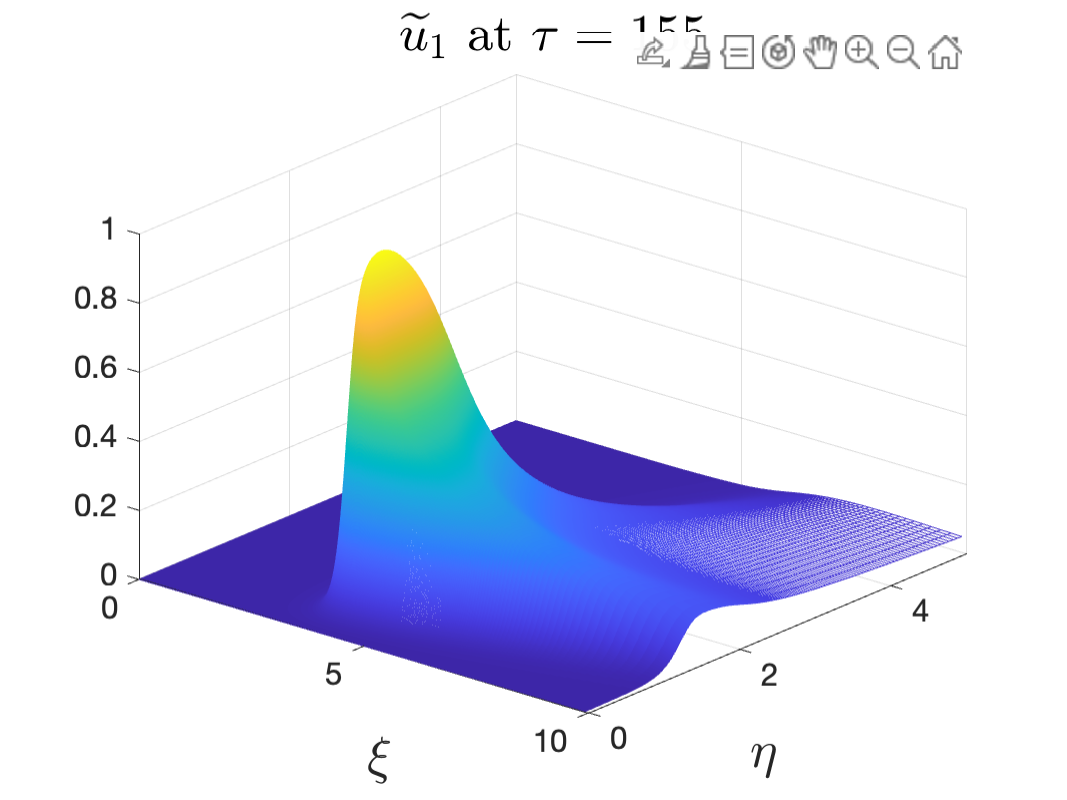}
        \caption{Rescaled profile for $\widetilde{u}_1$}
    \end{subfigure}
    \begin{subfigure}[b]{0.40\textwidth}
        \centering
        \includegraphics[width=1\textwidth]{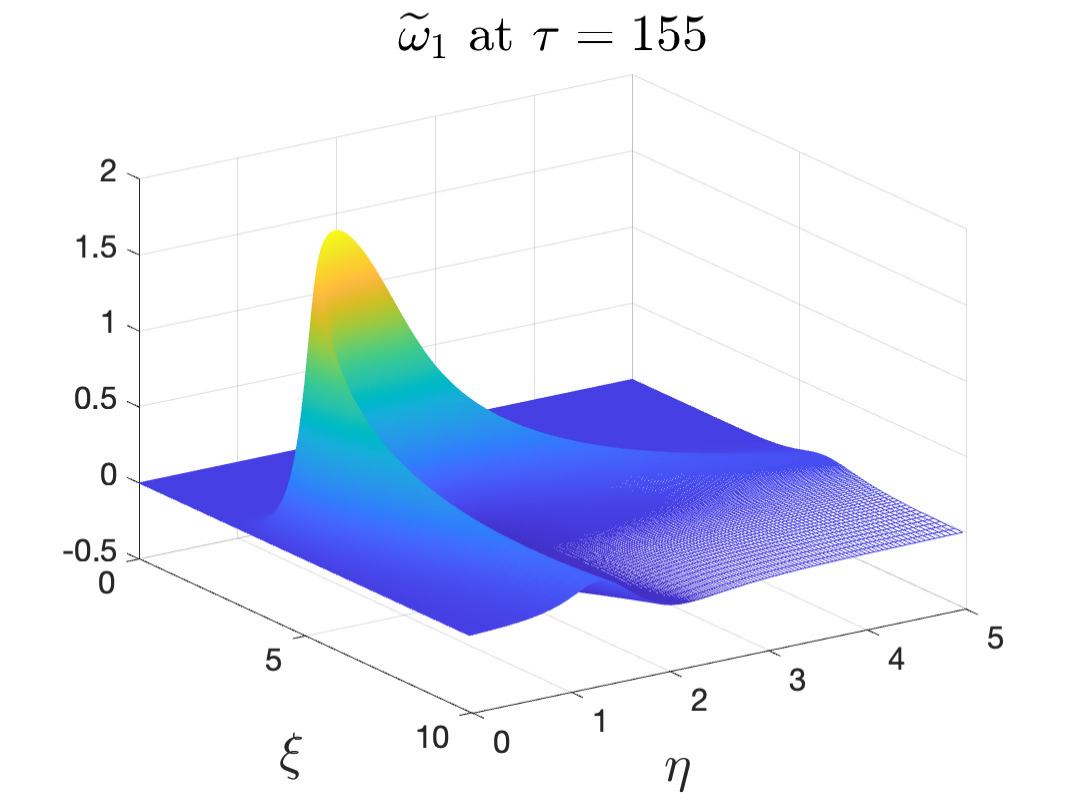}
        \caption{Rescaled profile for $\widetilde{\omega}_1$}
    \end{subfigure}
    \begin{subfigure}[c]{0.40\textwidth}
        \centering
        \includegraphics[width=1\textwidth]{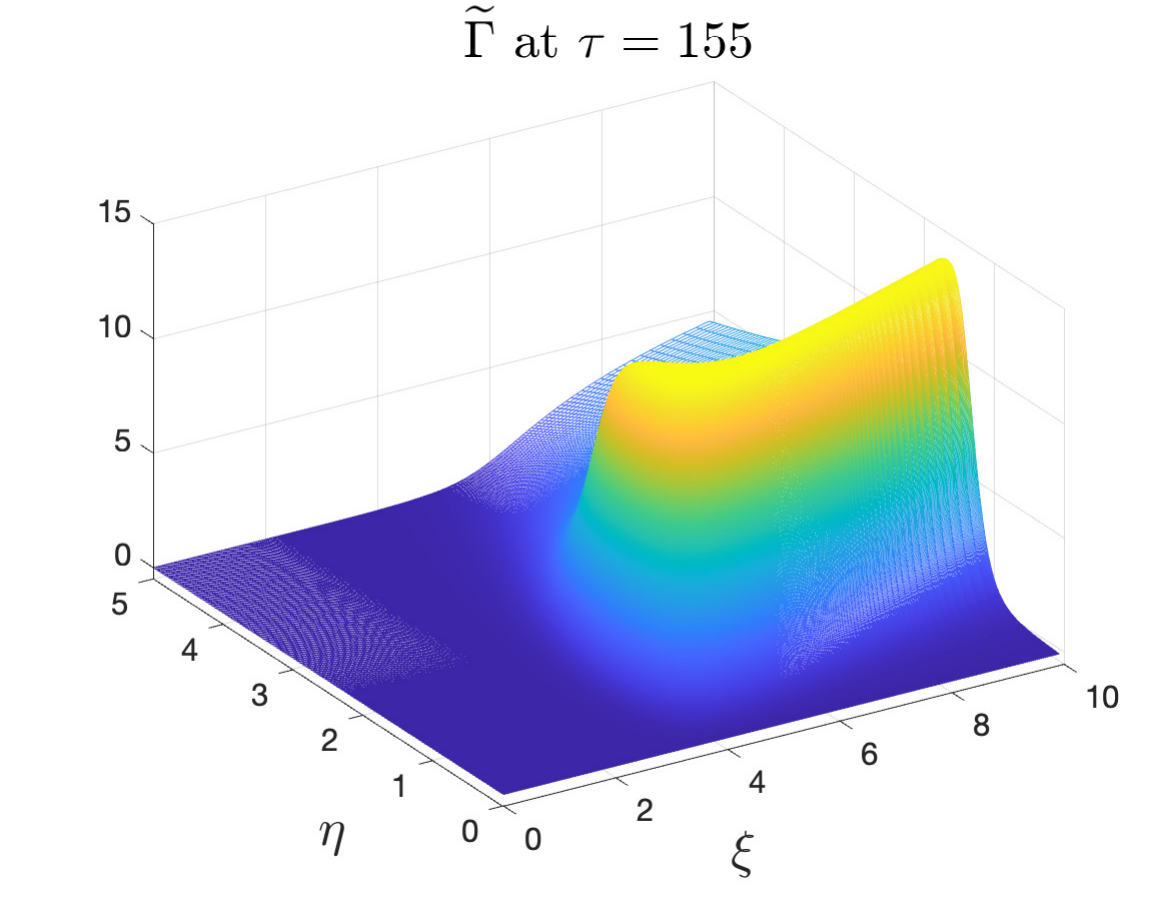}
        \caption{Rescaled profile for $\widetilde{\Gamma}$}
    \end{subfigure}
    \begin{subfigure}[d]{0.40\textwidth}
        \centering
        \includegraphics[width=1\textwidth]{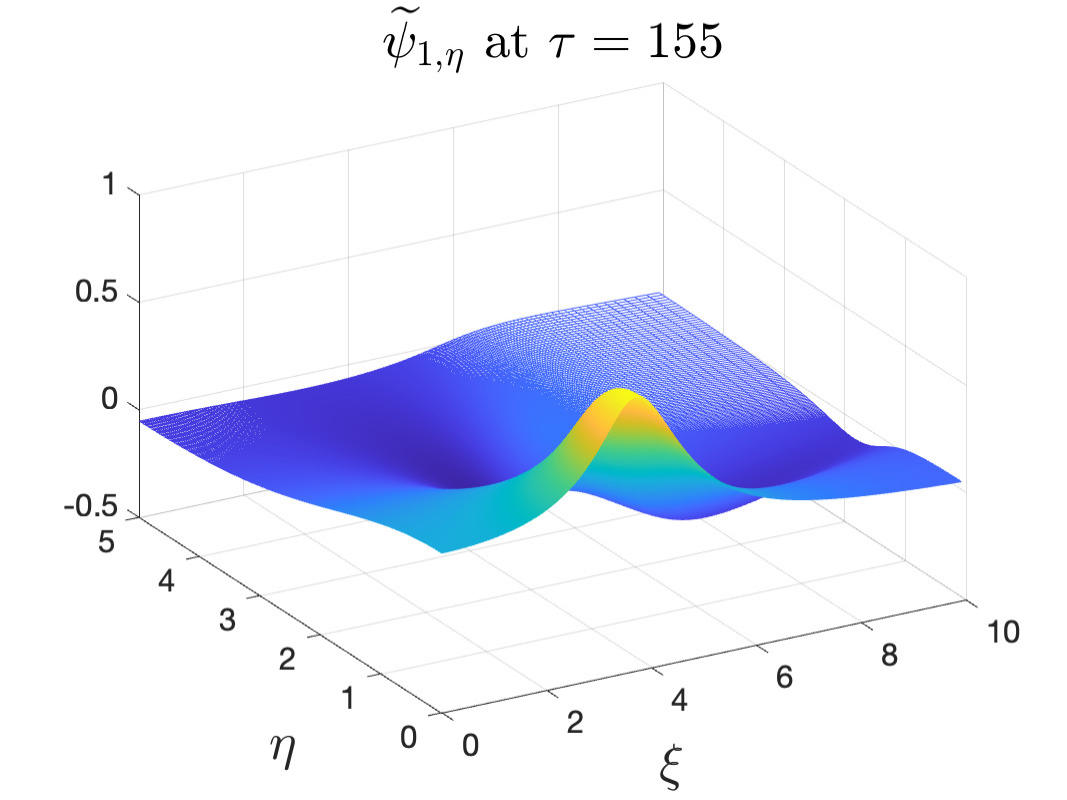}
        \caption{Rescaled profile for $\widetilde{\psi}_{1\eta}$}
    \end{subfigure}
    \caption[Profile]{The local view of rescaled profiles for the rescaled Navier-Stokes model at time $\tau = 155$. (a) $\widetilde{u}_1$; (b) $\widetilde{\omega}_1$; (c) $\widetilde{\Gamma}$; (d) $\widetilde{\psi}_{1,\eta}$.}  
     \label{fig:profiles}
\end{figure}

In Figure \ref{fig:profiles}, we present the solution profiles of $(\widetilde{u}_1$, $\widetilde{\om}_1$, $\widetilde{\Gamma}$, $\widetilde{\psi}_{1,\eta}$)
 at time $ \tau =155$. By this time, the maximum vorticity has increased by a factor of $1.4 \cdot 10^{30}$. We observe that the singular support of the profiles travels toward the origin with distance of order $O(10^{-15})$. 
 
The solution profiles look qualitatively similar to those for the generalized Navier--Stokes equations with solution dependent viscosity. Due to the relatively small viscosity $\nu_1$ for $\widetilde{\Gamma}$, we observe that the density $\widetilde{\Gamma}$ forms a shock like traveling wave  profile, propagating toward the symmetry axis $\xi = 0$. This sharp shock like profile induces a Delta function like source term for the $\widetilde{\omega}_1$ equation. The relatively large viscosity $\nu_2$ then regularizes this nearly singular source term and generates a regularized Delta function like profile for $\widetilde{\omega}_1$. We observe that $\widetilde{\psi}_{1,\eta}$ achieves its maximum value at $\eta = 0$ near $\xi = R_0$. As we commented earlier, this is a crucial property that overcomes the destabilizing effect of the transport along the $\eta$ direction, which pushes the solution upward away from $\eta=0$.

We observe that the rescaled profile $\widetilde{\omega}_1$ decays rapidly in the far field with boundary values of $O(10^{-21})$. The rescaled profile $\widetilde{\psi}_1$ also has a fast decay in the far field with boundary values of order $O(10^{-8})$. Although the boundary values of $\widetilde{\Gamma}$ are $O(1)$, the vortex stretching term $(\widetilde{\Gamma}^2/\xi^4)_\eta$ is extremely small at the far field boundary of order $O(10^{-35})$. Moreover,  the boundary values of the transport terms for $\widetilde{\Gamma}$ and $\widetilde{\omega}_1$ are of order $O(10^{-13})$ and $O(10^{-34})$, respectively.

\begin{figure}[!ht]
\centering
    \includegraphics[width=0.4\textwidth]{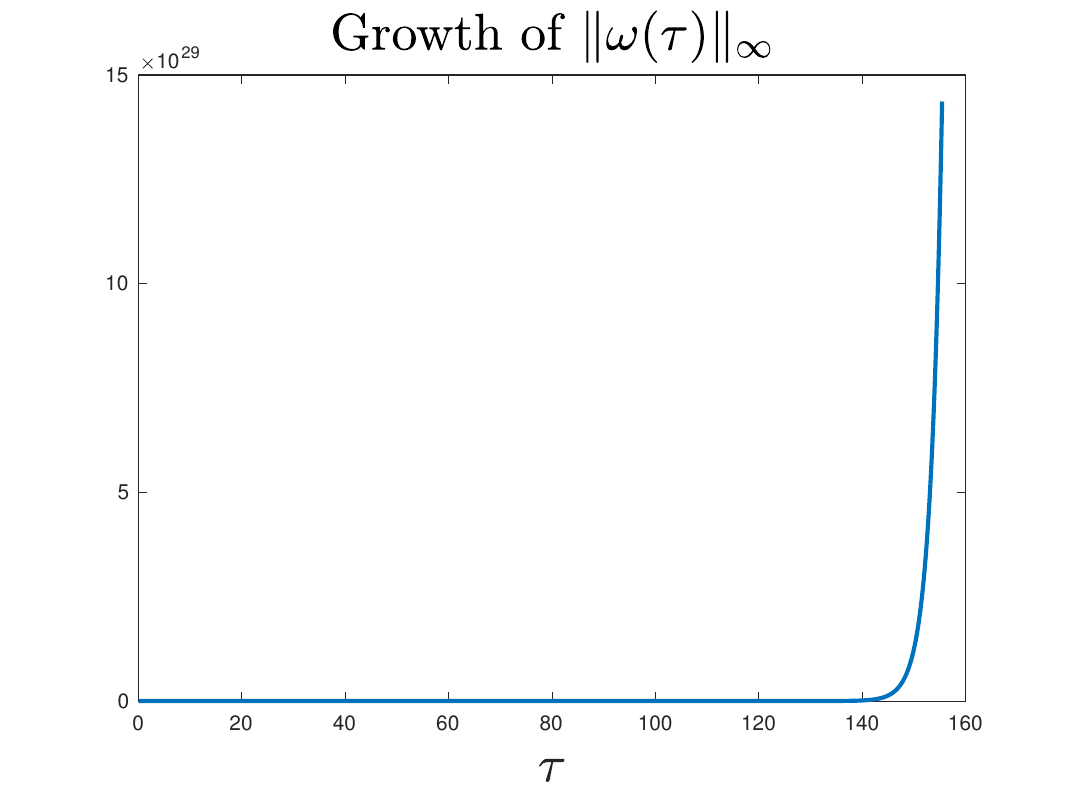}
    \includegraphics[width=0.4\textwidth]{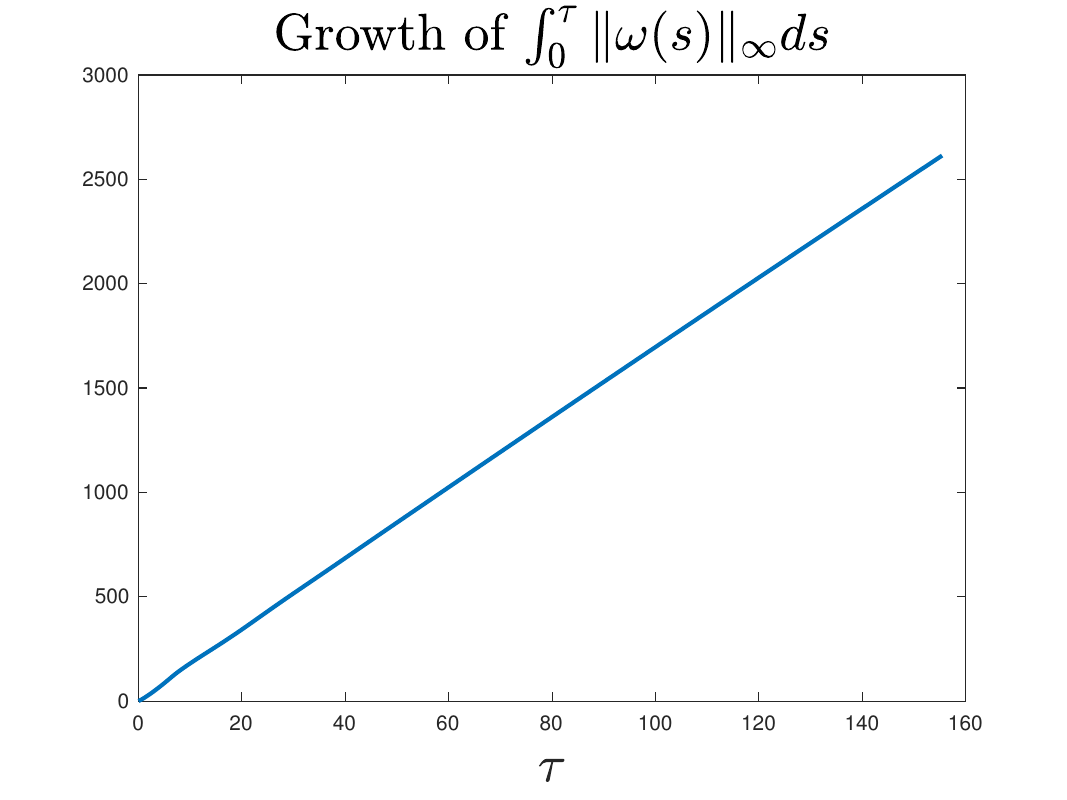} 
    \caption[Rapid growth]{  Growth of maximum vorticity for the rescaled Navier-Stokes model. Left plot: the amplification of maximum vorticity relative to its initial maximum vorticity, $\|\vom (\tau)\|_{L^\infty}/\|\vom (0)\|_{L^\infty}$ as a function of time. Right plot: the time integral of maximum vorticity,  $\|\widetilde{\vom}(0)\|_{L^\infty}^{-1}\int_0^\tau \|\widetilde{\vom} (s)\|_{L^\infty}ds$ as a function of time. The solution is computed using the $1024\times 1024$ grid. The final time instant is $\tau = 155$.} 
    \label{fig:rapid_growth_nse}
\end{figure}

\begin{figure}[!ht]
\centering
    \includegraphics[width=0.4\textwidth]{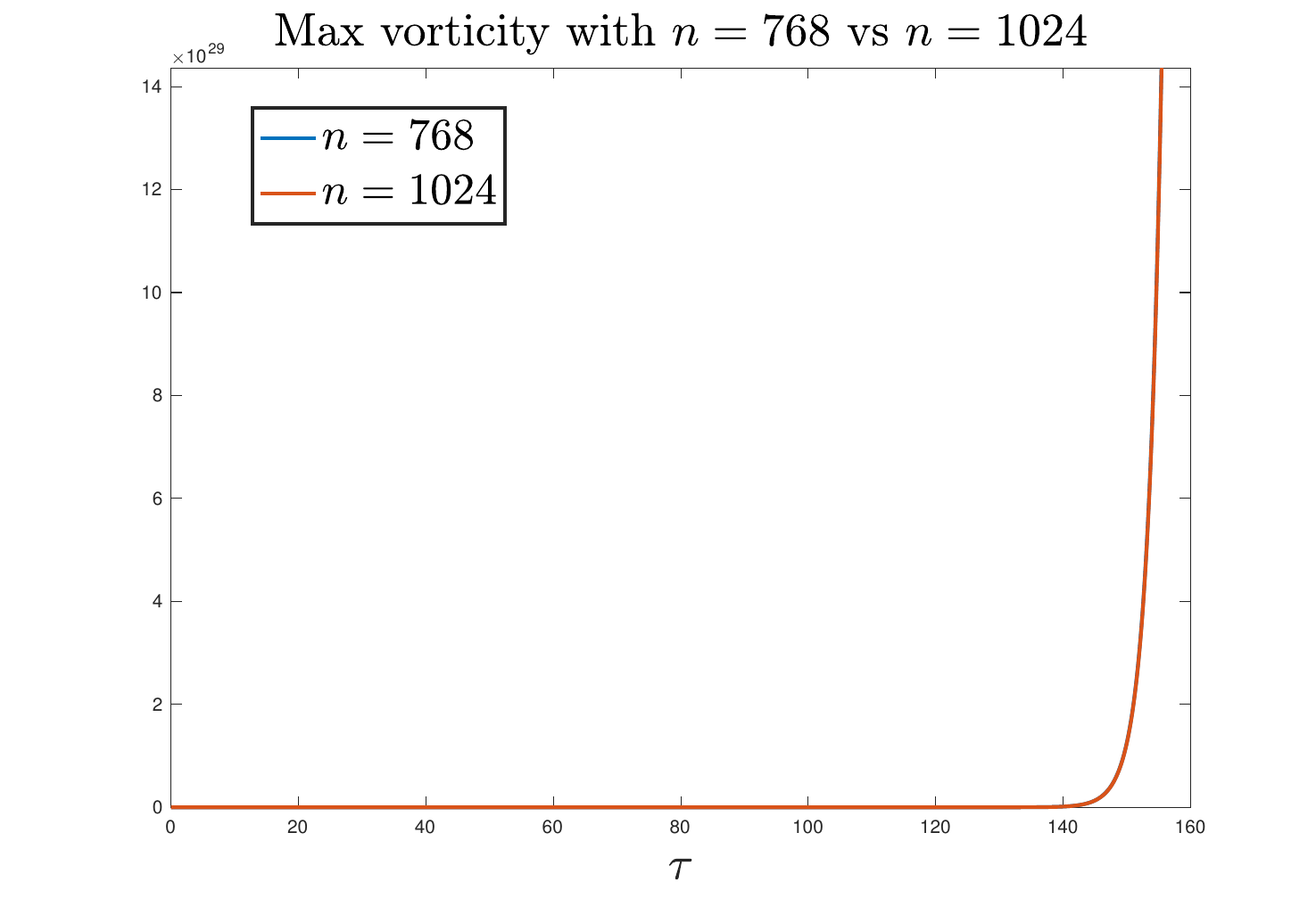}
    \includegraphics[width=0.4\textwidth]{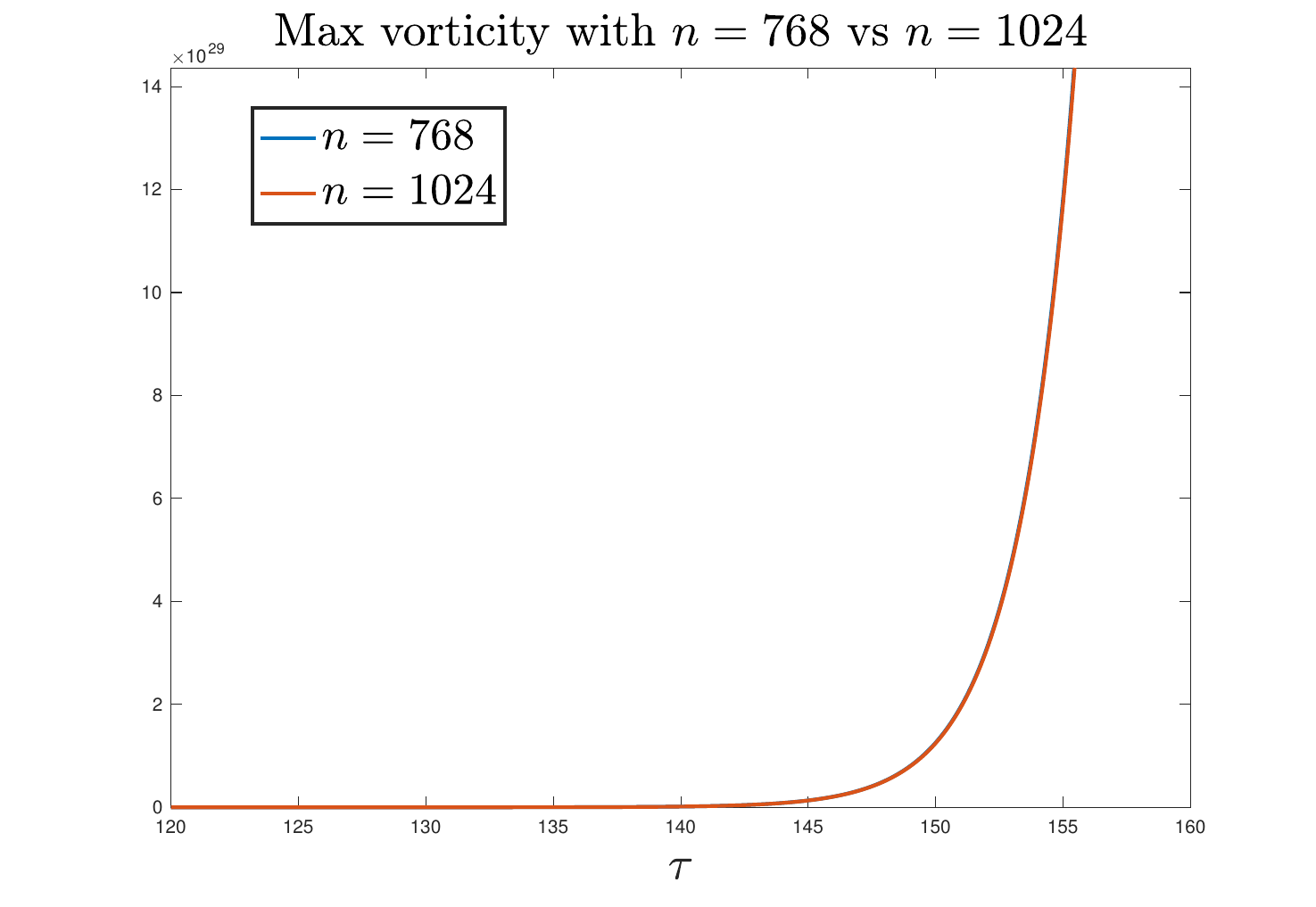} 
    \caption[Max vorticity compare]{  Left plot: Comparison of $\|\vom (\tau)\|_{L^\infty}/\|\vom (0)\|_{L^\infty}$ in time, $n=768$ (blue) vs $n=1024$ (red)  for the rescaled Navier-Stokes model. Right plot:  Zoomed-in version. Here $n$ stands for the numerical resolution of using an $n\times n$ grid.} 
    \label{fig:max_vort-comparison_nse}
\end{figure}

We observe that the solution develops rapid growth dynamically. In the left subplot of Figure \ref{fig:rapid_growth_nse}, we compute the relative growth of the maximum vorticity $\|\vom (t)\|_{L^\infty}/\|\vom (0)\|_{L^\infty}$ as a function of $\tau$.
We can see that the maximum vorticity grows extremely rapidly in time. We observe that $\|\omega (t)\|_{L^\infty}/\|\omega (0)\|_{L^\infty}$ has increased by a factor of $1.4\cdot 10^{30}$ by the end of the computation. To best of our knowledge, such a large growth rate of the maximum vorticity has not been reported for the $3$D incompressible Navier--Stokes equations in the literature. 

In the right subplot of Figure \ref{fig:rapid_growth_nse}, we plot that the time integral of the maximum vorticity in the rescaled time $\tau$, i.e. $\int_0^\tau \|\widetilde{\vom} (s) \|_{L^\infty}ds$ as a function of $\tau$. We observe that the growth rate is roughly linear with respect to $\tau$. This implies that the growth rate of the maximum vorticity is proportional $O\left ( \frac{1}{T-t}\right )$. 
The rapid growth of $\int_0^{t}\|\vom(s)\|_{L^\infty}\idiff s$ violates the well-known Beale-Kato-Majda blow-up criterion \cite{beale1984remarks}, which implies that the generalized axisymmetric Navier--Stokes equations develop a finite time singularity.

In Figure \ref{fig:max_vort-comparison_nse} (a), we compare the growth rate of the maximum vorticity using two different resolutions, $768\times 768$ vs $1024\times 1024$. A zoomed-in version is provided in Figure \ref{fig:max_vort-comparison_nse} (b). We can see that the two curves are almost indistinguishable with the $1024\times 1024$ resolution gives a slightly faster growth. This indicates that the maximum vorticity is well resolved by our computational mesh.

\subsection{Alignment of vortex stretching}

Due to the viscous regularization, the solution becomes smoother and is more stable. We are able to compute up to a time when $(R(t),Z(t))$ is very close to the origin with distance of order $O(10^{-15})$. This is something we could not achieve for the $3$D Navier--Stokes equations in a periodic cylinder \cite{Hou-nse-2022}.

\begin{figure}[!ht]
\centering
    \includegraphics[width=0.4\textwidth]{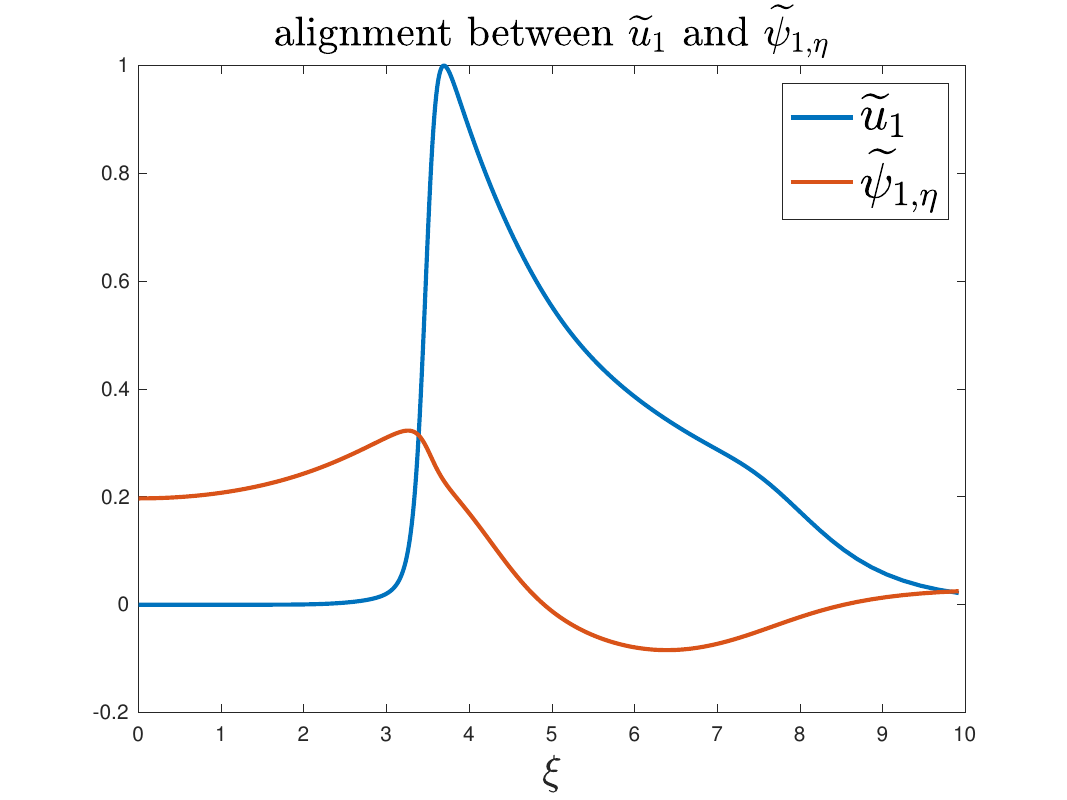}
    \includegraphics[width=0.4\textwidth]{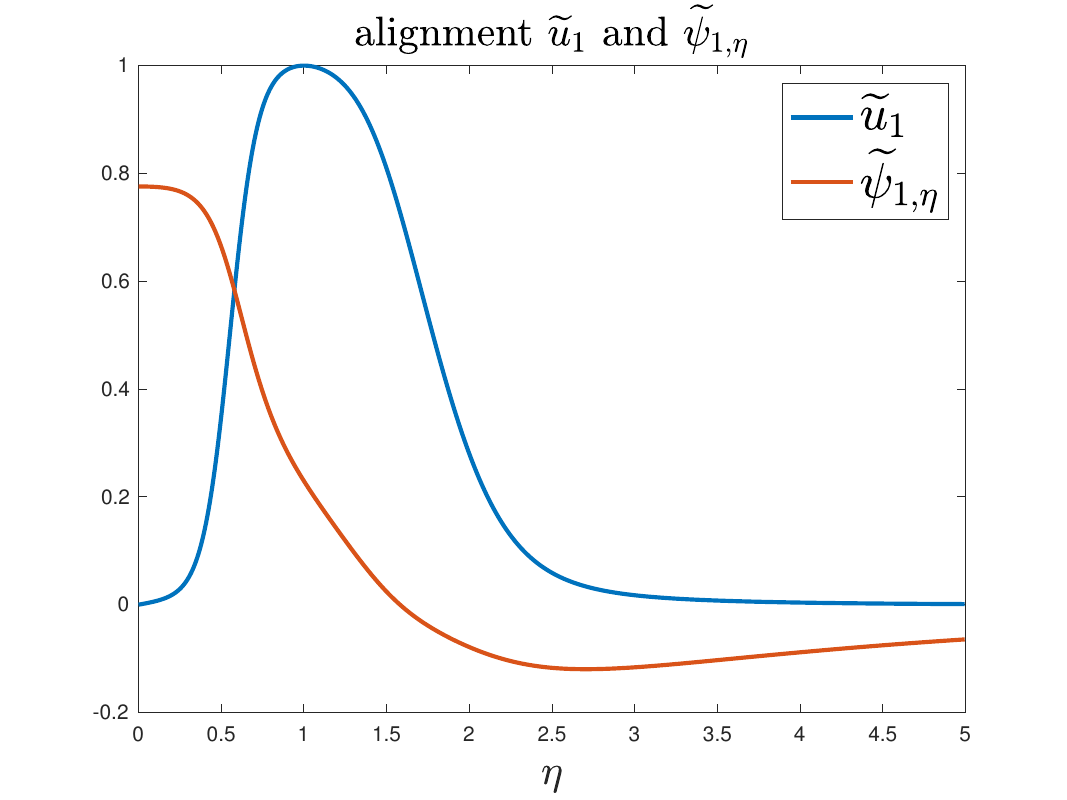} 
    \caption[Alignment-space]{ Left plot: Alignment between  $\widetilde{u}_1$ and $\widetilde{\psi}_{1,\eta}$ at $\eta = Z$ as a function of $\xi$  for the rescaled Navier-Stokes model at $\tau = 155$. We observe a sharper front along the $\xi$ direction due to the fact that we use a smaller viscosity $\nu_1$. Right plot:  Alignment between $\widetilde{u}_1$ and $\widetilde{\psi}_{1,\eta}$ at $\xi = R$ as a function of $\eta$ at $\tau = 155$.} 
    \label{fig:alignment-space}
\end{figure}

\begin{figure}[!ht]
\centering
    \includegraphics[width=0.4\textwidth]{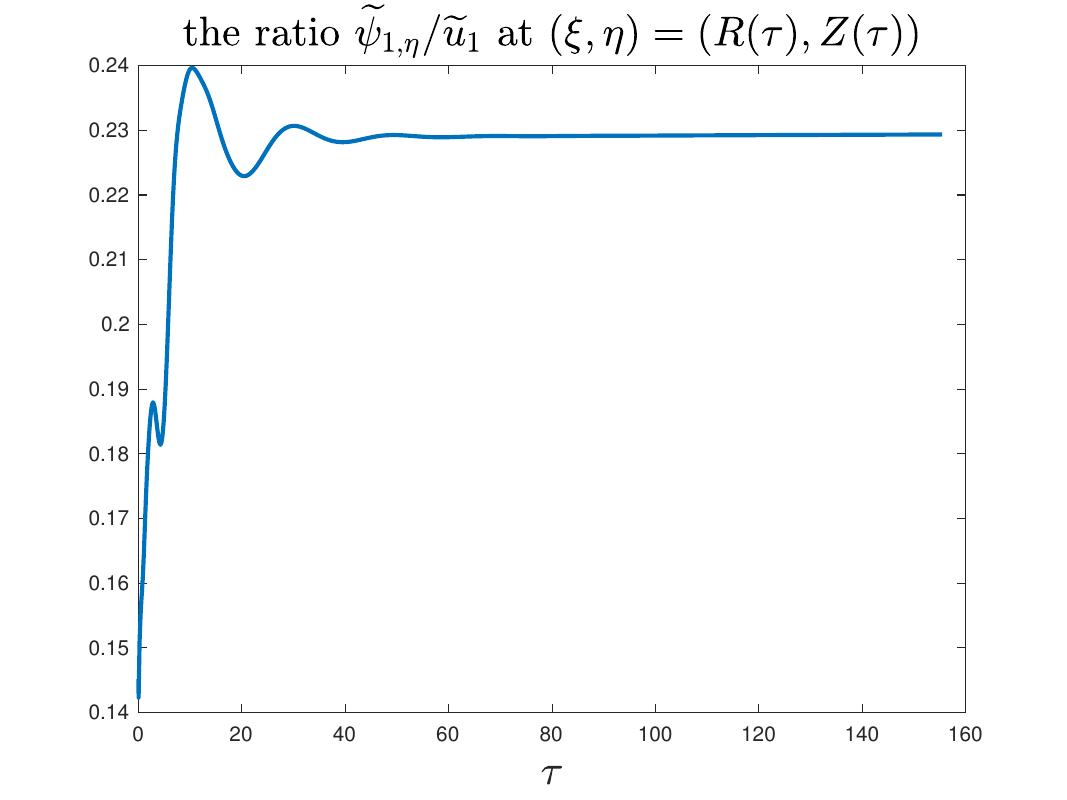}
 \includegraphics[width=0.4\textwidth]{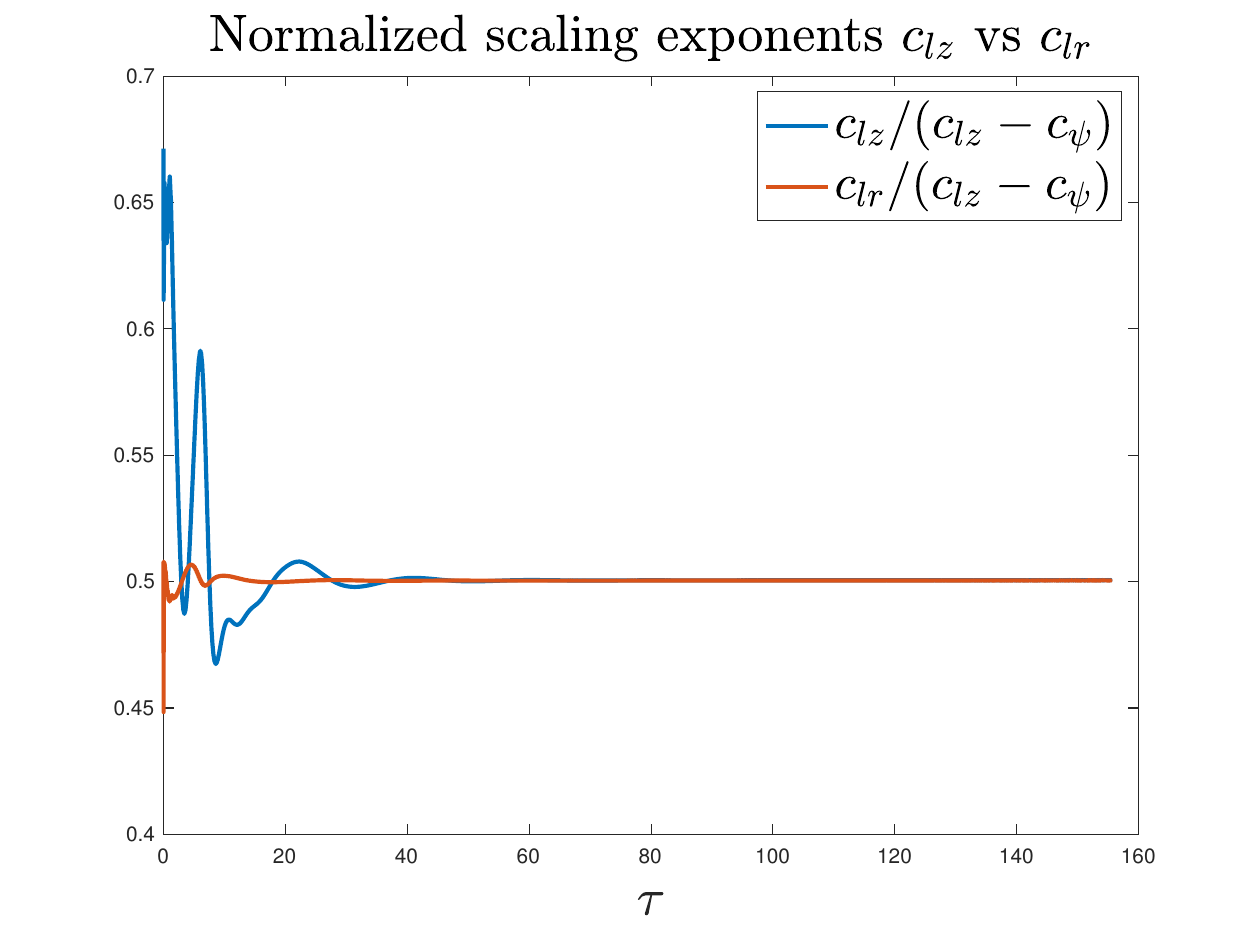} 
    \caption[Alignment-time]{ Left plot: The ratio between $\widetilde{\psi}_{1,\eta}$ and $\widetilde{u}_1$ at $(\xi,\eta)=(R,Z)$ as a function of $\tau$ for the rescaled Navier-Stokes model. Right plot: The normalized scaling exponent $c_{lz}/(c_{lz}-c_\psi)$ and $c_{lr}/(c_{lz}-c_\psi)$ as a function of $\tau$.} 
    \label{fig:alignment-time-nse}
\end{figure}

In Figure \ref{fig:alignment-space}(a)-(b), we demonstrate the alignment between $\widetilde{\psi}_{1\eta}$ and $\widetilde{u}_1$ at $\tau =155$. Although the maximum vorticity has grown so much by this time, the local solution structures have remained qualitatively the same in the late stage of the computation. 
We observe that the viscous effect actually enhances the nonlinear alignment of vortex stretching. Although we use two constant viscosity coefficients here,  we observe qualitatively the same phenomena as we did for the generalized Navier--Stokes equations with solution dependent viscosity. In particular, $\widetilde{\psi}_{1\eta}$ is relatively flat in a local region near the origin. This is an essential property that prevents the formation of a two-scale structure. 
Moreover, we observe the same qualitative positive feedback mechanism as we observed for the generalized Navier--Stokes equations with solution dependent viscosity.

In Figure \ref{fig:alignment-time-nse}(a), we observe that the positive alignment between $\widetilde{u}_1$ and $\widetilde{\psi}_{1
\eta}$ is almost flat in time with a mild increase in the late stage of the computation. This indicates that the rescaled Navier--Stokes model achieves a nearly self-similar scaling relationship.

In the Figure \ref{fig:alignment-time-nse} (b), we plot the normalized scaling exponents $c_{lz}/(c_{lz}-c_\psi)$ and $c_{lr}/(c_{lz}-c_\psi)$ as a function of $\tau$. We observe that they seem to approach the same value $0.5$ as $\tau$ increases. This is again due to the stabilizing effect of varying the space dimension to eliminate the scaling instability. If we kept $n=3$, the solution developed a two-scale solution structure. Using two different viscosity coefficients seems to play an essential role for us to obtain nearly parabolic scaling property in the sense that $c_{lz}/(c_{lz}-c_\psi)$ and $c_{lr}/(c_{lz}-c_\psi)$ approach to $0.5$ as $\tau$ increases. This is something that we could not have accomplished if we used the same viscosity coefficients for both $\Gamma$ and $\omega_1$ equations.


\subsubsection{The streamlines} 
\begin{figure}[!ht]
\centering
    \begin{subfigure}[b]{0.40\textwidth}
        \centering
        \includegraphics[width=1\textwidth]{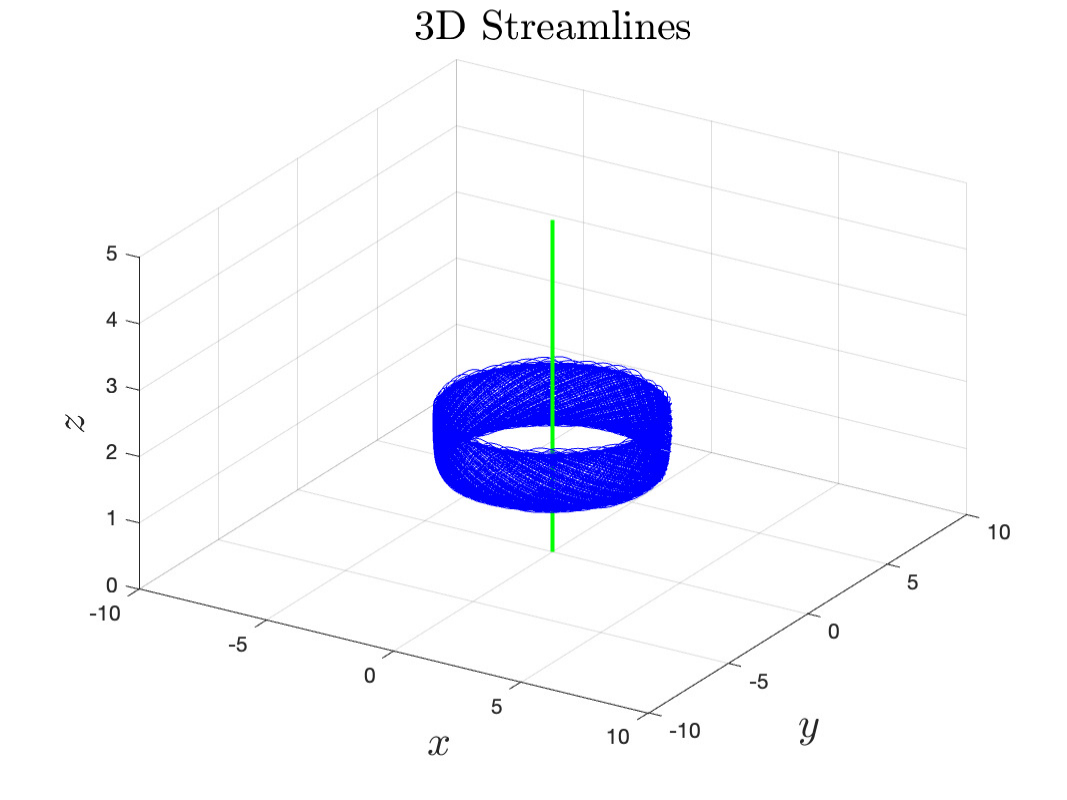}
        \caption{$r_0 = 4$, $z_0 = 2$}
    \end{subfigure}
    \begin{subfigure}[b]{0.40\textwidth}
        \centering
        \includegraphics[width=1\textwidth]{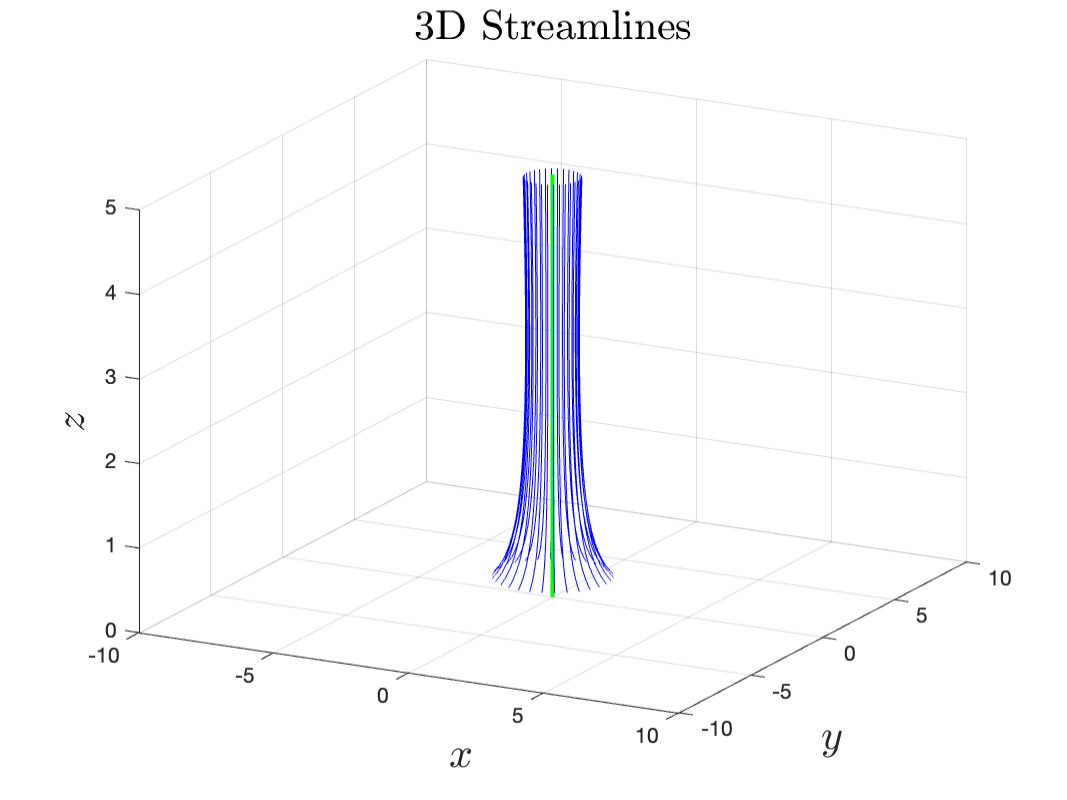}
        \caption{$r_0 = 2$, $z_0 = 0.25$}
    \end{subfigure}
    \begin{subfigure}[c]{0.40\textwidth}
        \centering
        \includegraphics[width=1\textwidth]{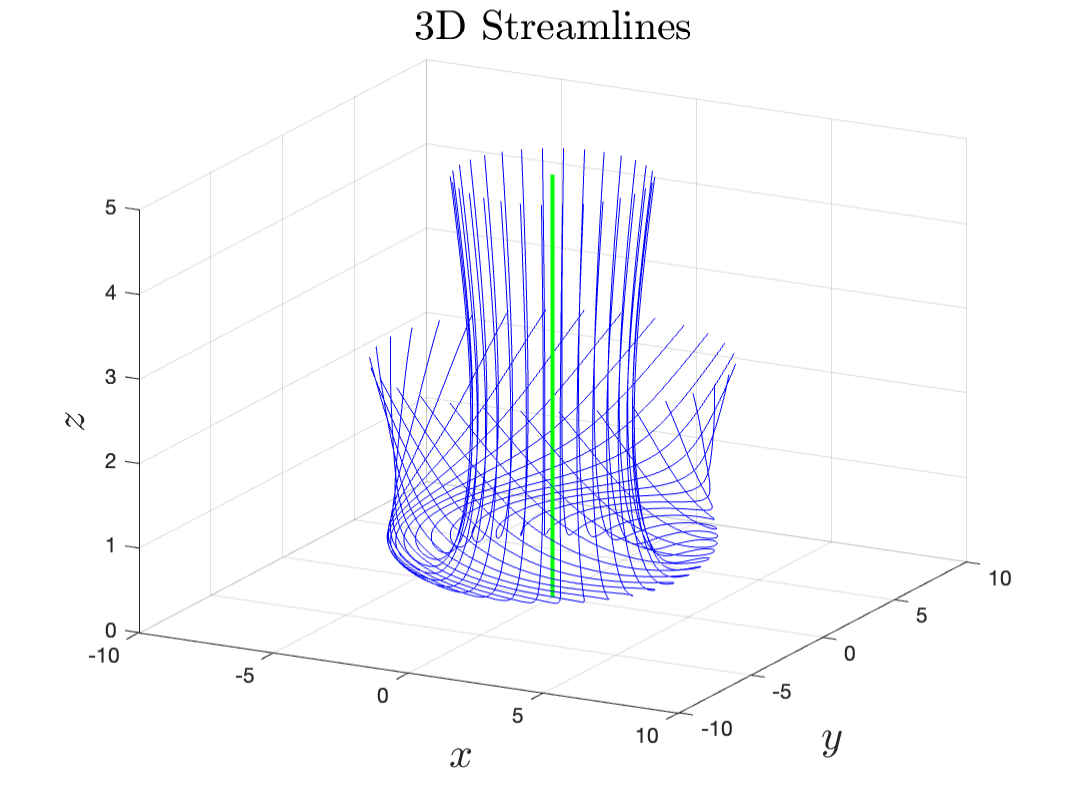}
        \caption{$r_0 = 6$, $z_0 = 2.8$}
    \end{subfigure}
    \begin{subfigure}[d]{0.40\textwidth}
        \centering
        \includegraphics[width=1\textwidth]{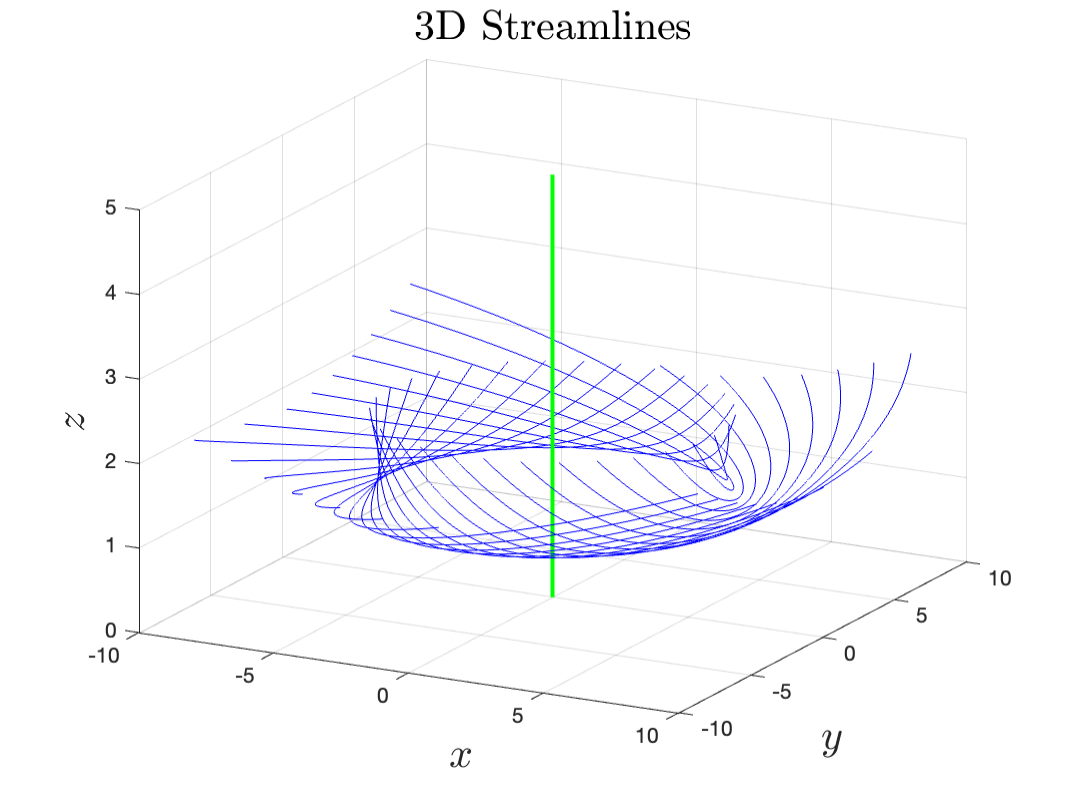}
        \caption{$r_0 = 6$, $z_0 = 2.2$}
    \end{subfigure}
    \caption[Global streamline]{The streamlines of $(u^r(t),u^\theta(t),u^z(t))$  for the rescaled Navier-Stokes model at time $\tau = 155$  with initial points given by (a) $(r_0,z_0) = (4,2)$, streamlines form a torus; (b) $(r_0,z_0) = (2,0.25)$, streamlines go straight upward; (c) $(r_0,z_0) = (6,2.8)$, streamlines first go downward, then travel inward, finally go upward; (d) $(r_0,z_0) = (6,2.2)$, streamlines spin downward and outward. The green pole is the symmetry axis.}  
     \label{fig:streamline_3D_nse}
\end{figure}

In this subsection, we investigate the features of the velocity field. We first study the velocity field by looking at the induced streamlines. In Figure~\ref{fig:streamline_3D_nse}, we plot the streamlines induced by the velocity field $\vu(t)$ at $\tau = 155$. By this time, the ratio between the maximum vorticity and the initial maximum vorticity, i.e. $\|\vom (t)\|_{L^\infty}/\|\vom (0)\|_{L^\infty}$, has increased by a factor of $1.4 \cdot 10^{30}$.

Surprisingly the induced streamlines look qualitatively the same as those obtained for the generalized Navier--Stokes equations with solution dependent viscosity. We will use a similar set of parameters to draw the streamlines to compare with the streamlines obtained for the generalized Navier--Stokes equations that we reported in the previous section. 
In  Figure~\ref{fig:streamline_3D_nse}(a) with $(r_0,z_0) = (4,1.5)$, we observe that the streamlines form the same type of torus spinning around the symmetry axis. In  Figure~\ref{fig:streamline_3D_nse}(b) with $(r_0,z_0) = (2,0.25)$, we obsevre that the streamlines go straight upward without any spinning. In  Figure~\ref{fig:streamline_3D_nse}(c) with $(r_0,z_0) = (6,2.8)$, the streamlines first go downward, then travel inward and finally go upward. Note that this starting point is different from the corresponding case for the generalized Navier-Stokes equations with solution dependent viscosity where we used $(r_0,z_0) = (6,0.5)$.  This explains why the two sets of streamlines look different. Finally, we plot in Figure~\ref{fig:streamline_3D_nse}(d) the streamlines that start with $(r_0,z_0) = (6,2)$. We can see that the streamlines first spin downward and then outward. It is also interesting to note that the solution behaves qualitatively the same as what we observed for the $3$D axisymmetric Navier--Stokes equations in a periodic cylinder \cite{Hou-nse-2022}. 

\subsubsection{The 2D flow} 
\begin{figure}[!ht]
\centering
    \includegraphics[width=0.40\textwidth]{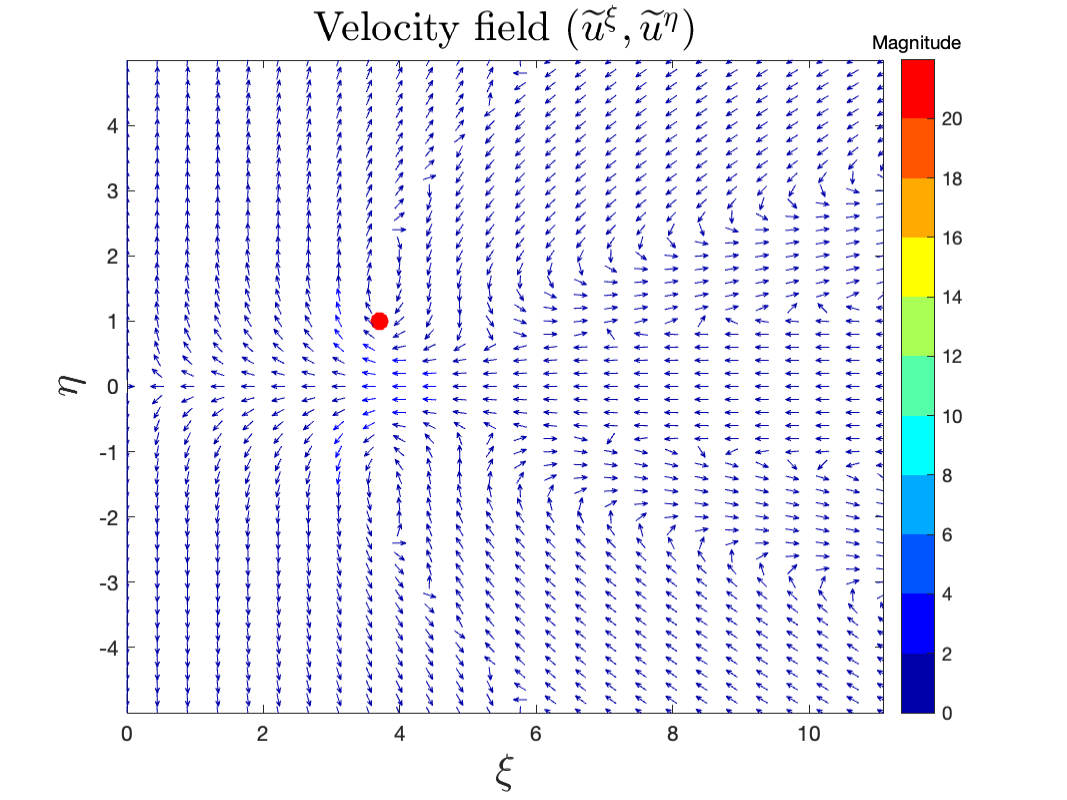}
    \includegraphics[width=0.40\textwidth]{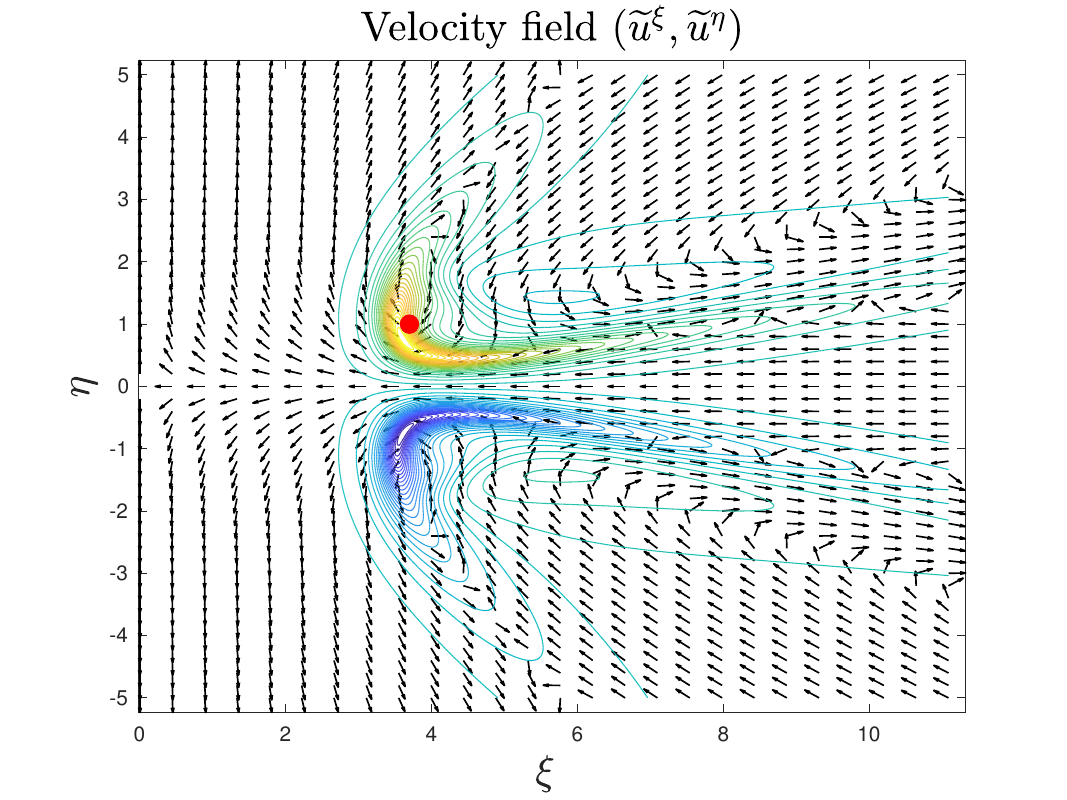} 
    \caption[Dipole]{The dipole structure of $\widetilde{\omega}_1$  and the induced dynamically rescaled velocity field $(\widetilde{u}^\xi,\widetilde{u}^\eta)$  for the rescaled Navier-Stokes model at $\tau=155$. Recall $\xi = r/C_{lr}$ and $\eta = z/C_{lz}$. The arrows in both plots are used to illustrate the direction of the velocity field. They are normalized to have the same length for better plotting quality and do not have any physical meaning here. Left plot: the velocity vector. Right plot: the velocity vector with the $\omega_1$ contour as background. The red dot is the position $(R,Z)$ where $\widetilde{u}_1$ achieves its global maximum.}  
     \label{fig:dipole_nse}
\end{figure}

To understand the phenomena in the most singular region as shown in Figure~\ref{fig:streamline_3D_nse}, we also study the $2$D velocity field $(u^r,u^z)$. In Figure~\ref{fig:dipole_nse}(a)-(b), we plot the dipole structure of $\widetilde{\omega}_1$ in a local symmetric region and the hyperbolic velocity field induced by the dipole structure in a local microscopic domain $[0,R_b]\times [0,Z_b]$ at  $\tau=155$.
The dipole structure for the rescaled Navier--Stokes model with constant viscosity looks qualitatively similar to that of the generalized Navier--Stokes equations with solution dependent viscosity. The negative radial velocity near $\eta=0$ induced by the vortex dipole pushes the solution toward $\xi=0$, then move upward away from $\eta=0$. This is the main driving mechanism for the flow to develop a hyperbolic structure. Since the value of $\widetilde{u}_1$ becomes very small near the symmetry axis $\xi=0$, the streamlines almost do not spin around the symmetry axis, as illustrated in Figure~\ref{fig:streamline_3D_nse}(b). 

\begin{figure}[!ht]
\centering
    \includegraphics[width=0.40\textwidth]{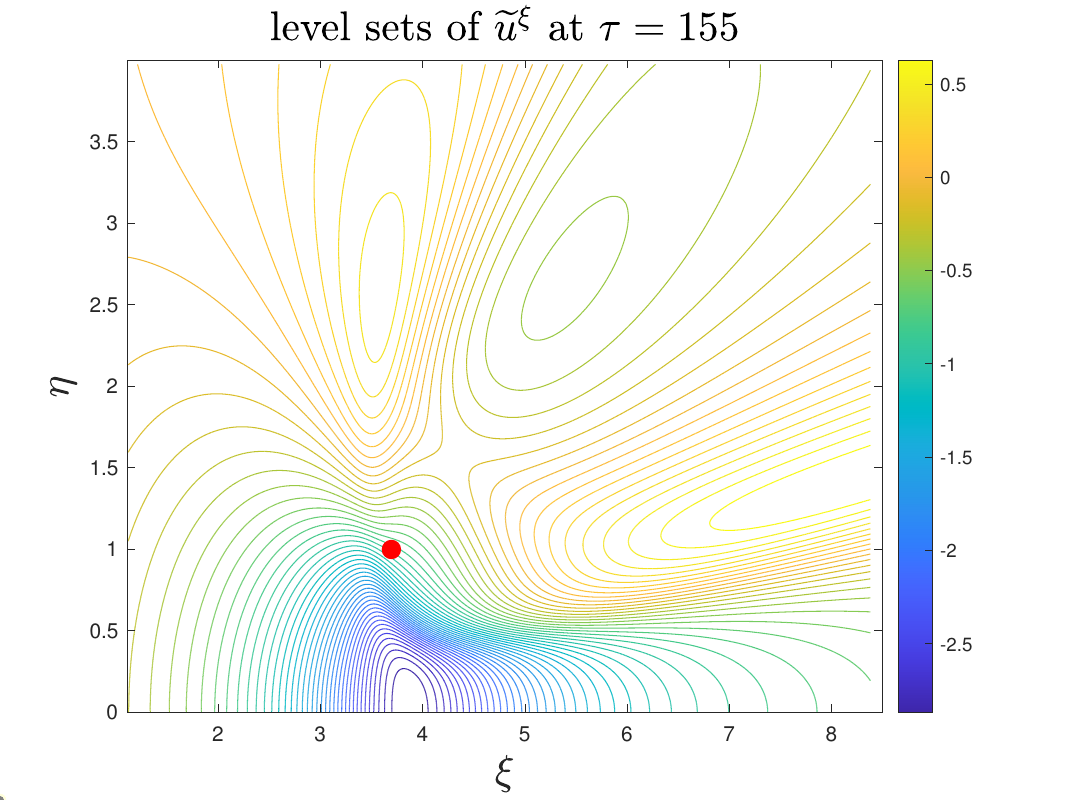}
    \includegraphics[width=0.40\textwidth]{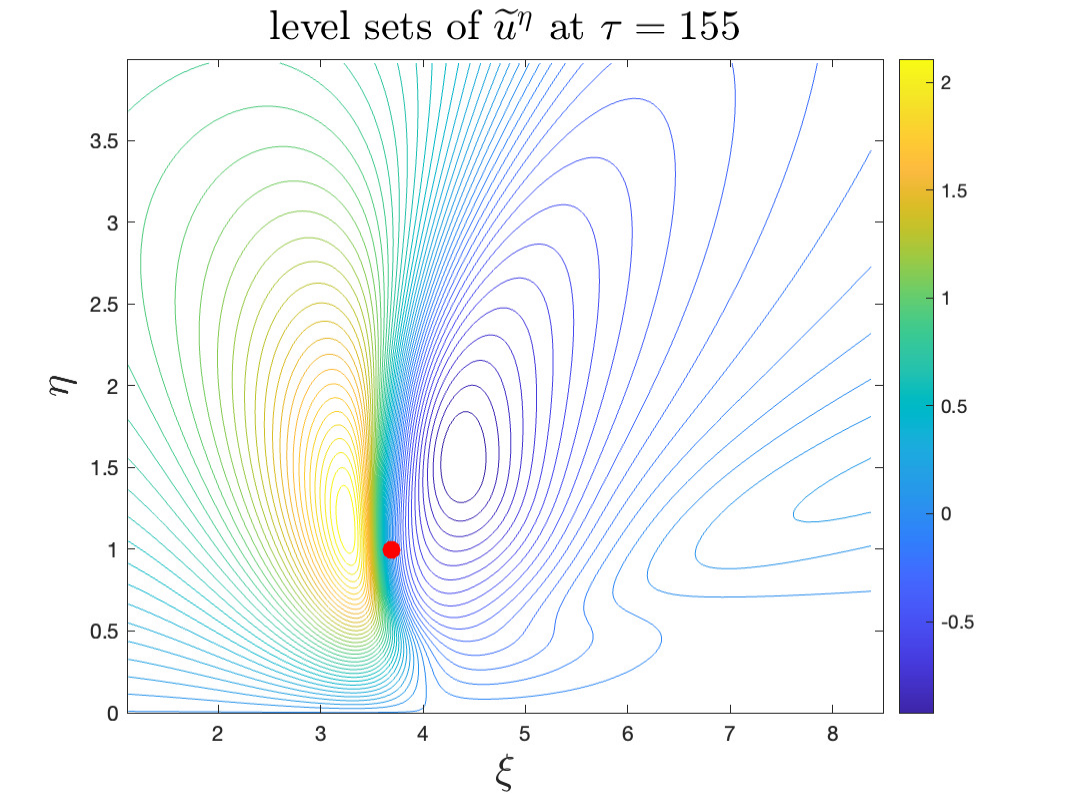} 
    \caption[Velocity level sets]{The level sets of $\tilde{u}^\xi$ (left) and $\tilde{u}^\eta$ (right) for the rescaled Navier-Stokes model at $\tau = 155$. The red point is the maximum location $(R,Z)$ of $\widetilde{u}_1$.}  
       \label{fig:velocity_levelset_nse}
\end{figure}

We can also understand this hyperbolic flow structure from the velocity contours in Figure~\ref{fig:velocity_levelset_nse} (a)-(b). Although the velocity contours look qualitatively the same as those for the generalized Navier-Stokes equations with solution dependent viscosity, a closer look shows that there are some subtle differences in the local solution structure, especially in the shape of vorticity contours of $\widetilde{\omega}_1$ if we compare Figure~ \ref{fig:dipole_nse-decay}(b) with Figure~\ref{fig:dipole_nse}(b). From Figure~\ref{fig:velocity_levelset_nse}(a), we observe that the radial velocity $u^r$ is negative and large in amplitude below the red dot $(R,Z)$, which pushes the flow toward the symmetry axis $\xi=0$. The axial velocity $u^z$ is negative and large in amplitude to the right hand side of $(R,Z)$, pushing the flow downward toward $\eta = 0$. On the left hand side of $(R,Z)$, it becomes large and positive on the left hand side of $(R,Z)$, which pushes the flow upward away from $\eta = 0$. This is very similar to the flow structure of the generalized Navier-Stokes equations with solution dependent viscosity.

As in the case of the axisymmetric Navier-Stokes equations in a periodic cylinder, the velocity field $(u^r(t),u^z(t))$ forms a closed circle right above $(R,Z)$. The corresponding streamlines are trapped in the circle region in the $\xi \eta$-plane, which is responsible for the formation of the spinning torus that we observed earlier.

\subsection{Scaling properties of the nearly self-similar blowup}

\begin{figure}[!ht]
\centering 
    \includegraphics[width=0.4\textwidth]{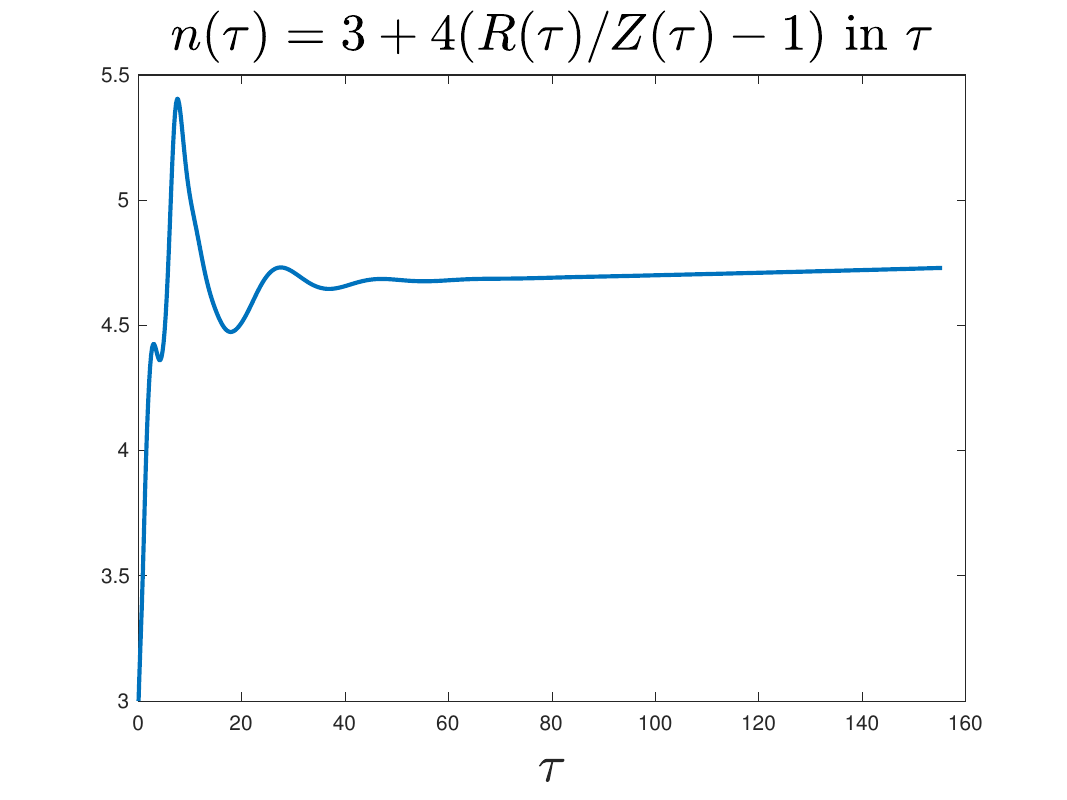}
 \includegraphics[width=0.4\textwidth]{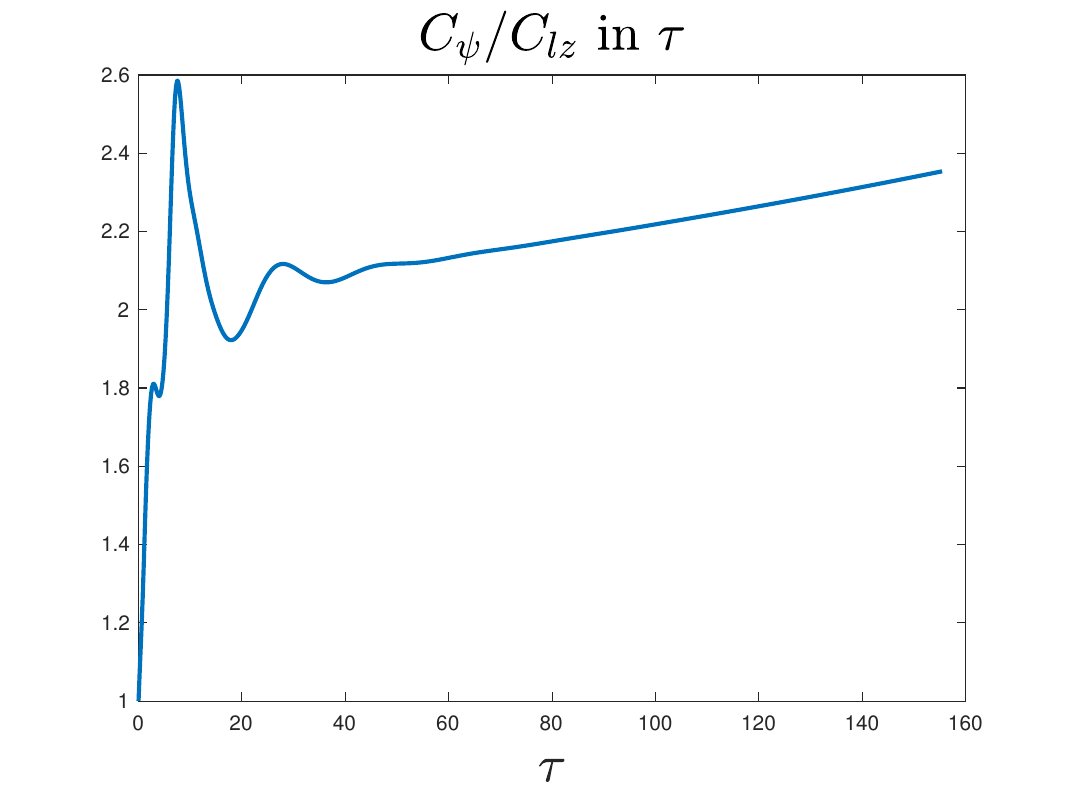} 
    \caption[Alignment-time]{ Left plot: The integration dimension for our generalized energy $n(\tau) = 3+4(R(\tau)/Z(\tau)-1)$  for the rescaled Navier-Stokes model in $\tau$ with ${n}(155)=4.73$. Right plot: $C_\psi(\tau)/C_{lz}(\tau)$ in $\tau$.} 
    \label{fig:alignment-time-2}
\end{figure}

In this subsection, we study the blowup scaling properties. We observe that all the scaling parameters $c_{lz}$, $c_{lr}$, $c_\psi$, $c_\omega$, and $c_\psi$ all converge to a constant value as $\tau$ increases. By the discussion in Section 2, we know that we should study the normalized scaling parameters defined by
$c_{lz}/(c_{lz}-c_\psi),\; c_{lr}/(c_{lz}-c_\psi), \; 
c_\psi/(c_{lz}-c_\psi)$, etc.
In Figure \ref{fig:alignment-time-2} (a), we plot the space dimension $\tilde{n}(\tau) = 3+ 4(R(\tau)/Z(\tau)-1)$ as a function of $\tau$. As we can see, $n(\tau)$ remains relatively flat in the late stage with $n(155) = 4.73$. In Figure \ref{fig:alignment-time-2}(b), we plot $C_\psi(\tau)/C_{lz}(\tau)$ as a function of $\tau$. We observe that $C_\psi(\tau)/C_{lz}(\tau)$ roughly has a linear growth with respect to $\tau$ with a very small slope $\epsilon$ since the growth of $C_\psi(\tau)/C_{lz}(\tau)$ is very small over a long time. 

From the discussion in Section $2$, we know that 
\[
\tau = c_0 \log\left ( \frac{1}{T-t} \right ), \;C_{lz} = (T-t)^{\widehat{c}_{lz}}, \; 
C_{\psi} = (T-t)^{1-\widehat{c}_{lz}}\;.
\]
If we assume that 
\[
C_\psi(\tau)/C_{lz}(\tau) = 1 + \epsilon \tau,
\]
for $\tau$ large, then we would obtain
\[
(T-t)^{1-2\widehat{c}_{lz}} = 1 + \epsilon \tau \;,
\]
which implies
\[
\widehat{c}_{lz} = \frac{1}{2} +  \frac{\log(1 +\epsilon \tau)}{2\tau}.
\]
Thus, we obtain that the convergence of $\widehat{c}_{lz}$ to $1/2$ with a  logarithmic rate only. Moreover, if $\tau$ is not large and $\epsilon $ is small, we get 
$\widehat{c}_{lz} \approx 1/2 + \epsilon/2$. 
By substituting $\widehat{c}_{lz} = \frac{1}{2} +  \frac{\log(1 +\epsilon \tau)}{2\tau}$ into $C_{lz} = (T-t)^{\widehat{c}_{lz}}$, we further obtain
\[
\lambda(t) = \frac{C_{lz}}{\sqrt{T-t}} = 
\frac{1}{\sqrt{1+\epsilon \tau}} = 
\frac{1}{\sqrt{1+c_0\epsilon |\log(T-t)|}}\;.
\] 
Since $\widehat{c}_u = c_u/(c_{lz} - c_\psi) = -1$, we conclude that 
$C_u \approx (T-t)$ and $\| u_1\|_\infty = O(1/(T-t))$,
which implies 
\[
\| \vom\|_\infty = O\left ( \frac{1}{T-t} \right )\;.
\]

\subsection{Checking against various blowup criteria}

In this subsection, we apply various blowup criteria to confirm the finite time blowup of the rescaled Navier--Stokes model with two constant viscosity coefficients. Our studies show that the nearly self-similar blowup satisfies almost all the generalized blowup criteria that have been established for the $3$D axisymmetric Navier--Stokes equations with smooth initial data. 

\subsubsection{Non-blowup criteria based on enstrophy growth}

\begin{figure}[!ht]
\centering
    \includegraphics[width=0.4\textwidth]{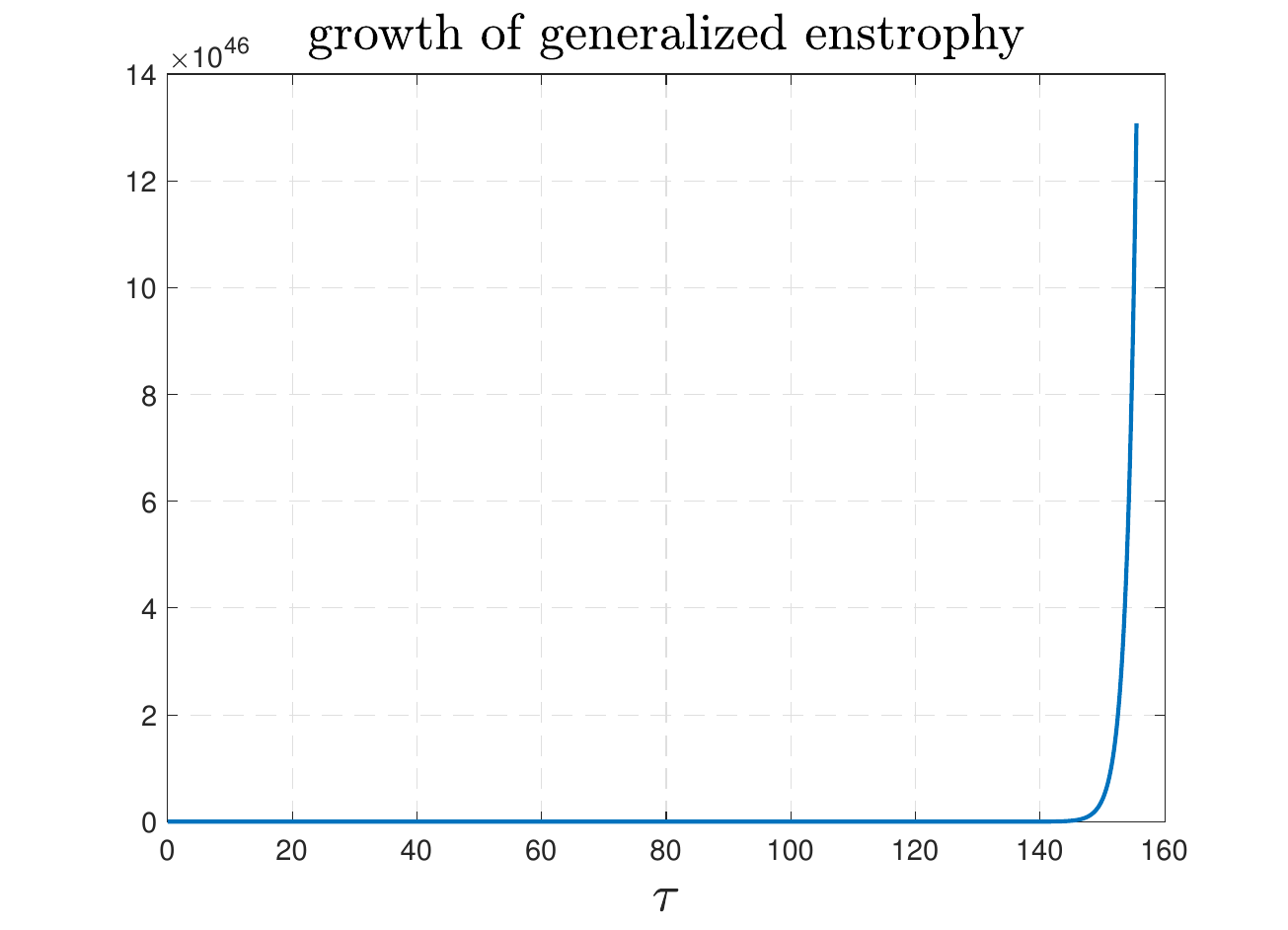}
 \includegraphics[width=0.4\textwidth]{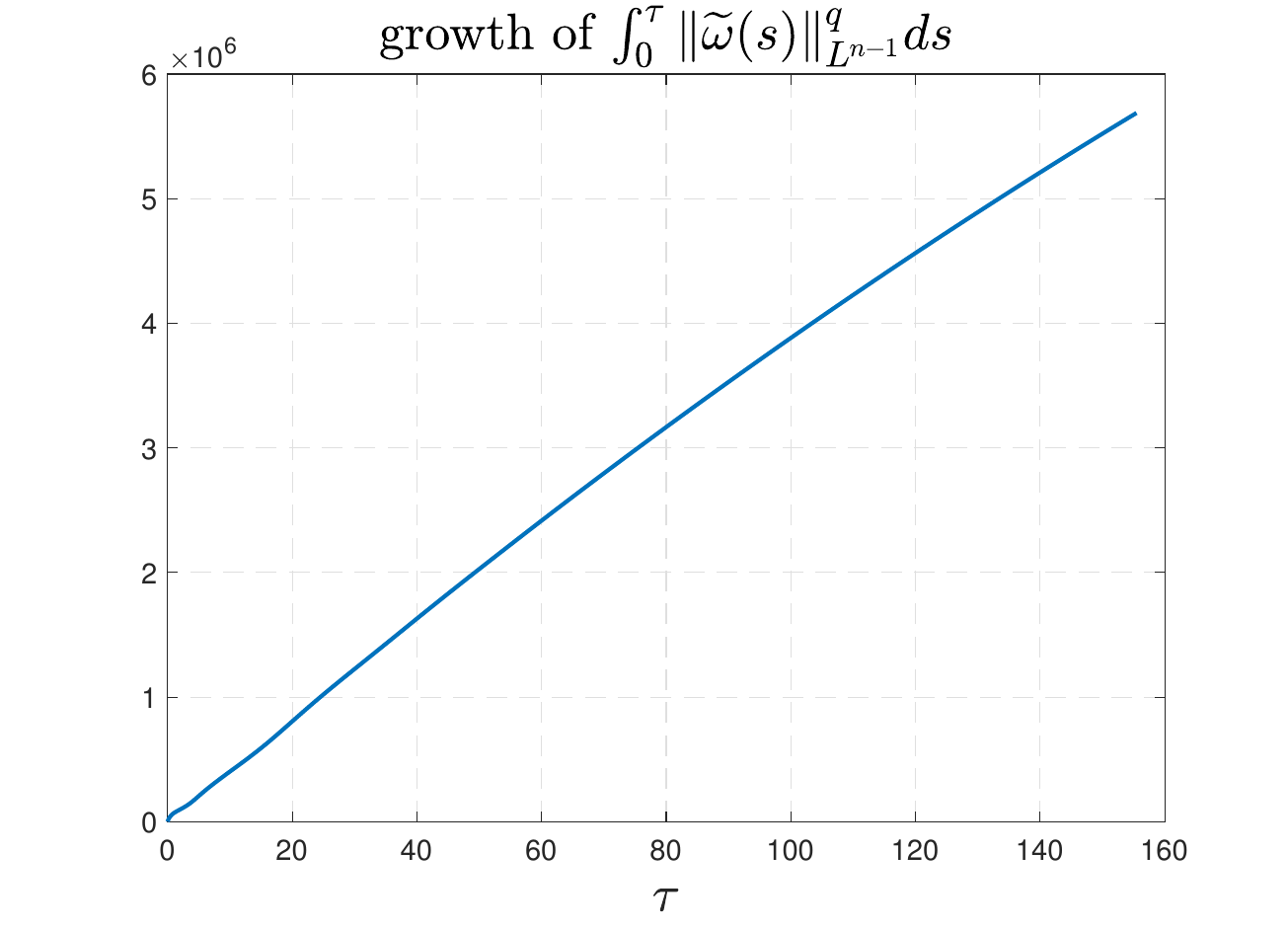} 
    \caption[Enstrophy-time]{ Left plot: The dynamic growth of the enstrophy $\int|\vom (t)|^{n-1} r^{n-2} dr dz$ for the rescaled Navier-Stokes model as a function of $\tau$. Right plot: The dynamic growth of $\int_0^\tau\|\widetilde{\vom}(s)\|_{L^{n-1}}^q ds$ with $q=\frac{2(n-1)}{n-2}$.} 
    \label{fig:enstrophy-growth}
\end{figure}

We first study the growth rate of a generalized enstrophy. 
For the $3$D axisymmetric Navier--Stokes, the enstrophy is defined as $\int |\vom (t)|^2 r dr dz$. In the $n$-dimensional setting, we define a generalized enstrophy, $\int|\vom (t)|^{n-1} r^{n-2} dr dz$. Using scaling analysis, one can show that if the rescaled Navier--Stokes model develops a self-similar blowup, $\int_0^T\|\vom (t)\|_{L^{n-1}}^q dt$ with $q=\frac{2(n-1)}{n-2}$ will blow up in finite time. Since $\|\vom (t)\|_\infty = O(1/(T-t))$ and $C_{lz}$ and $C_{lr}$ scaled like $(T-t)^{1/2}$, we expect that $\|\vom (t)\|_{L^{n-1}}^q$ roughly scales like $(T-t)^{-1}$. In Figure 
\ref{fig:enstrophy-growth} (a), we observe that  $\|\vom (t)\|_{L^{n-1}}^{n-1}$ develops rapid dynamic growth. In Figure 
\ref{fig:enstrophy-growth} (b), we plot $\int_0^T\|\vom (t)\|_{L^{n-1}}^q dt$ as a function of $\tau$. We observe that $\int_0^T\|\vom (t)\|_{L^{n-1}}^qdt$ grows slightly slower than linear growth with respect to $\tau$.

\begin{figure}[!ht]
\centering
    \includegraphics[width=0.4\textwidth]{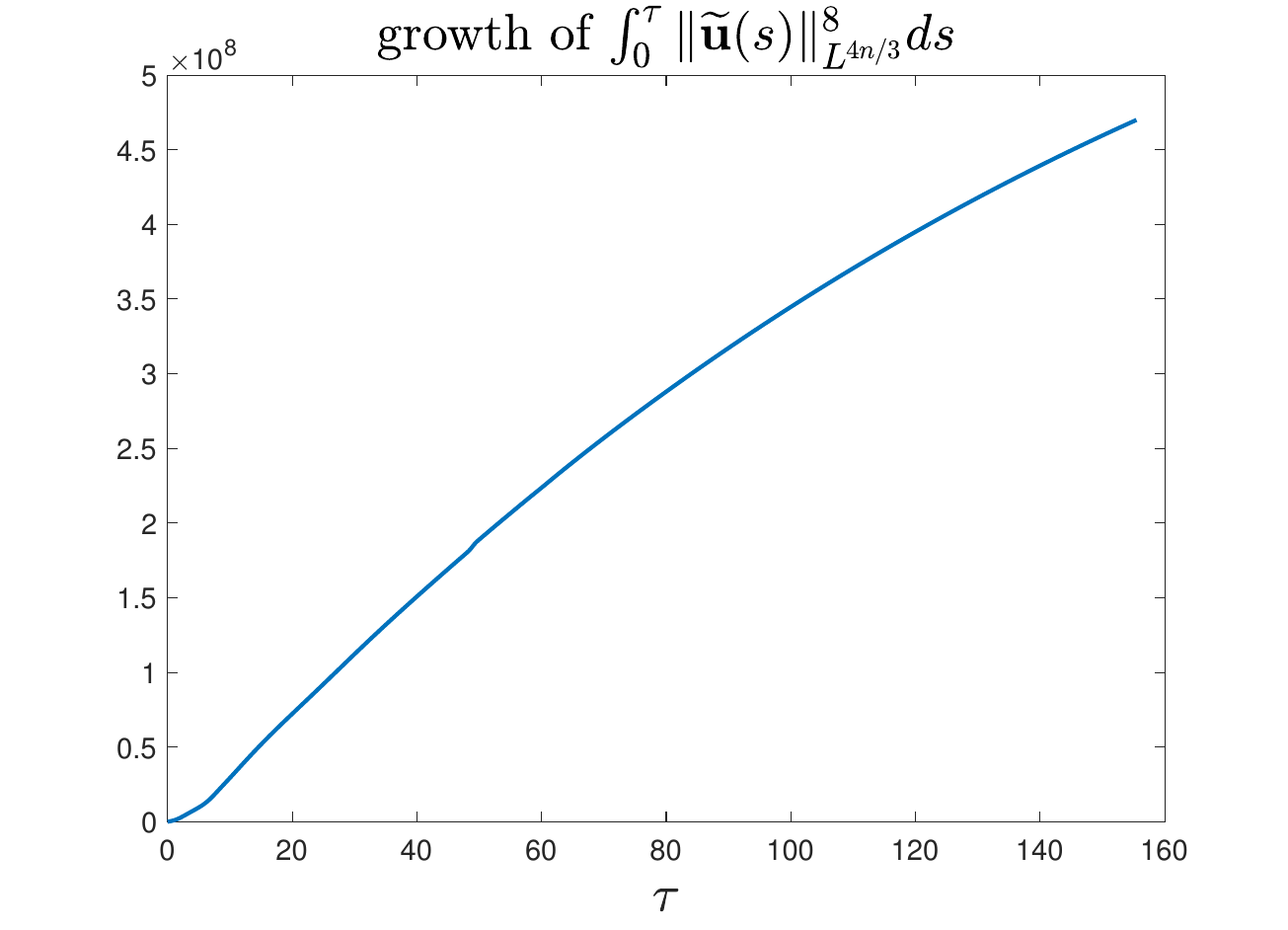}
 \includegraphics[width=0.4\textwidth]{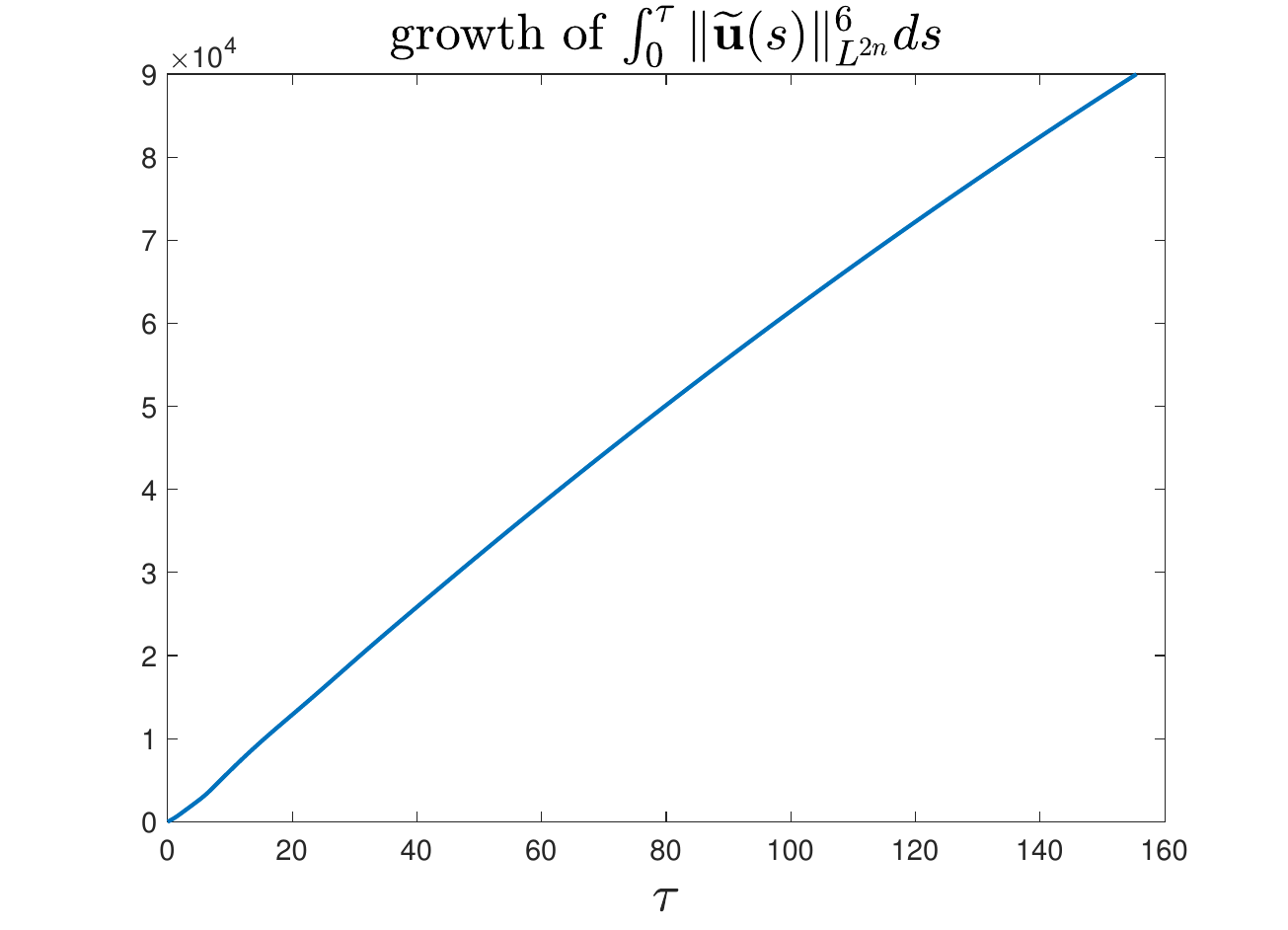} 
    \caption[Lpq-time]{ Left plot: The dynamic growth of $\int_0^\tau\|\widetilde{\bf u}(s)\|_{L^{4n/3}}^8 ds$ for the rescaled Navier-Stokes model as a function of $\tau$. 
    Right plot: The dynamic growth of $\int_0^\tau\|\widetilde{\bf u}(s)\|_{L^{2n}}^4 ds$ as a function of $\tau$. The nearly linear fitting implies that $\|{\bf u}(\tau)\|_{L^{2n,4}} \approx O(\tau^{1/4})$.} 
    \label{fig:Lpq-norm}
\end{figure}

\begin{figure}[!ht]
\centering
    \includegraphics[width=0.4\textwidth]{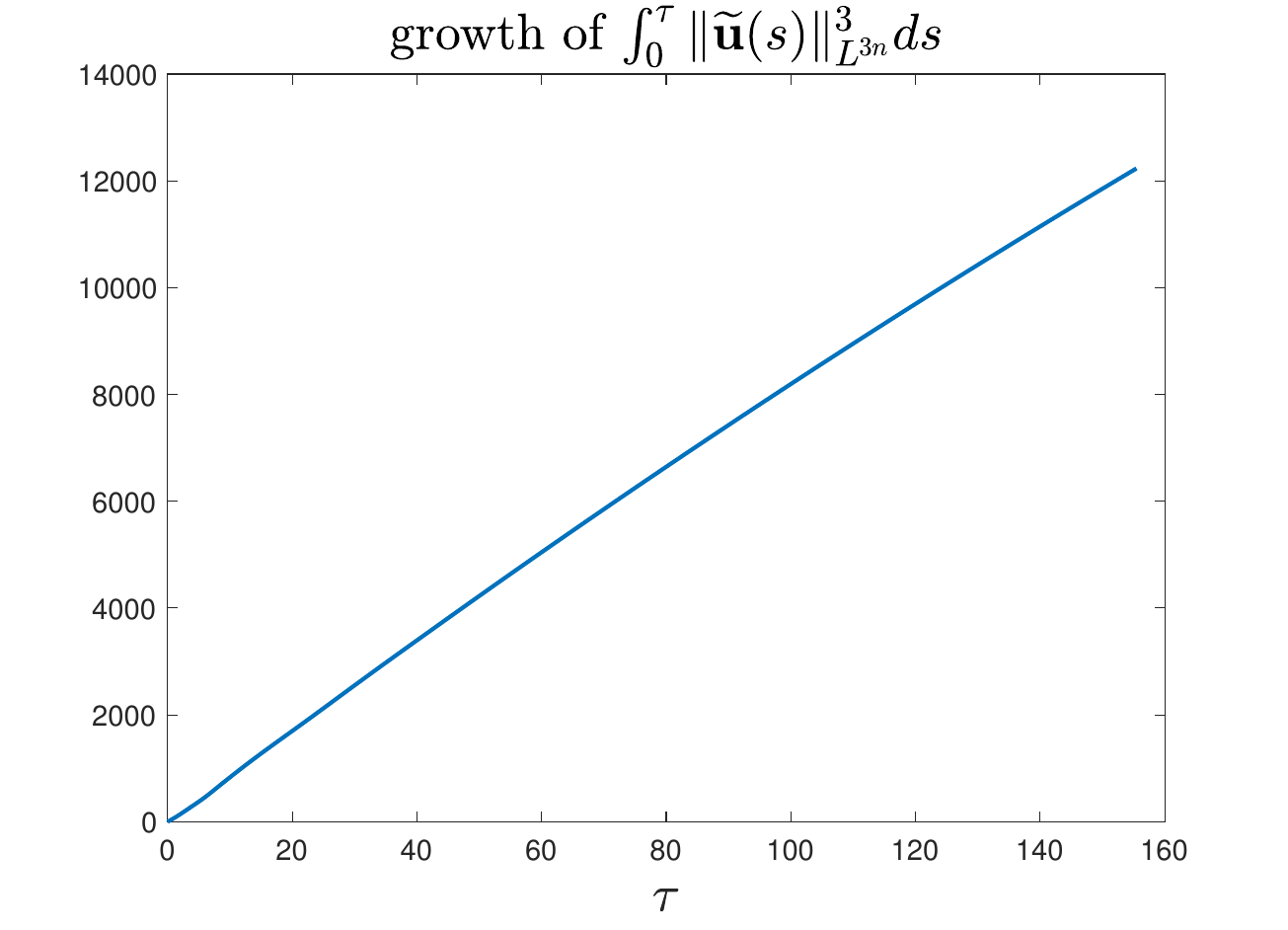}
 \includegraphics[width=0.4\textwidth]{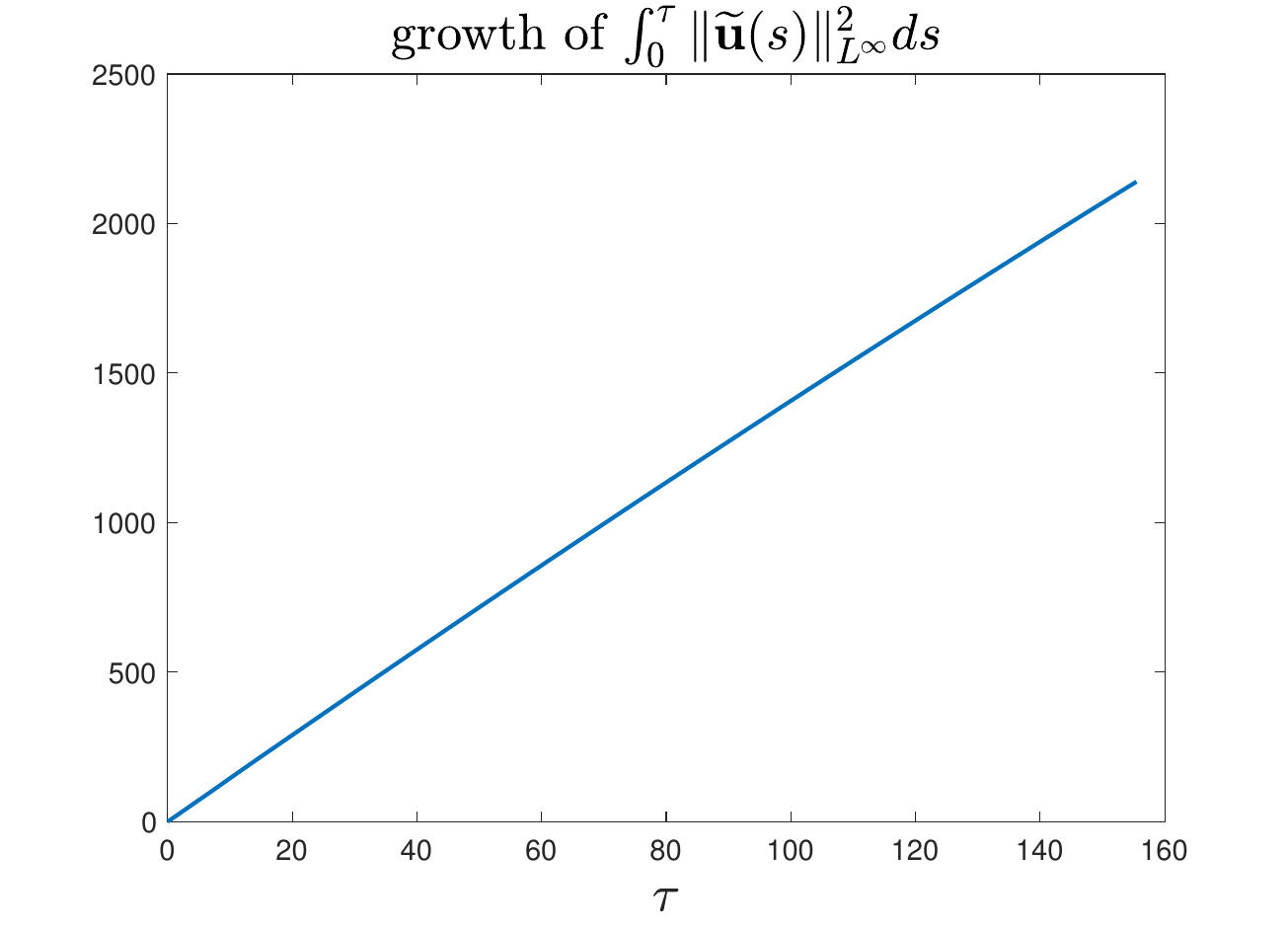} 
    \caption[Lpq-time-2]{ Left plot: The dynamic growth of $\int_0^\tau\|\widetilde{\bf u}(s)\|_{L^{3n}}^3 ds$  for the rescaled Navier-Stokes model as a function of $\tau$. 
    The almost linear fitting implies that $\|{\bf u}(\tau)\|_{L^{3n,3}} \sim O(\tau^{1/3})$. 
    Right plot: The dynamic growth of $\int_0^\tau\|\widetilde{\bf u}(s)\|_{L^\infty}^2 ds$ as a function of $\tau$. The almost linear fitting implies that $\|{\bf u}(\tau)\|_{L^{\infty}} \sim O(\tau^{1/2})$.} 
    \label{fig:Lpq-norm-2}
\end{figure}

\subsubsection{The Ladyzhenskaya-Prodi-Serrin regularity criteria} 

Next, we study the Ladyzhenskaya-Prodi-Serrin regularity criteria \cite{ladyzhenskaya1957,prodi1959,serrin1962}, which state that if a Leray-Hopf weak solution ${\bf u}$ for the $3$D Navier-Stokes equations \cite{leray1934,hopf1951} also lies in $L_t^q L_x^p$, with $3/p + 2/q \leq 1$, then the solution is unique and smooth in positive time. The endpoint result with $p=3$, $q=\infty$ has been proved in the work of Escauriaza-Seregin-Sverak in \cite{sverak2003}. In the $n$-dimensional setting, one can derive a similar result by studying the $L_t^q L_x^p$ norm of ${\bf u}$ with $n/p + 2/q \leq 1$.

In Figure \ref{fig:Lpq-norm}, we plot the dynamic growth of $\|{\bf u}\|_{L^{4n/3,8}}^8$ and $\|{\bf u}\|_{L^{2n,4}}^4$. We also plot $\|{\bf u}\|_{L^{3n,3}}^3$ in Figure \ref{fig:Lpq-norm-2}(a).
We can see that they all grow rapidly in time. For larger $p$ with $p=2n$ and $p=3n$, the growth rate is almost linear in $\tau$. This suggests that 
$\|{\bf u}\|_{L^{2n,4}} \sim O(\tau^{1/4})$ and
$\|{\bf u}\|_{L^{3n,3}} \sim O(\tau^{1/3})$.

In Figure \ref{fig:Lpq-norm-2}(b), we plot the dynamic growth of $\int_0^t \| {\bf u}(s)\|_{L^\infty}^2$.
The $L^{\infty,2}$ norm of the maximum velocity is one of the endpoint cases in the the Ladyzhenskaya-Prodi-Serrin regularity criteria with $p=\infty$ and $q = 2$. We observe that this quantity grows almost perfectly linear in $\tau$. This suggests that $\|{\bf u}(t)\|_{L^\infty}$ roughly scales like $1/(T-t)^{1/2}$, which provides further evidence for the finite time singularity of the rescaled Navier--Stokes model  with constant viscosity.

\begin{figure}[!ht]
\centering
    \includegraphics[width=0.4\textwidth]{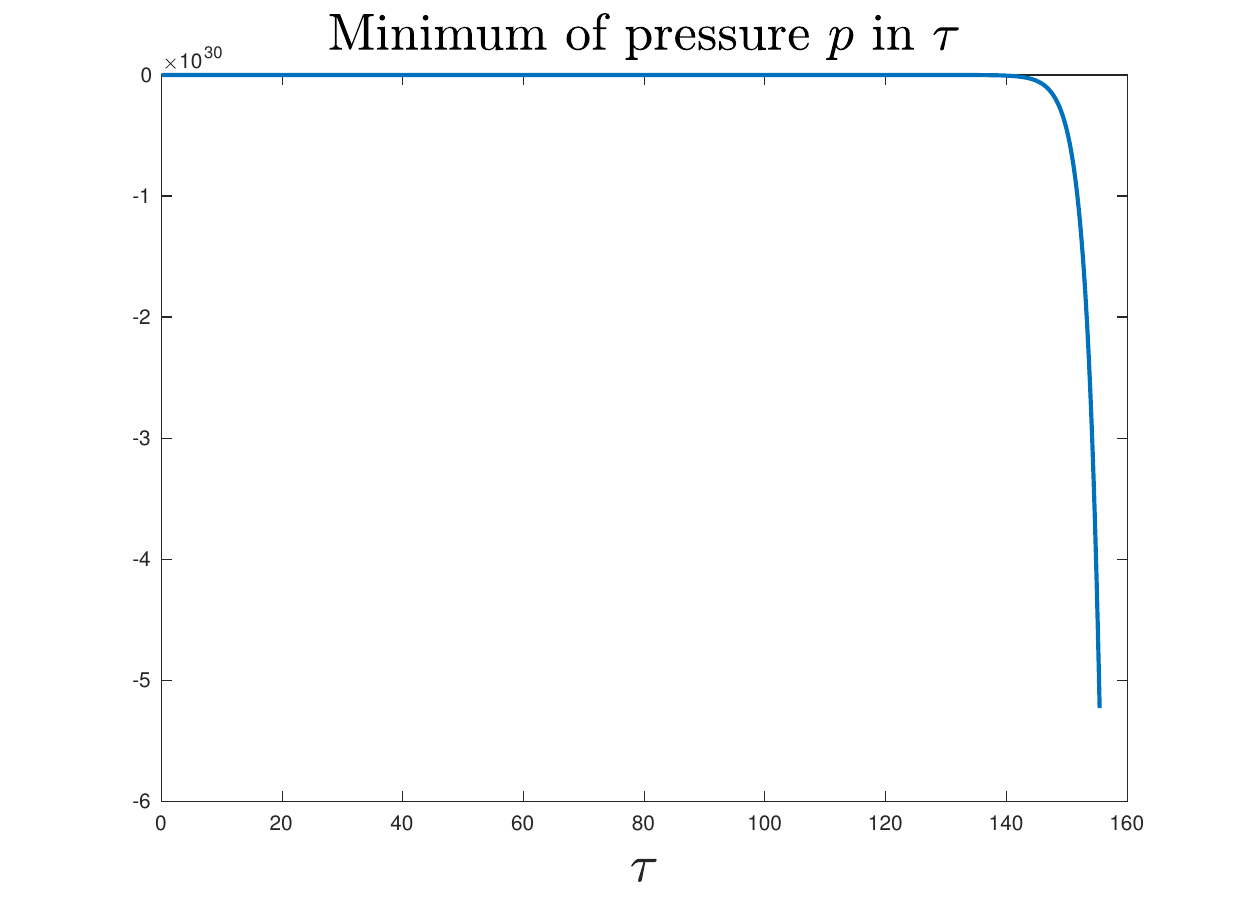}
 \includegraphics[width=0.4\textwidth]{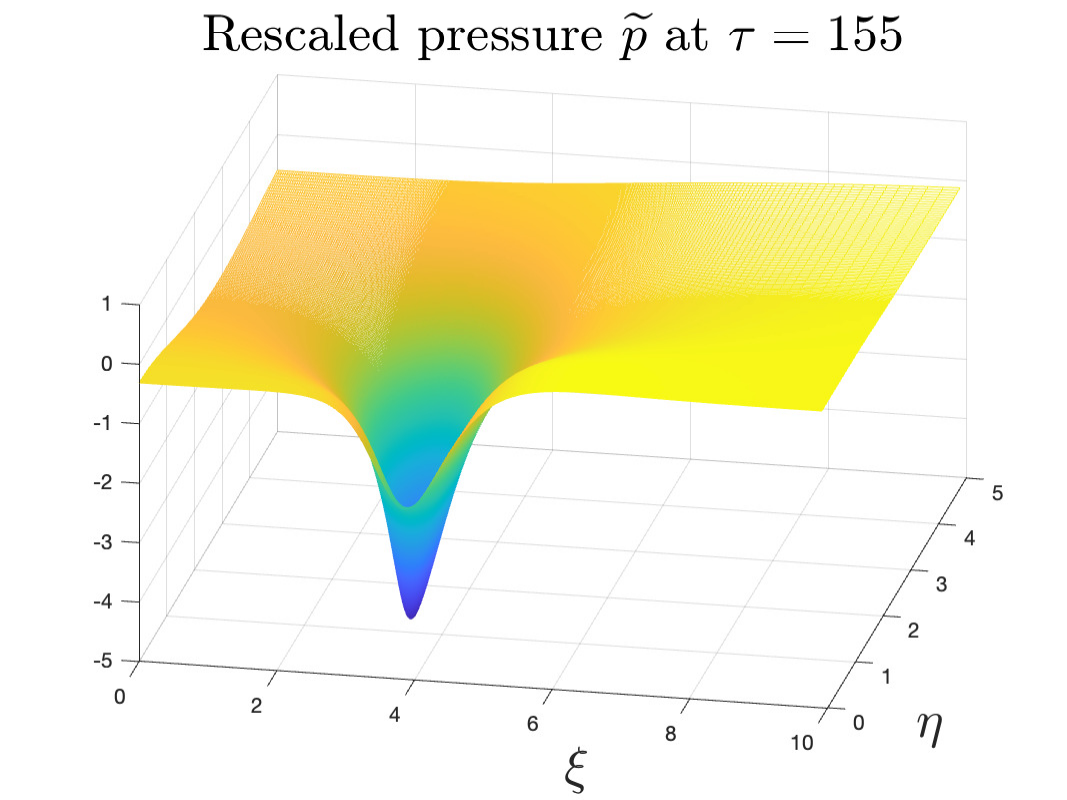} 
    \caption[Pressure-time-2]{ Left plot: The minimum of the original pressure $p$ for the rescaled Navier-Stokes model as a function of $\tau$. 
    Right plot: The profile  of the rescaled pressure $\widetilde{p}$ at $\tau=155$.
    } 
    \label{fig:pressure-profile}
\end{figure}

\begin{figure}[!ht]
\centering
    \includegraphics[width=0.4\textwidth]{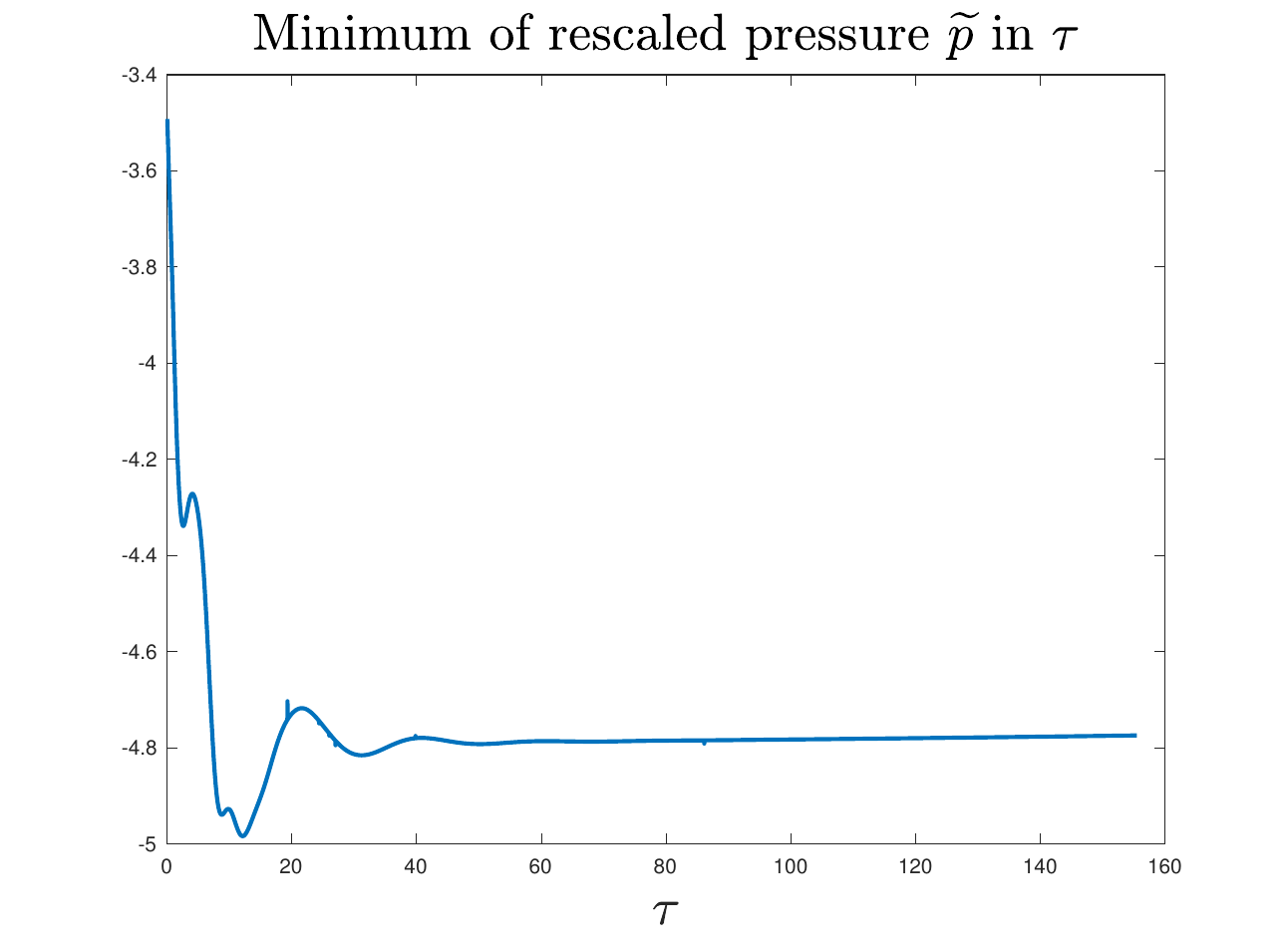}
 \includegraphics[width=0.4\textwidth]{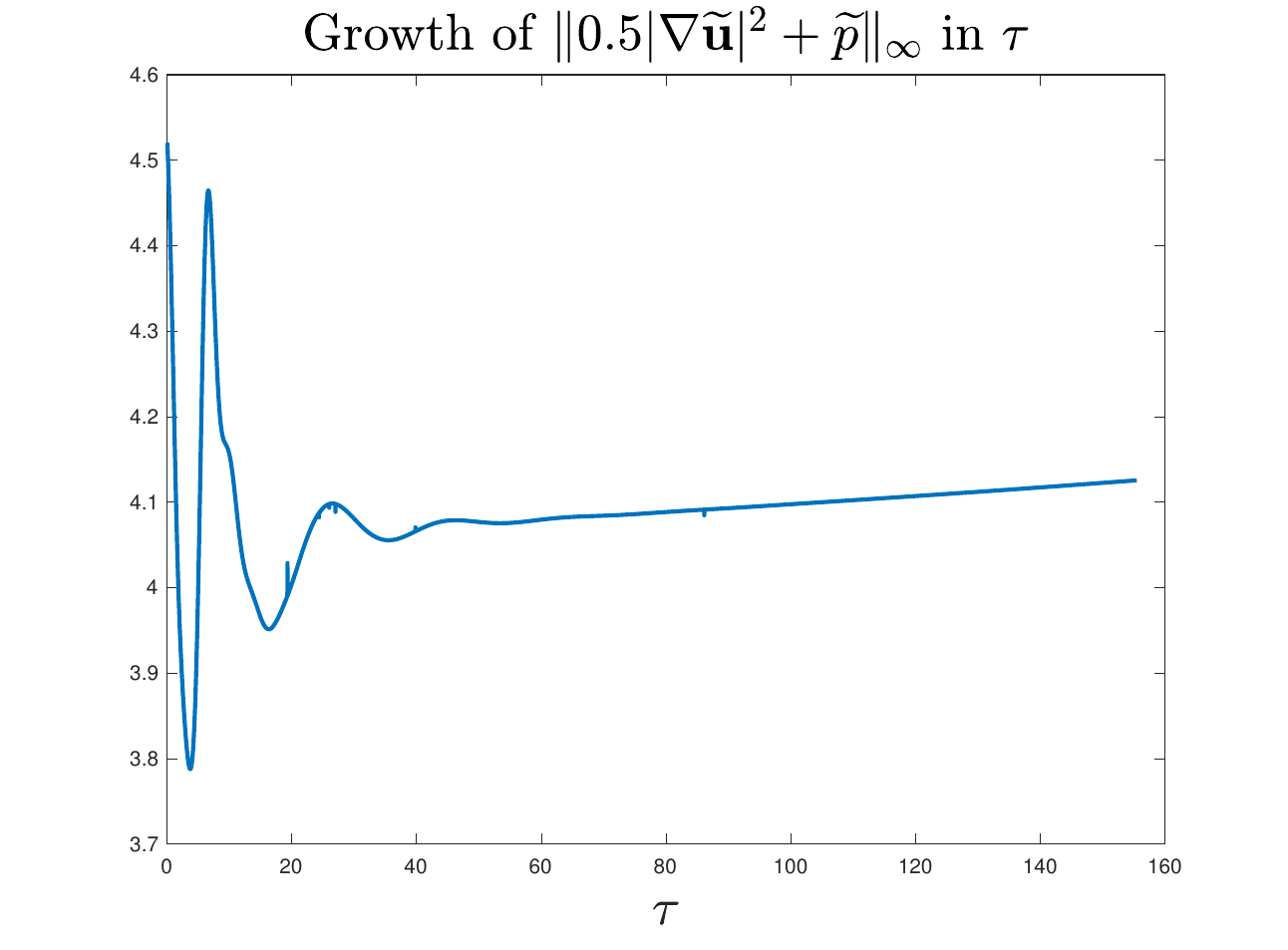} 
    \caption[Lpq-time-2]{ Left plot: The minimum  of the rescaled pressure $\widetilde{p}$ for the rescaled Navier-Stokes model as a function of $\tau$. Right plot: The growth of $\|0.5|\nabla \widetilde{\bf u}|^2+\widetilde{p}\|_\infty$ as a function of $\tau$.} 
    \label{fig:maxpre_growth}
\end{figure}

\subsubsection{The blowup of the negative pressure} 

Another blowup criteria is based on the blowup of the negative pressure \cite{sverak2002}.
In Figure \ref{fig:pressure-profile}(a), we plot the minimum of the original pressure $p$ as a function of $\tau$. We observe that the minimum of the pressure approaches to negative infinity. In Figure \ref{fig:pressure-profile}(b), we plot the rescaled pressure profile. We observe that the minimum of the rescaled pressure $\widetilde{p}$ is negative with its global minimum close to the origin. 
In Figure \ref{fig:maxpre_growth} (a), we plot the dynamic growth of the global minimum of the rescaled pressure $\widetilde{p}$ as a function of time.  In Figure \ref{fig:maxpre_growth} (b), we plot the $ \| | \frac{1}{2} |\widetilde{{\bf u}}|^2 + \widetilde{p} \|_{L^\infty} $ as a function of time. We can see that both of these two quantities stay bounded and seem to approach to a constant value as $\tau \rightarrow \infty$. Since the original pressure variable $p$ scales like $\|u_1\|_\infty \widetilde{p}$, and $\|u_1\|_\infty = 1/(T-t)$, we conclude that the minimum of the pressure goes to minus infinity with a blowup rate $O(1/(T-t))$ as $t \rightarrow T$. Similarly, we conclude that 
\[
\| p \|_\infty = O\left ( \frac{1}{T-t} \right ), \quad 
\| \frac{1}{2} | \nabla {\bf u} | + p \|_\infty = O\left ( \frac{1}{T-t} \right ).
\]
The rapid growth of these two quantities provides additional evidence for the development of potentially singular solutions of the rescaled Navier--Stokes model with constant viscosity \cite{sverak2002}.

\subsubsection{The growth of the critical $L^n$ norm of the velocity in $n$ dimensions}
\begin{figure}[!ht]
\centering
    \includegraphics[width=0.4\textwidth]{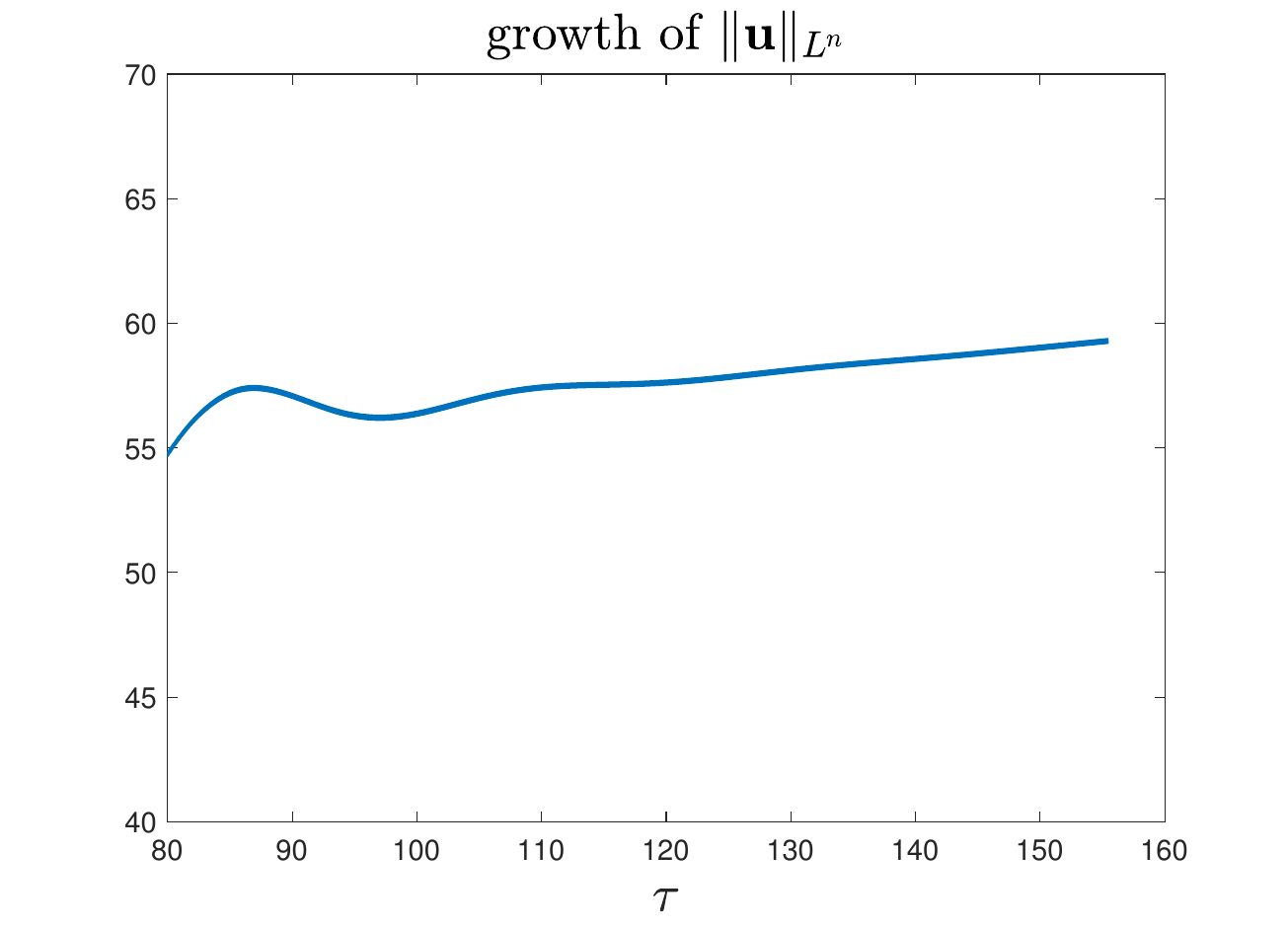}
 \includegraphics[width=0.4\textwidth]{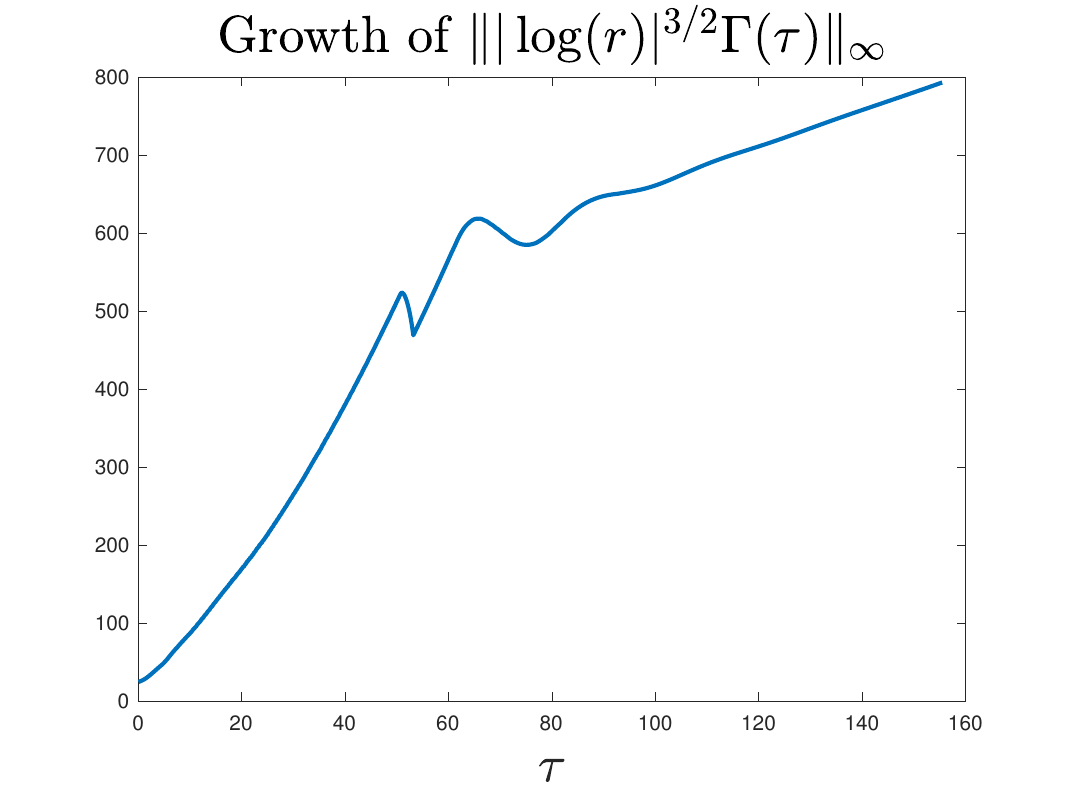} 
    \caption[VelL3-time]{ Left plot: The dynamic growth of $\|{\bf u}(\tau)\|_{L^n}$  for the rescaled Navier-Stokes model as a function of $\tau$. Right plot: The dynamic growth of $\||\log(r)|^{3/2} \Gamma(\tau, r,z)\|_{L^\infty}$ as a function of $\tau$. } 
    \label{fig:Vel-L3-Gamma}
\end{figure}

We now study the $L^{n}$ norm of the velocity field. As shown in \cite{sverak2003}, the $3$D Navier--Stokes equations cannot blow up at time $T$ if $\|{\bf u}(\tau)\|_{L^3}$ is bounded up to time $T$.  In $n$ dimensions, we should monitor the growth of $\|{\bf u} (\tau)\|_{L^n}$, which is scaling invariant.  In Figure \ref{fig:Vel-L3-Gamma} (a), we plot the dynamic growth of $\|{\bf u} (\tau)\|_{L^n}$  as a function of time in the late stage. We observe that $\|{\bf u} (t)\|_{L^n}$ experiences a mild logarithmic growth.  
Here we only plot the growth of $\|{\bf u} (\tau)\|_{L^n}$ in the late stage. 

We remark that the non-blowup criterion for $3$D Navier--Stokes using the  $\|{\bf u}\|_{L^{3}}$ estimate is based on a compactness argument. As a result, the bound on $\max_{0\leq t \leq T}\|{\bf u}(t)\|_{L^3}$ does not provide a direct estimate on the dynamic growth rate of the $3$D Navier--Stokes solution up to $T$.
In a recent paper \cite{tao2020}, Tao further examined the role of the $L^3$ norm of the velocity on the potential blow-up of the $3$D Navier-Stokes equations. He showed that as one approaches a finite blow-up time $T$, the critical $L^3$ norm of the velocity must blow up at least at a rate $\left (\log \log \log \frac{1}{T-t}\right )^c$ for some absolute constant $c$. This implies that even for a potential finite time blow-up of the Navier--Stokes equations, $\|{\bf u}(t)\|_{L^3}$ may blow up extremely slowly. Morever, the blow-up rate could be even slower for higher dimensions. We refer to \cite{palasek22} for a generalized result for dimension $n \geq 4$ by Palasek who showed that $\|{\bf u}(t)\|_{L^n}$ must blow up at least at a rate $\left (\log \log \log \log \frac{1}{T-t}\right )^c$ for some absolute constant $c$.

We also examine another non-blowup criteria based on the bound of $\| |\log(r)|^{3/2} \Gamma (t)\|_{L^\infty(r \leq r_0)}$ by D. Wei in  \cite{wei2016criticality} (see also a related paper by Lei and Zhang in \cite{lei2017criticality}). 
In Figure \ref{fig:Vel-L3-Gamma}(b), we plot the dynamic growth of $\| |\log(r)|^{3/2} \Gamma(\tau)\|_{L^\infty(\Omega(\tau)}$ over our expanding computational domain $\Omega(\tau)$. Note that since we only expand the domain by a factor of $Z(t)^{-1/5}$, the actual domain in the original physical space is actually shrinking in time. We observe that this quantity grows roughly linearly in $\tau$ in the late stage.  This implies that the non-blowup condition stated in \cite{wei2016criticality,lei2017criticality} is also violated. 

\begin{figure}[!ht]
\centering
    \includegraphics[width=0.4\textwidth]{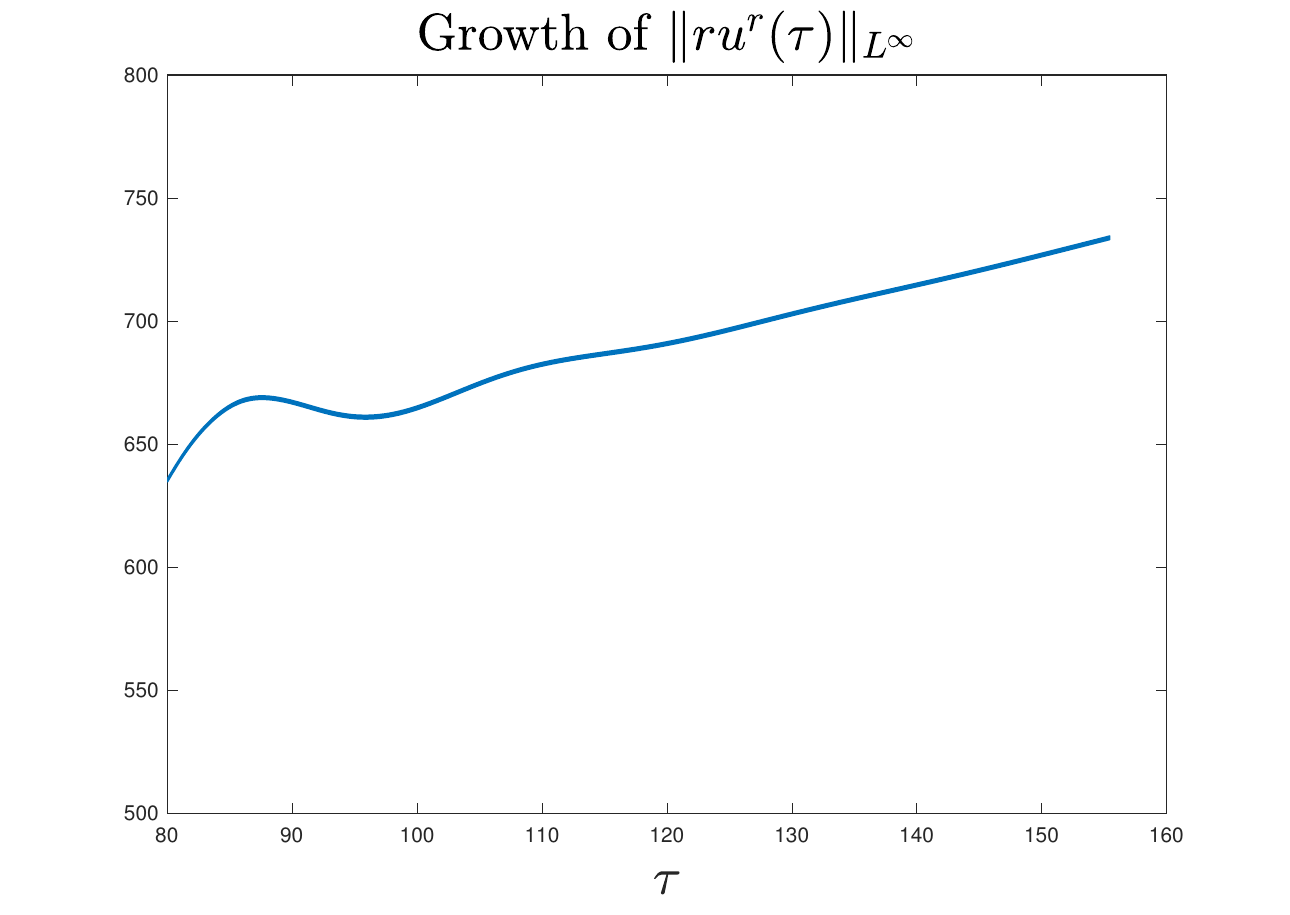}
 \includegraphics[width=0.4\textwidth]{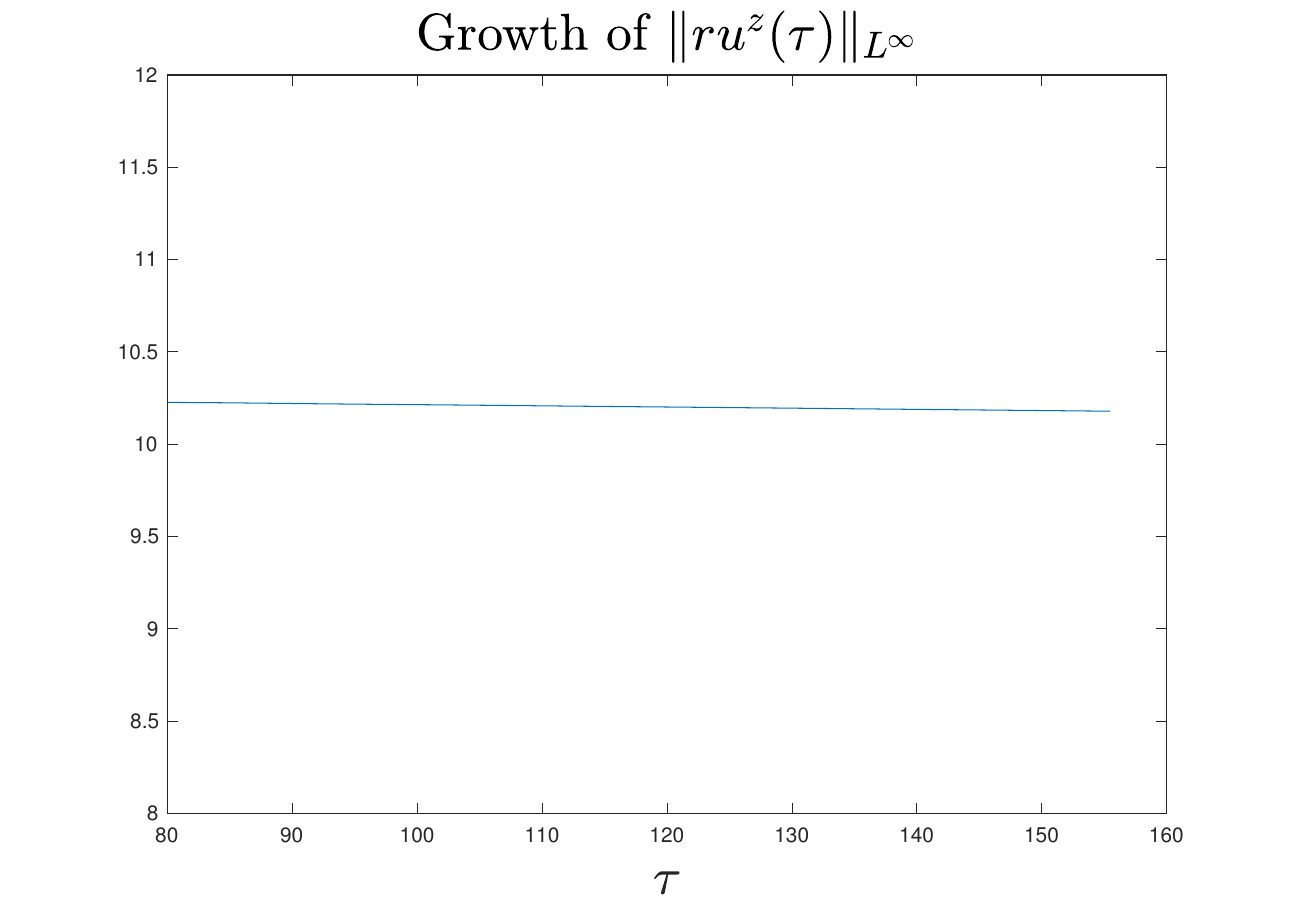} 
    \caption[VelL3-time]{ Left plot: The dynamic growth of $\|r u^r(\tau)\|_{L^\infty}$ for the rescaled Navier-Stokes model as a function of $\tau$. Right plot: The dynamic growth of $\|r u^z(\tau)\|_{L^\infty}$ as a function of $\tau$. } 
    \label{fig:maxu2r}
\end{figure}

Another important non-blowup result is the lower bound on the growth rate of the maximum velocity for the axisymmetric Navier--Stokes equations. The results in \cite{chen2008lower,chen2009lower,sverak2009} imply that the $3$D axisymmetric Navier--Stokes equations cannot develop a finite time singularity if the maximum velocity field is bounded by 
$\| {\bf u} (t)\|_{L^\infty} \leq C(T-t)^{1/2},$
provided that $ |r{\bf u} (t,r,z)| $ remains bounded for $r \geq r_0$ for some $r_0>0$.
These results are based on some compactness argument. In Figure \ref{fig:maxu2r}, we plot the growth $\| r u^r\|_{L^\infty}$ and $\| r u^z\|_{L^\infty}$ as a function of $\tau$.  We observe that $\| r u^r\|_{L^\infty}$ develops a mild linear growth in the late stage of the computation, which violates the non-blowup conditions stated in \cite{chen2008lower,chen2009lower,sverak2009}. 

\subsection{Balance between the source term and the diffusion term}

In this subsection, we study the balance between the nearly singular Delta function like source term and the diffusion term in the $\widetilde{\omega}_1$ equation. Since the density $\widetilde{\Gamma}$ satisfies a conservative advection diffusion equation, the nonlinear growth of maximum vorticity is mainly driven by the nearly singular source term in the $\widetilde{\omega}_1$ equation. It is important to monitor whether the source term and the diffusion term remain balanced throughout the computation.

\begin{figure}[!ht]
\centering
    \includegraphics[width=0.4\textwidth]{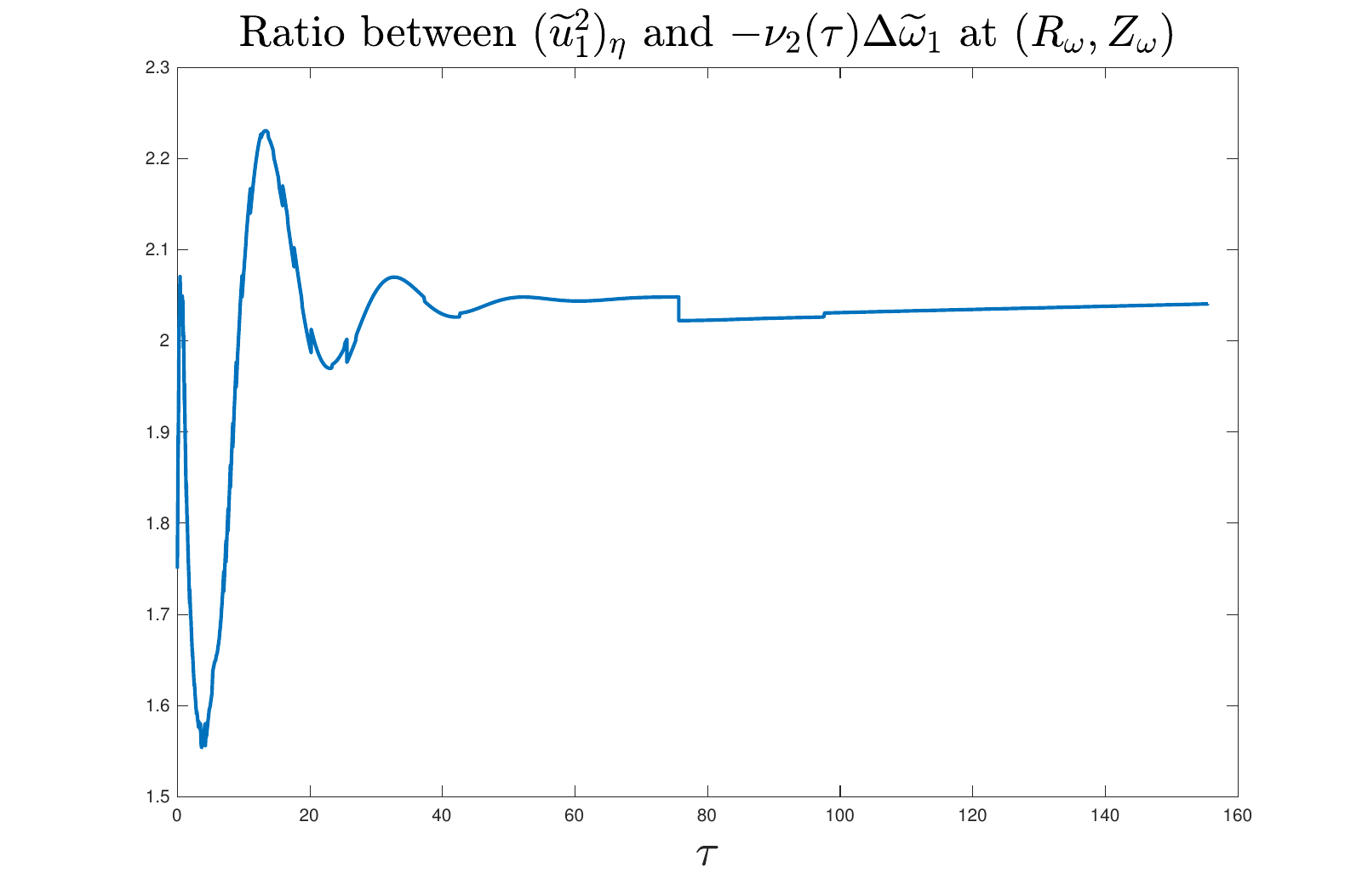}
 \includegraphics[width=0.4\textwidth]{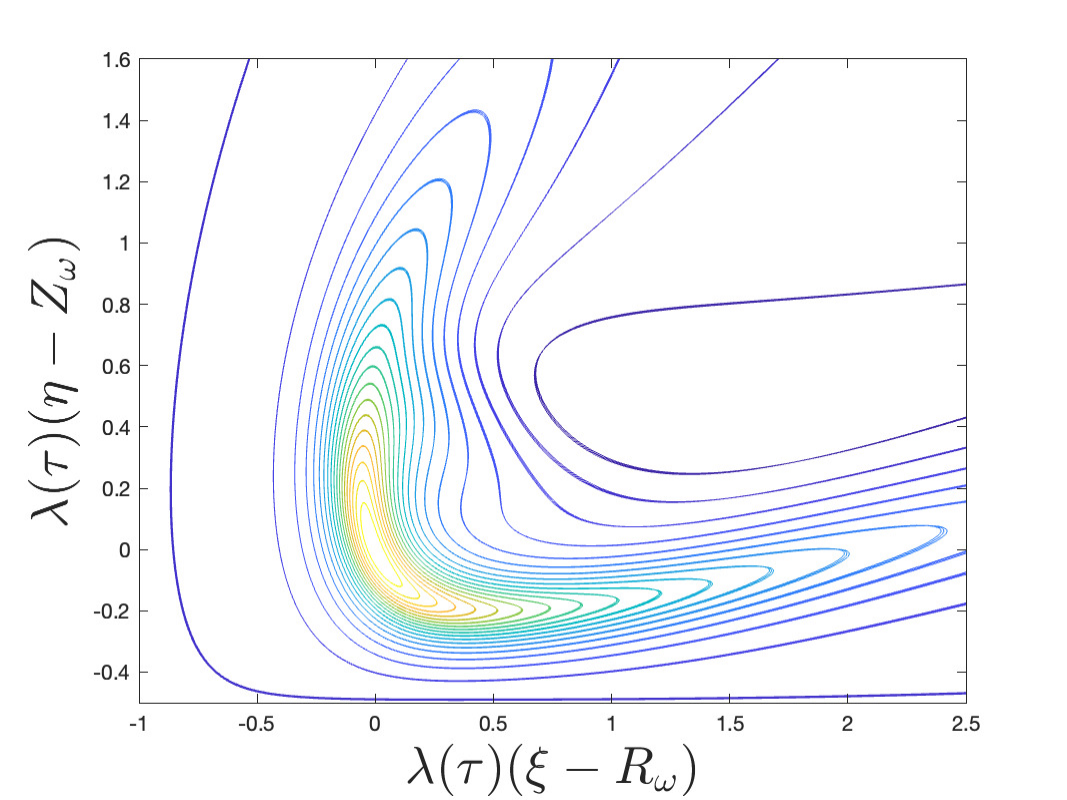} 
    \caption[vort-diff-time]{ Left plot: Ratio between the source term $(\widetilde{u}_1^2)_\eta$ and the diffusion term $\nu_2(\tau)\Delta \widetilde{\omega}_1$ at $(R_\omega,Z_\omega)$ for the rescaled Navier-Stokes model as a function of $\tau$, where $\nu_2(\tau) = \nu_2 C_\psi(\tau)/C_{lz}(\tau)$. Right plot: Contours of $\widetilde{\omega}_1$ with respect to $((\xi-R_\omega)\lambda(\tau),(\eta-Z_\omega)\lambda(\tau))$ at $\tau = 139, \;147,\; 155$ during which $\|\vom \|_\infty$ has increased by $1554$.} 
    \label{fig:ratio-vort-diff-w1}
\end{figure}

In Figure \ref{fig:ratio-vort-diff-w1}(a), we plot the ratio between the source term $(\widetilde{u}_1^2)_\eta$ 
and the diffusion term $-\nu_2(\tau)\Delta \widetilde{\omega}_1$
 at $(R_\omega,Z_\omega)$ where $\widetilde{\omega}_1$ achieves its maximum. Here $\nu_2(\tau) = \nu_2 C_\psi(\tau)/C_{lz}(\tau)$ with $\nu_2 = 0.006$. We observe that the ratio of these two terms has a mild increase in time and settles down to $2.04$ at $\tau=155$. This shows that the vortex stretching term dominates the diffusion throughout the computation. Since the vortex stretching comes from the $\omega_1$ equation only and there is no vortex stretching in the $\widetilde{\Gamma}$ equation, the balance between the vortex stretching term and the diffusion term in the $\widetilde{\omega}_1$ equation is crucial in maintaining the robust nonlinear growth of the maximum vorticity in time. 

Recall that $C_{lz}$ and $C_{lr}$ scale like $\lambda(t)\sqrt{T-t}$. In Figure \ref{fig:ratio-vort-diff-w1}(b), we plot the contours of $\widetilde{\omega}_1$ 
as a function of $(\lambda(\tau) (\xi-R_\omega ),\lambda(\tau) (\eta-Z_\omega))$ for three different time instants, $\tau=139, \; 147, \; 155$ using resolution $1024\times1024$. 
During this time interval, the maximum vorticity has increased by a factor of $1554$. We observe that these contours are almost indistinguishable from each other. This shows that ${\omega}_1$ actually enjoys a parabolic scaling property within the inner region centered at $(R_\omega, Z_\omega)$ with local scaling proportional to 
$C_{lz}/\lambda(t) \sim \sqrt{T-t}$ and $C_{lr}/\lambda(t) \sim \sqrt{T-t}$. This explains why we can achieve the balance between the source term $(\widetilde{u}_1^2)_\eta$ and the diffusion term $\nu_2(\tau)\Delta \widetilde{\omega}_1$ within this inner region centered at $(R_\omega,Z_\omega)$ with domain size shrinking to zero at a rate $\lambda(\tau)$.


\section{Concluding Remarks}
\label{sec:conclude}

We proposed the generalized Navier--Stokes equations with solution-dependent space dimension and provided strong numerical evidence that they develop a nearly self-similar blowup with smooth initial data and solution-dependent viscosity. One of our main contributions is that we provided the first rigorous derivation of the axisymmetric Navier--Stokes equations in any integer dimension  $n>3$ to the best of our knowledge.
 Moreover, we showed that the generalized Navier--Stokes equations reduce to the original Navier--Stokes equations when the dimension is an odd integer.

Due to scaling instability, traditional numerical methods enable us to approximate the potential blowup of the 3D Navier--Stokes equations but not to get arbitrarily close to the blowup time. For other nonlinear PDEs, such as the nonlinear Schr\"odinger equation or the Keller--Segel system, some unstable modes can be eliminated by utilizing the symmetry properties of the solution and studying the spectral properties of the compact linearized operator around an explicit ground state (see e.g. \cite{merle2005blow,Collot2019SpectralAF}). In our case, since we do not have an explicit ground state and the linearized operator is not compact, our approach is to enlarge the solution space by increasing the space dimension beyond 3 and using it as an additional degree of freedom to eliminate scaling instability and obtain a one-scale blowup.

A key contribution of this work is the introduction of a two-scale dynamic rescaling formulation by rescaling the $r$ and $z$ directions independently. This additional scaling freedom enables us to rescale the location of the maximum vorticity to a fixed position dynamically. Furthermore, we enforce the advection scaling uniformly along the radial and axial directions by varying the space dimension, thus preventing the development of a two-scale solution structure.

An important outcome of using a solution-dependent viscosity is that the self-similar blowup profile satisfies the self-similar equation for the generalized Navier--Stokes equations with constant viscosity. Since the generalized axisymmetric Euler equations conserve total circulation, we expect the normalized scaling exponent $\widehat{c}_{lr} \rightarrow 1/2$ as $\nu_0 \rightarrow 0$. Indeed, our numerical results suggest that the space dimension approaches 3 as the background viscosity coefficient $\nu_0$ is reduced. Thus, studying the self-similar blowup of the generalized Navier--Stokes equations offers a promising approach to investigating the potential blowup of the 3D Navier--Stokes equations.

In our future work, we plan to solve the self-similar equations directly using the background viscosity coefficient $\nu_0$ as a continuation parameter. By analyzing a sequence of self-similar profiles and studying the stability of the limiting profile as $\nu_0 \rightarrow 0$, we aim to establish a self-similar blowup of the original 3D Euler equations with the scaling exponent $c_l=1/2$. If successful, this strategy could provide valuable insights into the potential blowup of the 3D Navier--Stokes equations by treating the viscous term as a small perturbation to the 3D Euler equations and analyzing the stability of the self-similar profile.

We also investigated the nearly self-similar blowup of the rescaled Navier--Stokes model with two constant viscosity coefficients. The rescaled model preserves many properties of the 3D Navier--Stokes equations for a given constant $\tilde{n}$. Several blowup criteria were applied to analyze the nearly self-similar blowup of the rescaled model, including the Ladyzhenskaya-Prodi-Serrin non-blowup criteria, the $L^{\tilde{n}}$ norm of velocity, and specialized criteria for axisymmetric Navier--Stokes equations. These criteria support the potential finite-time singularity of the rescaled model.


It would be of great interest to develop a computer-assisted proof to justify our numerical results. Since the generalized Navier--Stokes equations with solution-dependent viscosity exhibit an asymptotically self-similar blowup, the analytical methods developed for the Hou-Luo blowup scenario \cite{ChenHou2023a,ChenHou2023b} could potentially be extended to this new context.
The analysis of the rescaled Navier--Stokes model poses additional challenges due to the presence of logarithmic corrections. 
We have recently made some encouraging progress in analyzing the stable type I blowup of the semilinear heat equation \cite{Wang2024} and the complex Ginzburg--Landau equation \cite{Wang2024b} with logarithmic corrections by using the dynamic rescaling formulation and energy estimates without relying on spectral information or topological arguments. It would be extremely interesting to extend these techniques to other nonlinear PDEs, especially the incompressible Navier--Stokes equations.

\vspace{0.1in}
{\bf Acknowledgments.} The research was in part supported by NSF Grants DMS-$2205590$ and DMS-$2508463$, the Choi Family Gift Fund and the Mike-Yan Gift Fund. I would like to thank my student Changhe Yang for a number of stimulating discussions regarding the two-scale dynamic rescaling formulation.  I would also like to thank the two anonymous referees for their very constructive comments and suggestions, especially for the suggestion to use the $e_\theta$ vector to derive the axisymmetric Navier--Stokes equations in the odd dimensional case. Finally, I would like to thank my student, Xiang Qin, for his valuable suggestion in deriving the axisymmetric Navier--Stokes equations in even dimensions and for providing a visualization of the axisymmetric flow in four and five space dimensions.

\bibliographystyle{abbrv}
\bibliography{reference}

\begin{thebibliography}{10}

\bibitem{beale1984remarks}
J.~Beale, T.~Kato, and A.~Majda.
\newblock Remarks on the breakdown of smooth solutions for the $3$-{D} {E}uler
  equations.
\newblock {\em Commun. Math. Phys.}, 94(1):61--66, 1984.

\bibitem{bp1994}
O.~N. Boratav and R.~B. Pelz.
\newblock Direct numerical simulation of transition to turbulence from a
  high-symmetry initial condition.
\newblock {\em Phys. Fluids}, 6:2757--2784, 1994.

\bibitem{brenner2016euler}
M.~Brenner, S.~Hormoz, and A.~Pumir.
\newblock Potential singularity mechanism for the {E}uler equations.
\newblock {\em Phys. Rev. Fluids}, 1:084503, 2016.

\bibitem{buckmaster2019formation}
T.~Buckmaster, S.~Shkoller, and V.~Vicol.
\newblock Formation of shocks for {2D} isentropic compressible {E}uler.
\newblock {\em Communications on Pure and Applied Mathematics},
  75(9):2069--2120, 2022.

\bibitem{buckmaster2019formation2}
T.~Buckmaster, S.~Shkoller, and V.~Vicol.
\newblock Formation of point shocks for {3D} compressible {E}uler.
\newblock {\em Communications on Pure and Applied Mathematics},
  76(9):2073--2191, 2023.

\bibitem{caffarelli1982partial}
L.~Caffarelli, R.~Kohn, and L.~Nirenberg.
\newblock Partial regularity of suitable weak solutions of the
  {N}avier--{S}tokes equations.
\newblock {\em CPAM}, 35(6):771--831, 1982.

\bibitem{caflisch93}
R.~E. Caflisch.
\newblock Singularity formation for complex solutions of the 3{D}
  incompressible {E}uler equations.
\newblock {\em Physica D: Nonlinear Phenomena}, 67(1-3):1--18, 1993.

\bibitem{caflisch89}
R.~E. Caflisch and O.~Orrelana.
\newblock Singular solutions and ill-posedness for the evolution of vortex
  sheets.
\newblock {\em SIAM J. Appl. Math.}, 20(2):249--510, 1989.

\bibitem{chen2009lower}
C.~C. Chen, R.~M. Strain, T.~P. Tsai, and H.~T. Yau.
\newblock Lower bounds on the blow-up rate of the axisymmetric
  {N}avier--{S}tokes equations {II}.
\newblock {\em Commun. PDEs}, 34(3):203--232, 2009.

\bibitem{chen2008lower}
C.~C. Chen, R.~M. Strain, H.~T. Yau, and T.~P. Tsai.
\newblock Lower bound on the blow-up rate of the axisymmetric
  {N}avier--{S}tokes equations.
\newblock {\em Intern. Math. Res. Notices}, 2008:rnn016, 2008.

\bibitem{chen2023remarks}
J.~Chen.
\newblock Remarks on the smoothness of the $ {C}^{1,\alpha}$ asymptotically
  self-similar singularity in the 3{D} {E}uler and 2{D} {B}oussinesq equations.
\newblock {\em arXiv preprint arXiv:2309.00150}, 2023.

\bibitem{chen2019finite2}
J.~Chen and T.~Y. Hou.
\newblock Finite time blowup of $2${D} {B}oussinesq and $3${D} {E}uler
  equations with ${C}^{1,\alpha}$ velocity and boundary.
\newblock {\em CMP}, 383(3):1559--1667, 2021.

\bibitem{ChenHou2023a}
J.~Chen and T.~Y. Hou.
\newblock Stable nearly self-similar blowup of the 2{D} {B}oussinesq and 3{D}
  {E}uler equations with smooth data {I}: {A}nalysis.
\newblock {\em arXiv preprint: arXiv:2210.07191v3 [math.AP]}, 2022.

\bibitem{ChenHou2023b}
J.~Chen and T.~Y. Hou.
\newblock Stable nearly self-similar blowup of the 2{D} {B}oussinesq and 3{D}
  {E}uler equations with smooth data {II}: {R}igorous numerics.
\newblock {\em MMS}, 23(1):25--130, 2025.

\bibitem{ChenHou2025}
J.~Chen and T.~Y. Hou.
\newblock Singularity formation in 3{D} {E}uler equations with smooth initial
  data and boundary.
\newblock {\em PNAS}, 122(27):e2500940122, 2025,
  https://doi.org/10.1073/pnas.2500940122.

\bibitem{chen2020finite}
J.~Chen, T.~Y. Hou, and D.~Huang.
\newblock On the finite time blowup of the {D}e {G}regorio model for the $3${D}
  {E}uler equation.
\newblock {\em CPAM}, 74(6):1281--1350, 2021,
  https://doi.org/10.1002/cpa.21991.

\bibitem{chen2021finite3}
J.~Chen, T.~Y. Hou, and D.~Huang.
\newblock Asymptotically self-similar blowup of the {H}ou-{L}uo model for the
  $3${D} {E}uler equations.
\newblock {\em Annals of PDE}, 8(24), 2022,
  https://doi.org/10.1007/s40818-022-00140-7.

\bibitem{Wang2024b}
J.~Chen, T.~Y. Hou, V.~T. Nguyen, and Y.~Wang.
\newblock On the stability of blowup solutions to the complex
  {G}inzburg-{L}andau equation in ${R}^d$.
\newblock {\em arXiv preprint: arXiv:2407.15812 [math.AP], 2024}.

\bibitem{choi2014on}
K.~Choi, T.~Y. Hou, A.~Kiselev, G.~Luo, V.~Sverak, and Y.~Yao.
\newblock On the finite-time blowup of a $1${D} model for the $3${D}
  axisymmetric {E}uler equations.
\newblock {\em CPAM}, 70(11):2218--2243, 2017.

\bibitem{choi2015}
K.~Choi, A.~Kiselev, , and Y.~Yao.
\newblock Finite time blow up for a $1${D} model of $2${D} {B}oussinesq system.
\newblock {\em CMP}, 334(3):1667--1679, 2015.

\bibitem{Collot2019SpectralAF}
C.~Collot, T.-E. Ghoul, N.~Masmoudi, and V.~T. Nguyen.
\newblock Spectral analysis for singularity formation of the two dimensional
  {K}eller--{S}egel system.
\newblock {\em Annals of PDE}, 8:1--74, 2019.

\bibitem{cfm1996}
P.~Constantin, C.~Fefferman, and A.~Majda.
\newblock Geometric constraints on potentially singular solutions for the
  $3$-{D} {E}uler equations.
\newblock {\em Commun. PDEs}, 21:559--571, 1996.

\bibitem{cordoba2023blow}
D.~C{\'o}rdoba and L.~Mart{\'\i}nez-Zoroa.
\newblock Blow-up for the incompressible 3{D}-{E}uler equations with uniform
  ${C}^{1,1/2-\epsilon} \cap l^2$ force.
\newblock {\em arXiv preprint arXiv:2309.08495}, 2023.

\bibitem{cordoba2023finite}
D.~Cordoba, L.~Martinez-Zoroa, and F.~Zheng.
\newblock Finite time singularities to the 3{D} incompressible {E}uler
  equations for solutions in ${C}^{1,\alpha} \cap {C}^{\infty}({R}^3 \backslash
  \{ 0 \} ) \cap {L}^2$.
\newblock {\em arXiv preprint arXiv:2308.12197}, 2023.

\bibitem{dhy2005}
J.~Deng, T.~Y. Hou, and X.~Yu.
\newblock Geometric properties and non-blowup of $3${D} incompressible {E}uler
  flow.
\newblock {\em Commun. PDEs}, 30:225--243, 2005.

\bibitem{es1994}
W.~E and C.-W. Shu.
\newblock Small-scale structures in {B}oussinesq convection.
\newblock {\em Phys. Fluids}, 6:49--58, 1994.

\bibitem{elgindi2019finite}
T.~M. Elgindi.
\newblock Finite-time singularity formation for ${C}^{1,\alpha}$ solutions to
  the incompressible {E}uler equations on $\mathbb{R}^3$.
\newblock {\em Annals of Mathematics}, 194(3):647--727, 2021.

\bibitem{Elg19}
T.~M. Elgindi, T.~Ghoul, and N.~Masmoudi.
\newblock On the stability of self-similar blow-up for ${C}^{1,\alpha}$
  solutions to the incompressible {E}uler equations on ${R}^3$.
\newblock {\em arXiv:1910.14071}, 2019.

\bibitem{elgindi2023instability}
T.~M. Elgindi and F.~Pasqualotto.
\newblock From instability to singularity formation in incompressible fluids.
\newblock {\em arXiv preprint arXiv:2310.19780}, 2023.

\bibitem{elgindi2023invertibility}
T.~M. Elgindi and F.~Pasqualotto.
\newblock Invertibility of a linearized boussinesq flow: a symbolic approach.
\newblock {\em arXiv preprint arXiv:2310.19781}, 2023.

\bibitem{sverak2003}
L.~Escauriaza, G.~Seregin, and V.~Sverak.
\newblock ${L}_{3,\infty}$-solutions to the {N}avier--{S}tokes equations and
  backward uniqueness.
\newblock {\em Russian Mathematical Surveys.}, 58(2):211--250, 2003.

\bibitem{fefferman2006existence}
C.~Fefferman.
\newblock Existence and smoothness of the {N}avier--{S}tokes equation.
\newblock {\em The millennium prize problems}, pages 57--67, 2006.

\bibitem{gibbon2008}
J.~Gibbon.
\newblock The three-dimensional {E}uler equations: Where do we stand?
\newblock {\em Physica D}, 237:1894--1904, 2008.

\bibitem{gs1991}
R.~Grauer and T.~C. Sideris.
\newblock Numerical computation of $3${D} incompressible ideal fluids with
  swirl.
\newblock {\em Phys. Rev. Lett.}, 67:3511--3514, 1991.

\bibitem{Tsai2023}
S.~Gustafson, E.~Miller, and T.-P. Tsai.
\newblock Growth rates for anti-parallel vortex tube {E}uler flows in three and
  higher dimensions.
\newblock {\em arXiv:2303.12043 [math.AP]}, 2023.

\bibitem{hopf1951}
E.~Hopf.
\newblock \"uber die anfangswertaufgabe f\"ur die hydrodynamischen
  grundgleichungen.
\newblock {\em Math. Nachr.}, 4:213--231, 1951.

\bibitem{Hou-euler-2022}
T.~Y. Hou.
\newblock Potential singularity of the $3${D} {E}uler equations in the interior
  domain.
\newblock {\em Found Comput Math}, 23:2203--2249, 2023.

\bibitem{Hou-nse-2022}
T.~Y. Hou.
\newblock The potentially singular behavior of the $3${D} {N}avier--{S}tokes
  equations.
\newblock {\em Found Comput Math}, 23:2251--2299, 2023.

\bibitem{Hou-Huang-2022}
T.~Y. Hou and D.~Huang.
\newblock A potential two-scale traveling wave asingularity for $3${D}
  incompressible {E}uler equations.
\newblock {\em Physica D}, 435:133257, 2022.

\bibitem{Hou-Huang-2023}
T.~Y. Hou and D.~Huang.
\newblock Potential singularity formation of $3${D} axisymmetric {E}uler
  equations with degenerate viscosity coefficients.
\newblock {\em MMS}, 21(1):218--268, 2023.

\bibitem{lei2009stabilizing}
T.~Y. Hou and Z.~Lei.
\newblock On the stabilizing effect of convection in three-dimensional
  incompressible flows.
\newblock {\em CPAM}, 62(4):501--564, 2009.

\bibitem{hou2008dynamic}
T.~Y. Hou and C.~Li.
\newblock Dynamic stability of the three-dimensional axisymmetric
  {N}avier--{S}tokes equations with swirl.
\newblock {\em CPAM}, 61(5):661--697, 2008.

\bibitem{hl2006}
T.~Y. Hou and R.~Li.
\newblock Dynamic depletion of vortex stretching and non-blowup of the 3-{D}
  incompressible {E}uler equations.
\newblock {\em J. Nonlinear Sci.}, 16:639--664, 2006.

\bibitem{hl2008}
T.~Y. Hou and R.~Li.
\newblock Blowup or no blowup? the interplay between theory and numerics.
\newblock {\em Physica D.}, 237:1937--1944, 2008.

\bibitem{Wang2024}
T.~Y. Hou, V.~T. Nguyen, and Y.~Wang.
\newblock L2-based stability of blowup with log correction for semilinear heat
  equation.
\newblock {\em arXiv preprint: arXiv:2404.09410v1 [math.AP], 2024}.

\bibitem{HouZhang2024}
T.~Y. Hou and S.~Zhang.
\newblock Potential singularity of the axisymmetric {E}uler equations with
  ${C}^\alpha$ initial vorticity for a large range of $\alpha$.
\newblock {\em MMS}, 22(4):1326--1364, 2024.

\bibitem{Huang2023b}
D.~Huang, X.~Qin, X.~Wang, and D.~Wei.
\newblock On the exact self-similar finite-time blowup of the {H}ou-{L}uo model
  with smooth profiles.
\newblock {\em arXiv preprint: arXiv:2308.01528v1 [math.AP]}, 2023.

\bibitem{Huang2023a}
D.~Huang, X.~Qin, X.~Wang, and D.~Wei.
\newblock Self-similar finite-time blowups with smooth profiles of the
  generalized {C}onstantin-{L}ax-{M}ajda model.
\newblock {\em ARMA}, 248(22), 2024,
  https://doi.org/10.1007/s00205-024-01971-3.

\bibitem{kenig2006global}
C.~E. Kenig and F.~Merle.
\newblock Global well-posedness, scattering and blow-up for the
  energy-critical, focusing, non-linear {S}chr{\"o}dinger equation in the
  radial case.
\newblock {\em Inventiones mathematicae}, 166(3):645--675, 2006.

\bibitem{kerr1993}
R.~M. Kerr.
\newblock Evidence for a singularity of the three-dimensional incompressible
  {E}uler equations.
\newblock {\em Phys. Fluids A}, 5:1725--1746, 1993.

\bibitem{kiselev2018}
A.~Kiselev.
\newblock Small scales and singularity formation in fluid dynamics.
\newblock In {\em Proceedings of the International Congress of Mathematicians},
  volume~3, 2018.

\bibitem{ladyzhenskaya1957}
A.~Kiselev and O.~Ladyzhenskaya.
\newblock On the existence and uniqueness of the solution of the nonstationary
  problem for a viscous, incompressible fluid.
\newblock {\em Izv. Akad. Nauk SSSR. Ser Mat.}, 21(5):655--690, 1957.

\bibitem{kryz2016}
A.~Kiselev, L.~Ryzhik, Y.~Yao, and A.~Zlatos.
\newblock Finite time singularity for the modified {SQG} patch equation.
\newblock {\em Ann. Math.}, 184:909--948, 2016.

\bibitem{kiselev2014small}
A.~Kiselev and V.~Sverak.
\newblock Small scale creation for solutions of the incompressible two
  dimensional {E}uler equation.
\newblock {\em Annals of Mathematics}, 180:1205--1220, 2014.

\bibitem{sverak2009}
G.~Koch, N.~Nadirashvili, G.~Seregin, and V.~Sverak.
\newblock Liouville theorems for the {N}avier--{S}tokes equations and
  applications.
\newblock {\em Acta Mathematica}, 203(1):83--105, 2009.

\bibitem{landman1988rate}
M.~J. Landman, G.~C. Papanicolaou, C.~Sulem, and P.-L. Sulem.
\newblock Rate of blowup for solutions of the nonlinear {S}chr{\"o}dinger
  equation at critical dimension.
\newblock {\em Phys. Rev. A (3)}, 38(8):3837--3843, 1988.

\bibitem{lei2017criticality}
Z.~Lei and Q.~Zhang.
\newblock Criticality of the axially symmetric {N}avier--{S}tokes equations.
\newblock {\em Pacific Journal of Mathematics}, 289(1):169--187, 2017.

\bibitem{leray1934}
J.~Leray.
\newblock Sur le mouvement d'un liquide visqueux emplissant l'espace.
\newblock {\em Acta Math.}, 63(1):193--248, 1934.

\bibitem{lin1998new}
F.~Lin.
\newblock A new proof of the {C}affarelli--{K}ohn--{N}irenberg theorem.
\newblock {\em CPAM}, 51(3):241--257, 1998.

\bibitem{liuwang2006}
J.~Liu and W.~Wang.
\newblock Convergence analysis of the energy and helicity preserving scheme for
  axisymmetric flows.
\newblock {\em SINUM}, 44(6):2456--2480, 2006.

\bibitem{luo2014potentially}
G.~Luo and T.~Y. Hou.
\newblock Potentially singular solutions of the $3${D} axisymmetric {E}uler
  equations.
\newblock {\em Proceedings of the National Academy of Sciences},
  111(36):12968--12973, 2014.

\bibitem{luo2014toward}
G.~Luo and T.~Y. Hou.
\newblock Toward the finite-time blowup of the $3${D} axisymmetric {E}uler
  equations: a numerical investigation.
\newblock {\em Multiscale Modeling \& Simulation}, 12(4):1722--1776, 2014.

\bibitem{majda2002vorticity}
A.~Majda and A.~Bertozzi.
\newblock {\em Vorticity and incompressible flow}, volume~27.
\newblock Cambridge University Press, 2002.

\bibitem{martel2014blow}
Y.~Martel, F.~Merle, and P.~Rapha{\"e}l.
\newblock Blow up for the critical generalized {K}orteweg--de {V}ries equation.
  {I}: Dynamics near the soliton.
\newblock {\em Acta Mathematica}, 212(1):59--140, 2014.

\bibitem{mclaughlin1986focusing}
D.~McLaughlin, G.~Papanicolaou, C.~Sulem, and P.~Sulem.
\newblock Focusing singularity of the cubic {S}chr{\"o}dinger equation.
\newblock {\em Physical Review A}, 34(2):1200, 1986.

\bibitem{merle2005blow}
F.~Merle and P.~Raphael.
\newblock The blow-up dynamic and upper bound on the blow-up rate for critical
  nonlinear {S}chr{\"o}dinger equation.
\newblock {\em Annals of mathematics}, 161:157--222, 2005.

\bibitem{merle1997stability}
F.~Merle and H.~Zaag.
\newblock Stability of the blow-up profile for equations of the type $u_t= {
  \Delta u } + | u|^{ p- 1} u $.
\newblock {\em Duke Math. J}, 86(1):143--195, 1997.

\bibitem{merle2015stability}
F.~Merle and H.~Zaag.
\newblock On the stability of the notion of non-characteristic point and
  blow-up profile for semilinear wave equations.
\newblock {\em CMP}, 333(3):1529--1562, 2015.

\bibitem{Tsai2022}
E.~Miller and T.-P. Tsai.
\newblock On the regularity of axisymmetric, swirl-free solutions of the
  {E}uler equation in four and higher dimensions.
\newblock {\em arXiv:2204.13406 [math.AP]}, 2022.

\bibitem{palasek22}
S.~Palasek.
\newblock A minimum critical blowup rate for the high-dimensional
  {N}avier--{S}tokes equations.
\newblock {\em Journal of Mathematical Fluid Mechanics}, 24(4):1--28, 2022.

\bibitem{frisch06}
W.~Pauls, T.~Matsumoto, U.~Frisch, and B.~J.
\newblock Nature of complex singularities for the 2{D} {E}uler equation.
\newblock {\em Physica D: Nonlinear Phenomena}, 219(1):40--59, 2006.

\bibitem{prodi1959}
G.~Prodi.
\newblock Un teorema di unicit\`a per le equazioni di {N}avier--{S}tokes.
\newblock {\em Ann. Math. Pura Appl.}, 4(48):173--182, 1959.

\bibitem{protas2022}
B.~Protas.
\newblock Systematic search for extreme and singular behaviour in some
  fundamental models of fluid mechanics.
\newblock {\em Philosophical Transactions A}, 380:20210035, 2022.

\bibitem{sverak2002}
G.~Seregin and V.~Sverak.
\newblock Navier--{S}tokes equations with lower bounds on the pressure.
\newblock {\em Arch. Rat. Mech. Anal.}, 9(1):65--86, 2002.

\bibitem{serrin1962}
J.~Serrin.
\newblock On the interior regularity of weak solutions of the
  {N}avier--{S}tokes equations.
\newblock {\em Arch. Ration. Mech. Anal.}, 9:187--191, 1962.

\bibitem{caflisch09}
M.~Siegel and R.~E. Caflisch.
\newblock Calculation of complex singular solutions to the 3{D} incompressible
  {E}uler equations.
\newblock {\em Physica D: Nonlinear Phenomena}, 238(23):2368--2379, 2009.

\bibitem{soffer90}
A.~Soffer and M.~I. Weinstein.
\newblock Multichannel nonlinear scattering for nonintegrable equations.
\newblock {\em Comm. Mth. Phys.}, 133:119--146, 1990.

\bibitem{tao2016nse}
T.~Tao.
\newblock Finite time blowup for an averaged three-dimensional
  {N}avier--{S}tokes equation.
\newblock {\em J. Amer. Math. Soc.}, 29:601--674, 2016.

\bibitem{tao2020}
T.~Tao.
\newblock Quantitative bounds for critically bounded solutions to the
  {N}avier--{S}tokes equations.
\newblock {\em arXiv:1908.04958v2 [math.AP]}, 2020.

\bibitem{wang23}
Y.~Wang, C.~Y. Lai, J.~Gomez-Serrano, and T.~Buckmaster.
\newblock Asymptotic self-similar blow-up profile for three-dimensional
  axisymmetric {E}uler equations using neural networks.
\newblock {\em Phys. Rev. Lett.}, 130:244002, 2023.

\bibitem{wei2016criticality}
D.~Wei.
\newblock Regularity criterion to the axially symmetric {N}avier--{S}tokes
  equations.
\newblock {\em J. Math. Anal. Appl.}, 435(1):402--413, 2016.

\end{thebibliography}

\end{document}